\documentclass[11pt,letterpaper,oneside,reqno]{amsart}
\usepackage[left=1in,right=1in,top=1in,bottom=1in,
           footskip=0.25in, letterpaper]{geometry}
\usepackage{graphicx,float,hyperref, tikz} 
\usepackage{bbm,amsmath,amsthm,amssymb,amsfonts,geometry,mathtools,enumerate, framed}
\usepackage[abbrev,msc-links,backrefs]{amsrefs}
\usepackage{enumitem,varwidth,tasks}

\usepackage{stmaryrd}
\usepackage{xcolor}
\usepackage{algorithm, algorithmic}
\usepackage{thmtools}

\hypersetup{colorlinks=true,citecolor=blue,linktoc=page}
\hypersetup{%
  colorlinks,
  linkcolor={red!50!black},
  citecolor={green!50!black},
  urlcolor={blue!50!black}
}
  
\usepackage{geometry}
\usepackage{setspace}
\theoremstyle{plain}
\newtheorem{theorem}{Theorem}[section]
\newtheorem{lemma}[theorem]{Lemma}
\newtheorem{claim}[theorem]{Claim}
\newtheorem{proposition}[theorem]{Proposition}
\newtheorem{observation}[theorem]{Observation}
\newtheorem{corollary}[theorem]{Corollary}

\newtheorem{question}[theorem]{Question}
\newtheorem{remark}[theorem]{Remark}
\newtheorem{definition}[theorem]{Definition}

\newcommand{\Bin}{\ensuremath{\textrm{Bin}}}
\newcommand{\Var}{\ensuremath{\textrm{Var}}}

\newcommand{\E}{\mathbb{E}}
\newcommand{\cA}{\mathcal{A}}
\newcommand{\cF}{\mathcal{F}}
\newcommand{\Prob}{\mathbb{P}}
\newcommand\cB{\mathcal B}

\newcommand\cE{\mathcal E}
\newcommand\cG{\mathcal G}

\newcommand\cK{\mathcal K}

\newcommand\cM{\mathcal M}

\newcommand\eps{\varepsilon}

\newcommand\ZZ{\mathbb{Z}}

\newcommand\PP{\mathbb{P}}

\newcommand\NN{\mathbb{N}}

\newcommand\Erw{\mathbb{E}}
\newcommand\ex{\Erw}

\newcommand\RR{\mathbb{R}}

\renewcommand{\subset}{\subseteq}
\usepackage{color}
\definecolor{red}{rgb}{1,0,0}
\definecolor{Gold}{rgb}{1,0.843,0.2}
\definecolor{DarkBlue}{rgb}{0,0,0.6}
\definecolor{magenta}{rgb}{1,0,1}

\DeclareMathOperator{\polylog}{polylog}

\let\epsilon=\varepsilon
\let\phi=\varphi
\def \1{\mathbbm{1}}

\title{Almost every $3$-uniform Hypergraph is Class 1}
\author{Arnab Chatterjee, Asaf Ferber, Marcelo Sales and Mason Shurman}
\thanks{The authors are supported by an Air force grant FA9550-23-1-0298.}
\thanks{In addition, A.F is also supported by NSF grant DMS-1953799, NSF Career DMS-2146406, and a Sloan's fellowship. }

\address{Department of Mathematics, University of California, 
    Irvine, CA, USA}
\email{\{asaff|mshurman\}@uci.edu, mtsales@caltech.edu, arnab.chatterjee@tu-dortmund.de}

\begin{document}

\begin{abstract}
The chromatic index $\chi'(H)$ of a hypergraph $H$ is the smallest number of colors needed to color its edges so that no two adjacent edges share the same color. The chromatic index is trivially bounded below by the maximum degree $\Delta(H)$, and when $\chi'(H)= \Delta(H)$, we say that $H$ is Class 1. It is well known that almost every 2-graph is Class 1, and we extend this result to 3-graphs. We show for $n^{-2/3-1/200}\leq p \leq 0.99$ that the $3$-uniform binomial random hypergraph $H\sim H_{n,p}^{(3)}$ is Class 1 with high probability. We also show that this result doesn't hold for uniformities above $3$. For a broad range of $p$, when $k \geq 4$ and $k \nmid n$, the hypergraph $H \sim H_{n,p}^{(k)}$ is with high probability not Class 1.
\end{abstract}

\maketitle

\section{Introduction}
How many colors are needed to color the edges of a graph or hypergraph so that
no two intersecting edges receive the same color?  Equivalently, how many
matchings are needed to partition its edge set?  This question is one of the
classical problems in graph theory, with connections to combinatorial
optimization, scheduling, communication networks, and design theory
(e.g.~\cites{januario2016edge,gandham2008link,baranyai1975factorization}).

Formally, a proper edge-coloring of a hypergraph \(H\) is an assignment of
colors to the edges of \(H\) such that any two intersecting edges receive
different colors. Thus every color class is a matching. The \emph{chromatic
index} \(\chi'(H)\) is the smallest integer \(\ell\) such that \(E(H)\) can be
partitioned into \(\ell\) matchings.

The most basic lower bound comes from the maximum degree. Indeed, all edges
incident to a fixed vertex must receive distinct colors, and therefore
\[
    \chi'(H)\geq \Delta(H).
\]
For ordinary graphs, this trivial lower bound is remarkably close to the truth:
Vizing's theorem~\cite{Vizing} states that every graph \(G\) satisfies
\[
    \chi'(G)\in\{\Delta(G),\Delta(G)+1\}.
\]
This gives the classical distinction between \emph{Class \(1\)} graphs, for
which \(\chi'(G)=\Delta(G)\), and \emph{Class \(2\)} graphs, for which
\(\chi'(G)=\Delta(G)+1\).

At first sight, Vizing's theorem seems to almost settle the edge-coloring
problem for graphs. However, the remaining distinction between Class \(1\) and
Class \(2\) is itself highly nontrivial: deciding whether a graph is Class \(1\)
is NP-hard~\cite{Holyer_NP}. Nevertheless, many important graph classes are
known to be Class \(1\). K\"onig's line-coloring theorem~\cite{konig1916graphen}
shows that every bipartite graph is Class \(1\), while Erd\H{o}s and
Wilson~\cite{erdos1977chromatic} proved that almost all graphs are Class \(1\).
More recently, Class \(1\) results have been obtained for sufficiently
pseudorandom regular graphs~\cite{ferber20181}, and related decomposition
results for pseudorandom graphs further illustrate how delicate exact
edge-coloring questions can be even in the graph setting~\cite{ferber2024regular}.
These results naturally motivate the broader question of when the lower bound
\(\chi'(H)\geq \Delta(H)\) is attained exactly.

The situation changes substantially for hypergraphs. For \(k\geq 3\), there is
no analogue of Vizing's theorem for \(k\)-uniform hypergraphs, which we also call \(k\)-graphs. A construction of Alon and Kim~\cite{alon1997degree} shows that,
for infinitely many \(k\), there exist \(k\)-graphs \(H\) with
\[
    \chi'(H)
    \geq
    \left(k-2+\frac{1}{k-1}\right)\Delta(H).
\]
They also conjectured that this bound should be asymptotically sharp for unrestricted $k$-graphs. The best current bound on the conjecture is due to Frederickson, Hao and Kelly~\cite{frederickson2026improved}. An important feature of the Alon and Kim construction is that the $k$-graph is intersecting, so all the color classes have size exactly one. We will see that under additional constraints this bound can be substantially improved.

Answering a question of Erd\H{o}s and Hanani regarding packings and covering of sets, R\"{o}dl~\cites{rodl1985packing, frankl1985near} developed a semi-random method now called \emph{nibbling}. The method, in its standard form, allows us to obtain an almost perfect matching in an almost-regular hypergraph as long as its maximum codegree is much smaller than its maximum degree. More precisely, let $\Delta_2(H)$ be the maximum codegree of $H$, i.e., the maximum $\deg_H(u,v)$ for every $u,v \in V(H)$. The authors in~\cite{frankl1985near} showed that if $H$ is a $k$-graph with $\deg(v)=(1 + o(1))D$ for all vertices $v$ and $\Delta_2(H)=o(D)$, then $H$ contains a matching of size $(1-o(1))\frac{n}{k}$. Since its introduction, the nibble method has been developed and refined in many directions, becoming a fundamental technique for constructing large matchings, independent sets and approximate decompositions in hypergraphs (see~\cites{alon1997nearly, kahn1996asymptotically, kim1995ramsey}).

The nibble method is well suited to study the chromatic index of hypergraphs under the same codegree restriction introduced above. The obvious idea is to simply apply nibbling repeatedly to obtain a large set of nearly-perfect matchings. An even better strategy is to build an auxiliary hypergraph where near-perfect matchings correspond to good colorings in the original graph. In this direction, Pippenger and Spencer showed that an almost $D$-regular $k$-graph $H$ with codegree $o(D)$ satisfies $\chi'(H)=(1+o(1))D$. Their original proof does not provide any explicit estimate on the $o(1)$ error term. The first explicit quantitative improvements were obtained by Molloy and Reed~\cite{molloy2000near}, who gave an explicit bound in terms of the maximum codegree. More precisely, they showed that if $H$ is a $k$-graph with $\Delta(H)\leq D$ and $\Delta_2(H)\leq C$ with $C=o(D)$, then
\begin{align*}
    \chi'(H)\leq D+D\left(\frac{D}{C}\right)^{-1/k}\log^4\left(\frac{D}{C}\right),
\end{align*}
More recently, this bound was improved even further in the work of Kang, K\"{u}hn, Methuku and Osthus~\cite{kang2023new} and Gould and Kelly~\cite{gouldkelly2025}. These quantitative improvements still leave an error term of polynomial order in $D$. This might indicate a possible natural barrier on the pure nibble method, as suggested in~\cite{gouldkelly2025}.

In a certain sense, the condition $\Delta_2(H)=o(D)$ can be viewed as a pseudorandom condition on the hypergraph $H$. Indeed, let $H$ be an instance of the binomial random hypergraph $H_{n,p}^{(k)}$, i.e., the hypergraph obtained by selecting each edge of the complete $k$-graph on $n$ vertices independently at random with probability $p$. Then it is well known that for an appropriate range of $p$, the maximum codegree $\Delta_2(H)=O(pn^{k-2})$ is much smaller than $\Delta(H)=\Theta(pn^{k-1})$. 

Because nibbling, the most fruitful technique in this area, only requires pseudorandomness, comparatively few papers have studied the chromatic index of random hypergraphs specifically. Kurauskas and Rybarczyk~\cite{kurauskas2015chromatic} proved that, even if the uniformity $k$ increases slowly with respect to $n$, a random uniform hypergraph admits an edge coloring with $(1+o(1))\Delta(H)$ colors with high probability. Ferber and Vu~\cite{ferber2018packing} studied a natural counterpart to the coloring problem: they used a technique called \emph{online sprinkling} to show that $H \sim H^{(k)}_{n,p}$ contains $(1-o(1))\Delta(H)$ edge-disjoint perfect matchings. Despite these works, the best previously known asymptotics for the chromatic index of the binomial random hypergraph $H_{n,p}^{(k)}$ came from the pseudorandom nibbling. The error term is different for various regimes of $p$ but is always polynomial. For example, when $p=1/2$, the best known result is due to Gould and Kelly:
\begin{theorem}[\cite{gouldkelly2025}]
For any $\eps>0$, if $H \sim  H^{(3)}_{n,1/2}$, then with high probability,
    $$\chi'\left(H \right) \leq \Delta\left( H\right)+\Delta\left( H\right)^{1+\eps}n^{-1/2}.$$
\end{theorem}


These results leave open a natural question: For random hypergraphs, can the error term be eliminated
completely? Are random hypergraphs Class 1? When $k \geq 4$, the answer is no. A simple averaging argument (see Section \ref{sec:no_vizing}) shows that when $k \geq 4$ and $k \nmid n$, $\chi'\left( H_{n,p}^{(k)}\right)> \Delta \left( H_{n,p}^{(k)} \right)$. Our main result shows that when $k = 3$ and $n^{-2/3-1/200}\leq p\leq 0.99$, the answer is yes.

\begin{theorem}\label{thm:main}
   For $p:=p(n)$ satisfying $n^{-2/3-1/200}\leq p\leq 0.99$, a random $3$-graph $H\sim H_{n,p}^{(3)}$ with high probability satisfies
   \begin{align*}
       \chi'(H)=\Delta(H).
   \end{align*}
\end{theorem}

We make two remarks about the range of $p$ in Theorem~\ref{thm:main}. First, the upper bound $p\leq 0.99$ can be relaxed to $p\leq 1-o(1)$, although our current approach still requires $p$ to be bounded away from $1$. We chose the stated bound for simplicity of presentation. Second, while the methods developed in this paper do not extend below the threshold $p=n^{-2/3-1/200}$, different techniques can be used to prove that $\chi'(H)=\Delta(H)$ throughout the range $\polylog(n)/n\leq p\leq n^{-2/3}$. This complementary regime will be addressed in forthcoming work.

Finally, we'll formally state that there is no direct analogue to the result of Erd\H{o}s and Wilson~\cite{erdos1977chromatic} and Theorem~\ref{thm:main} for uniformity $k\geq 4$. This is proved in Section \ref{sec:no_vizing}.

\begin{proposition}\label{prop:no_vizing_4}
    Let $k\geq 4$ and $n$ be integers such that $n$ is not a multiple of $k$. For $p:=p(n)$ satisfying $p\geq C\log n/n^{k-3}$ the following holds for sufficiently large $n$. The random $k$-graph $H\sim H_{n,p}^{(k)}$ with high probability satisfies $\chi'(H)>\Delta(H)$.
\end{proposition}

\section{Proof overview and probabilistic tools}\label{sec:overviewandprelim}

\subsection{Proof overview}

The following simple observation guides our proof: 
\begin{observation}\label{obs:firstobservation}
With high probability,
\(H\) has a unique vertex \(x_\star\) of maximum degree, and
\[
    \deg_H(x_\star)-\deg_H(y)
    \geq
    \frac{\sqrt{pn^2}}{\log^2 n}
\]
for every \(y\neq x_\star\) (see Theorem~\ref{thm:deggap}). Further, a proper edge-coloring uses exactly
    $\Delta(H)=\deg_H(x_\star)$
colors if and only if every color class contains an edge incident to
\(x_\star\). 
\end{observation}

With this degree gap in hand, we might try the following naive proof strategy: Reserve a reservoir \(E_\star\) of edges incident to
\(x_\star\), with \(|E_\star|\) smaller than the degree gap, and define
\[
    \widetilde H:=H\setminus E_\star .
\]
Since the reservoir is smaller than the gap, \(x_\star\) remains the unique
maximum-degree vertex of \(\widetilde H\). Now apply an existing
chromatic-index theorem for hypergraphs with small codegrees (see \cite{kang2023new,gouldkelly2025}). Such a theorem
colors \(\widetilde H\) using nearly \(\Delta(\widetilde H)\) matchings. Since
\(x_\star\) has degree \(\Delta(\widetilde H)\) in \(\widetilde H\), exactly
\(\Delta(\widetilde H)\) of these matchings must contain an edge with
\(x_\star\). The additional matchings, coming from the error term in the
chromatic-index theorem, do not contain \(x_\star\). \

To obtain an optimal edge-coloring of $H$, we must repair each of the matchings not containing $x_\star$. Ideally, we could simply insert a reserved edge from $E_\star$ into each of them. We would then put the unused edges of $E_\star$ into singleton color classes. This would yield a proper edge-coloring of $H$ where 
each color class contains $x_\star$.

Unfortunately, this strategy does not work. During the repair procedure, each reserved edge can only be used once. Thus, for our strategy to succeed, we need more reserved edges than there are deficient matchings. 
The difficulty is that the number of deficient matchings is equal to the error term in the previously-known chromatic-index
bounds, which is much larger than the degree gap. 

If our starting hypergraph was sparser, the error term from the available chromatic index bounds would be smaller. If the error term was reduced to less than the degree gap, our idea would work.
 With this in mind, we modify our approach. For simplicity, we will start by assuming that $3 \mid n$, so perfect matchings exist.

Suppose that we found a large set of perfect matchings $\cK$ that covered almost all of $H \sim H^{(3)}_{n,p}$. Let the remainder graph $R$ contain all the edges not in $\cK$. Then \(R\) is sparse, has small codegrees, and still satisfies
\[
    \deg_R(x_\star)-\deg_R(y)
    \geq
    \frac{\sqrt{pn^2}}{\log^2 n}
    \qquad\text{for every }y\neq x_\star .
\]
Our previous strategy is now viable: we reserve edges containing \(x_\star\) inside
\(R\), apply a hypergraph chromatic-index theorem to the rest of \(R\), and
then use a greedy algorithm to repair the deficient
matchings (see Proposition \ref{prop:recoloring}). Because \(R\) is sparse, the error in the chromatic-index
theorem is now much smaller than the degree gap, so the repair step works. This yields a proper edge-coloring of $R$ with $\Delta(R)$ colors. We combine this with $\cK$ to obtain a proper edge-coloring of $H$ with $\Delta(H)$ colors.

\subsubsection{Overview of the Iterative Process}
The core of this paper is an iterative procedure to find this large set of perfect matchings $\cK$. We have to use a method other than nibbling: while the matchings obtained from nibbling in a naive way are near-perfect, they are still not large enough to suitably sparsify the graph. The problem of packing a large set of perfect matchings in $\cK$ has previously been attacked using online sprinkling (see \cite{ferber2018packing}). This technique actually works when $p$ is very small and $H$ is sparse, and is the main idea in our forthcoming work for the sparse case. However, it fails when $H$ is dense. Therefore, in Section \ref{sec:stageI}, we introduce a different sparsification step, called Stage I, that removes perfect matchings by ``projecting'' the 3-graph into bipartite 2-graphs. Unfortunately, even after an application of Stage I, the graph is too dense to use our finishing strategy. We would like to apply Stage I repeatedly, but the output of Stage I is not pseudorandom enough to be used as input for the next application of Stage I. This motivates our actual approach: a 2-stage iterative absorption scheme used to apply Stage I multiple times. 

We will expose the edges of $H^{\left(3\right)}_{n,p}$ in a large but constant number of phases. Every edge in $K^{\left(3\right)}_n$ gets included in $G^1$ with probability $q_1$. Every edge not in $G^1$ is included in $G^2$ with probability $q_2$. In general, every edge not included in $G^1 \cup \ldots \cup G^{i-1}$ is included in $G^i$ with probability $q_i$. By appropriately choosing the $q_i$, we will have:
\[G^1 \cup G^2 \cup \ldots \cup G^{2t} \sim H^{\left(3\right)}_{n,p}.\]
Our plan, very broadly, is to do the following: Find a large set of perfect matchings $\cK_1$ in $G^{1}$ that cover all the edges in $G^1$. Define \[R^1 = G^1 \setminus \cK_1.\]
Then find a large set of perfect matchings $\cK_2$ in $G^2 \cup R^1$ that cover $R^1$. 
Set $R^2 = G^2 \cup R^1 \setminus \cK_2$. In step $i$, we find a set of perfect matchings $\cK_i$ in $G^i \cup R^{i-1}$ that cover $R^{i-1}$. We then set
\[R^i= G^i \cup R^{i-1} \setminus \cK_i.\] 
Because we stipulate that $\cK_i$ must cover $R^{i-1}$, we have that $R^i \subseteq G^i$. Therefore, if we pick the exposure probabilities $q_i$ such that the final probability $q_{2t} \approx 1/n$, we can decompose $H \sim H^{\left(3\right)}_{n,p}$ into a set of perfect matchings
\[\cK:= \cK_1 \cup \ldots \cup \cK_s\]
and a sparse remainder graph
\[R:= R^s,\]
with edge density slightly larger than $1/n$. Combining this with the finishing argument described above proves Theorem \ref{thm:main}.

The main difficulty here is finding these sets of perfect matchings $\cK_i \subseteq G^i \cup R^{i-1}$ that cover $R^{i-1}$. If $G^{i+1}$ is sparser than $G^i$, these matchings also must cover most of $G^i$, since otherwise $R^i$ will be too large and impossible to cover in the next step.

We will use two different methods to find these $\cK_i$, which we refer to as Stage I and Stage II. For either of these stages to work, $R^{i-1}$ must be very slightly sparser than $G^i$.

We can only use Stage I when the remainder graph $R^{i-1}$ is ``nice'', meaning that the codegrees are all approximately the same. In Stage I, we find a set of matchings $\cK_i$ that covers all of $R^{i-1}$ and almost all of $G^i$. The remainder graph $R^i$ is much sparser than both $R^{i-1}$ and $G^i$. This means that $G^{i+1}$ can be sparser than $G^i$, and we can set 
\[q_{i+1}=  \nu q_i\] for some common ratio $\nu \approx n^{-1/1000}$. Recall our goal is to make $G^s$ have density about $1/n$, so every time we apply Stage I, we get slightly closer. 
Unfortunately, we can not simply apply Stage I repeatedly, because the output remainder $R^i$ does not have well-concentrated codegrees. 

Stage II fixes this problem. To use Stage II, we do not require that $R^{i-1}$ has ``nice'' codegrees. In Stage II, we find a set of matchings $\cK_i$ that covers all of $R^{i-1}$ such that the new remainder $R^i$ has well-concentrated codegrees.
However, $\cK_i$ will not cover most of $G^i$. This means that $R^i$ has about the same density as $R^{i-1}$ and $G^i$, so $G^{i+1}$ has about the same density as $G^i$. Specifically, we use
\[q_{i+1}= 2q_i.\] 

Alternating between Stage I and Stage II, we slowly reduce the density of the $R^i$ and $G^i$ until $q_{2t} \approx 1/n$. We now describe how to execute Stage I and Stage II.

\subsubsection{Stage I: Reducing Density} The main idea of Stage I is to reduce the problem of finding perfect matchings in a 3-graph to finding perfect matchings in bipartite 2-graphs. Consider the following projection: Partition the vertices of $G^i \cup R^{i-1}$ into three equally sized parts $X$, $Y$, and $Z$. Pair the vertices in $X$ and $Y$ and regard each pair $xy$ as a single vertex in the projection. Every edge $xyz$ containing such a pair becomes an edge in the bipartite graph from $xy$ to $z$. Every edge that doesn't contain a pair is discarded. A perfect matching in this projected bipartite graph corresponds to a perfect matching in the original 3-graph. 
A \emph{bipartite structure} generalizes this projection idea. A bipartite structure $\Xi=\left(X,Y,Z,M,E\right)$ consists of an equipartition of the vertices $X \cup Y \cup Z$, a set $M \subseteq X \times Y$ of pairs between $X$ and $Y$, and a set $E \subseteq M \times Z$. The bipartite graph corresponding to $\Xi$ has vertex set $M \cup Z$ and only includes edges that are in $E$.

We begin Stage I by using a standard application of the probabilistic method to find a set of bipartite structures that partition the edges of $G^i \cup R^{i-1}$ (see Lemma~\ref{lem:decomposition_bipartite_structures}). Within each of the bipartite projections, the image of $G^i$ will be pseudorandom, and the image of $R^{i-1}$ will be almost-regular (see Lemma \ref{lem:bipartification}).

In the next part of Stage I, we will use a randomized algorithm adapted from Kang et al. \cite{kang2023thresholds} (see Lemma~\ref{lem:randomizedalgorithm}) to find perfect matchings in the bipartite graphs that cover the images of $R^{i-1}$, along with most of $G^i$ (see Lemma~\ref{lem:global_covering}). 
Projecting back to the 3-graph, we obtain a collection of perfect matchings $\cK_i$ in $G^i \cup R^{i-1}$ that cover all of $R^{i-1}$ and almost all of $G^i$, as desired.

\subsubsection{Stage II: Fixing the Codegrees}
Recall that the goal of Stage II is to find a set of perfect matchings $\cK_i$ in $G^i \cup R^{i-1}$ such that $G^i \cup R^{i-1} = \cK \cup R^i$ for some remainder graph $R^i$ with concentrated codegrees. 

We will use the previous chromatic index bounds \cite{kang2023new, molloy2000near} to partition $R^{i-1}$ into 
\[r:= \Delta\left(R^{i-1}\right)+\Delta\left(R^{i-1}\right)n^{-1/3}\]
disjoint matchings $E_1, E_2, \ldots, E_r$. By a simple Markov argument, nearly all of these matchings are almost-perfect (see Proposition~\ref{prop:partition_remainder_into_matchings}). We randomly split each matching into $k = \lceil \nu^{-3} \rceil \approx n^{3/1000}$ submatchings $E_{i,1}, \ldots, E_{i,k}$. For each of these submatchings $E_{i,j}$, we find a set of edges in $G^i$ that extend $E_{i,j}$ to a perfect matching. This creates our set of perfect matchings $\cK_i$ covering $R^{i-1}$.

To ensure that the new remainder graph $R^i= G^i \cup R ^{i-1} \setminus \cK_i$ has well-concentrated codegrees, we need to be careful about how we construct $\cK_i$. We first show that there are many ways to extend the $E_{i,j}$ to perfect matchings using any edges, not necessarily just those in $G^i$. In most of these collections of extensions, each pair appears in approximately the same number of edges. In other words, the codegrees are concentrated. We discard any set of extensions without this property. We then use a fractional Kahn-Kalai argument to show that with high probability, one of the remaining collections of extensions lies in $G^i$.
Since $G^i$ is distributed as a random hypergraph, and $R^i$ is obtained by removing a collection of extensions with concentrated codegrees from $G^i$, $R^i$ will have concentrated codegrees.

\subsubsection{What happens when $3 \nmid n$?}
So far in this proof overview, we have assumed that $3 \mid n$, so perfect matchings exist. If $3 \nmid n$, we execute the same strategy, finding \emph{maximum matchings} instead of perfect matchings. A maximum matching contains $\lfloor n/3 \rfloor$ edges, and all but one or two of the vertices. We must be somewhat careful which vertices we exclude. 
Recall that every matching must contain the vertex of maximum degree. 
Whenever a different vertex is excluded from a maximum matching, its degree increases relative to the maximum degree vertex. We need to pick the excluded vertices in a way that preserves the gap between the largest and second-largest degrees. 

To control this effect, after exposing $G^1$, let $S$ consist of the $n/100$ vertices of smallest degree in $G^1$. Whenever we choose a maximum matching, we exclude 1 or 2 vertices of $S$, chosen at random. Each vertex in $S$ will thus be excluded from about a $O\left(1/n\right)$ proportion of the maximum matchings. Because most of $H$ is contained in $G^1$, this will not increase the degrees enough to destroy the gap.


The paper is organized as follows. In the rest of Section \ref{sec:overviewandprelim}, we collect useful probabilistic tools and facts about random hypergraphs. We then provide a short proof of Proposition~\ref{prop:no_vizing_4} in Section~\ref{sec:no_vizing}. In Section~\ref{sec:pseudo}, we introduce the notion of pseudorandom bipartite graphs and pseudorandom 3-graphs to be used in the paper and prove several standard properties of them. In Section \ref{sec:bipartification}, we introduce the bipartification projection discussed earlier and show that pseudorandom 3-graphs project into pseudorandom bipartite graphs. Section \ref{sec:randomized} is devoted to describing the randomized algorithms used in the proofs of Stages I and II. Sections~\ref{sec:stageI} and~\ref{sec:stageII} develop Stages I and II of our iterative absorption, respectively. We present a proof of the sparsification process described earlier in Lemma~\ref{lem:mainiterative} in Section \ref{sec:iter}. Finally, in Section \ref{sec:recoloring}, we complete our strategy and provide a proof of Theorem \ref{thm:main}. 

\subsection{Probabilistic tools}

In this subsection we collect a series of probabilistic inequalities that are going to be used throughout the paper. In our proofs, we extensively use the following well-known Chernoff bounds for the upper and lower tails of binomial random variables (see \cite{JRLrandomgraphs}).
\begin{theorem}[Chernoff Bound]\label{chernoff}
Suppose $X_1, \ldots X_n$ are independent $0/1$-random variables. Let $X=\sum_{i=1}^n X_i$ and let $\mathbb{E}(X)= \mu$. Then, for $t>0$, we have
\begin{align*}
    \Prob(|X-\mu| > t)\leq 2\exp\left(-\frac{t^2}{2(\mu+t/3)}\right)
\end{align*}
In particular, if $t=\delta \mu$ for $\delta\in (0,1)$, then it holds
\begin{align*}
    \Prob(|X-\mu| > \delta \mu)\leq 2e^{-\delta^2\mu/3}
\end{align*}
\end{theorem}

We will also make use of the following version of the Chernoff bounds for hypergeometric random variables. We say that $X\sim \text{HG}(N,K,n)$ if $X$ is distributed as the number of successes in $n$ draws without replacement from a population of size $N$ containing exactly $K$ desired elements.

\begin{theorem}\label{thm:hyperchernoff}
    Let $X\sim \text{HG}(N,K,n)$ be a hypergeometric random variable and let $\E(X)=\mu=\frac{nK}{N}$. Then, for $\delta\in (0,1)$, we have
    \begin{align*}
    \Prob(|X-\mu| > \delta\mu)\leq 2e^{-\delta^2\mu/3}
\end{align*}
\end{theorem}

The next inequality is a well-known estimate of the tail bounds for martingales with bounded differences. 

\begin{theorem}[Azuma-Hoeffding]\label{thm:azuma}
Let $\{X_i\}_{i=0}^n$ be a martingale, and for each $1\leq i\leq n$, let $c_i$ be a positive constant such that $|X_i-X_{i-1}|\leq c_i$. Then for every $t>0$,
\begin{align*}
    \PP(|X_n-X_0|>t)\leq 2\exp\left(-\frac{t^2}{2\sum_{i=1}^n c_i^2}\right)
\end{align*}
\end{theorem}

A corollary of Theorem~\ref{thm:azuma} is Hoeffding's inequality, which is stated next for convenience.

\begin{theorem}[Hoeffding's Inequality]\label{thm:Hoeffding}
    Let $X_1, \ldots, X_n$ be independent random variables, such that $a_i \leq X_i \leq b_i$ almost surely for some constants $a_i$ and $b_i$. Let $S_n = X_1 + X_2 +\ldots +X_n.$ Then for any $t>0$,
    $$ \Prob\left( |S_n-\mathbb{E}(S_n)| \geq t \right) \leq 2 \exp\left(\frac{-t^2}{\sum_{i=1}^n(b_i-a_i)^2} \right)$$
\end{theorem}

 In some of our applications, Hoeffding's inequality is not strong enough, so we use the following inequality due to Freedman \cite{freedman1975tail}. We give a two-sided consequence of a formulation from \cite{tropp2011freedman}.

\begin{theorem}[Freedman] \label{thm:freedman}
Consider a real-valued martingale $\{ Y_k \}_{k\in \NN}$ with respect to a filtration $\{\cF_{k}\}_{k\in \NN}$. Let $\{ X_k \}_{k\geq 1}$ be the difference sequence, i.e., $X_k=Y_k-Y_{k-1}$. Assume that the sequence $\{X_k\}_{k\geq 1}$ is such that
$|X_k| \leq R$ almost surely for $k\geq 1$. Define the predictable quadratic variation process of the martingale $\{Y_k\}_{k\in \NN}$ as
\begin{align*}
W_k := \sum_{j=1}^k \E(X_j^2\mid \cF_{j-1}),
\end{align*}
for $k\geq 1$. Then, for all $t \geq 0$ and $\sigma^2 > 0$, it holds that
\begin{align*}
\Prob\left( \exists k \geq 0 : |Y_k| \geq t \ \text{ and }\ 
	W_k \leq \sigma^2 \right)
	\leq 2\exp \left( - \frac{t^2 }{2(\sigma^2 + Rt/3)} \right).
\end{align*}
\end{theorem}
Finally, we will use the following consequence of a concentration inequality of McDiarmid ~\cite{mcdiarmid2002concentration}.

\begin{theorem}[McDiarmid's inequality for random permutations]
\label{thm:mcdiarmid}
Let $\pi$ be a uniformly chosen random permutation, and let $Z(\pi)$ be a non-negative random variable satisfying the following two conditions:
\begin{enumerate}
    \item Swapping any two elements in $\pi$ changes the value of $Z$ by at most $c.$
    \item If $Z(\pi) \geq s$, there is a set of $rs$ coordinates in $\pi$ that certify $Z(\pi) \geq s$: If $\pi'$ agrees with $\pi$ on these coordinates, then $Z(\pi') \geq s$
\end{enumerate}swapping any two elements in $\pi$ changes the value of $Z$ by at most $c.$ Then, if $m$ is a median of $Z$, for every $t \geq 0$,
\[
    \Prob(|Z-m|\geq t)
    \leq
    4\exp\left(
        -\frac{t^2}{16rc^2(m+t)}
    \right).
\]
\end{theorem}






\section{A proof of Proposition~\ref{prop:no_vizing_4}}\label{sec:no_vizing}

In this section we provide a proof that Theorem~\ref{thm:main} does not hold for $k\geq 4$ as long as $n$ does not satisfy some divisibility requirements.

\begin{proof}[Proof of Proposition~\ref{prop:no_vizing_4}]
Note that the number of edges of $H$ follows a binomial distribution $\Bin\left(\binom{n}{k},p\right)$ with mean $\E(e(H))=p\binom{n}{k}$. Therefore, by the Chernoff bound (Theorem~\ref{chernoff}), we obtain that
\begin{align*}
    \PP\left(\big|e(H)-\E(e(H))\big|>\sqrt{pn^k \log n}\right)\leq 2\exp\left(-\frac{pn^k\log n}{3p\binom{n}{k}}\right)<\frac{1}{n}.
\end{align*}
Let $\cE$ be the event that $e(H)>p\binom{n}{k}-\sqrt{pn^k\log n}$. Then by the above inequality $\PP(\cE)\geq 1-1/n$.

Similarly, for a vertex $x\in [n]$, the degree $\deg_H(x)$ follows a binomial distribution $\Bin\big(\binom{n-1}{k-1},p\big)$. Therefore, by the Chernoff bound,
\begin{align*}
    \PP\left(\big|\deg_H(x)-\E(\deg_H(x))\big|>\sqrt{pn^{k-1} \log n}\right)\leq 2\exp\left(-\frac{pn^{k-1}\log n}{3p\binom{n-1}{k-1}}\right)\leq\frac{2}{n^2},
\end{align*}
for $k\geq 4$. Let $\cF$ be the event that $\Delta(H)\leq p\binom{n-1}{k-1}+\sqrt{pn^{k-1}\log n}$. By a union bound and the computation above, we obtain that $\PP(\cF)\geq 1-2/n$. 

We claim that if $H$ is an instance satisfying $\cE\cap \cF$, then the conclusion of the statement holds. Indeed, since $k$ does not divide $n$, every matching in $H$ contains at most $(n-1)/k$ edges. Therefore, since $H$ satisfies $\cE\cap \cF$ it holds that
\begin{align*}
    \chi'(H)\geq \frac{e(H)}{\frac{n-1}{k}}\geq \frac{k}{n}\left(\frac{p\binom{n}{k}-\sqrt{pn^k\log n}}{1-\frac{1}{n}}\right)\geq \left(1+\frac{1}{n}\right)\left(p\binom{n-1}{k-1}-\frac{k}{n}\sqrt{pn^k \log n}\right).
\end{align*}
Since $p\geq C\log n/n^{k-3}$, we obtain that
\begin{align*}
    \frac{p}{n}\binom{n-1}{k-1}-\left(\frac{k}{n}+\frac{k}{n^2}\right)\sqrt{pn^k\log n}>\sqrt{pn^{k-1}\log n}
\end{align*}
for sufficiently large $C$. Hence, it holds that
\begin{align*}
    \chi'(H)> p\binom{n-1}{k-1}+\sqrt{pn^{k-1}\log n}\geq \Delta(H),
\end{align*}
which concludes the proof of the proposition.
\end{proof}

\section{Pseudorandom properties}\label{sec:pseudo}

In this section, we describe the pseudorandom properties required from our hypergraphs and graphs during the proof. We first start with the definition of pseudorandomness for bipartite graphs used in the paper.

\begin{definition}\label{def:bipartite_pseudorandom}
    Let $\epsilon, p \in [0,1]$ and let $B$ be a bipartite graph with parts $X, Y$ with $|X|=|Y|=n$. We say that $B$ is \emph{$(\epsilon,p, n)$-pseudorandom} if the following hold:
    \begin{enumerate}[label=(A\arabic*)]
        \item\label{P1:degree} For every $z\in V(B)$, we have $\deg_B(z)=(1\pm\epsilon)pn$;
        \item\label{P2:largesets} For every $S\subseteq X$, $T\subseteq Y$ such that $|S|$, $|T|\geq \epsilon n$, it holds
        \begin{align*}
            e_B(S,T)=(1\pm \epsilon)p|S||T|;
        \end{align*}
        \item\label{P3:smallS} For every $S\subseteq X$, $T\subseteq Y$ such that $|S|<\epsilon n$ and $|T|=2|S|$, it holds
        \begin{align*}
            e_B(S,T)\leq 4\epsilon pn|S|;
        \end{align*}
        \item\label{P4:smallT} For every $S\subseteq X$, $T\subseteq Y$ such that $|T|<\epsilon n$ and $|S|=2|T|$, it holds
        \begin{align*}
            e_B(S,T)\leq 4\epsilon pn |T|.
        \end{align*}
    \end{enumerate}
    We will often refer to a bipartite graph $B$ satisfying \ref{P1:degree} as an \emph{$(1\pm \epsilon)pn$-regular} bipartite graph.
\end{definition}

\begin{remark}\label{rmk:pseudo_bipartite}
    We observe from the definition that the property of being $(\epsilon,p,n)$-pseudorandom is monotone. That is, if $\epsilon,\delta\in [0,1]$ and $\epsilon\leq \delta$, then $B$ being $(\epsilon,p,n)$-pseudorandom implies that $B$ is $(\delta,p,n)$-pseudorandom. That is the main reason why in both conditions \ref{P3:smallS} and \ref{P4:smallT} we chose the constant to be $4$.
\end{remark}

Given a bipartite graph $B$ and a weight function $\phi: B\rightarrow [0,1]$, we define the \emph{random inhomogeneous bipartite graph $B_{\phi}$} as the random subgraph of $B$ where each edge is chosen independently with probability $\PP(\{x,y\}\in B_{\phi})=\phi(x,y)$. Our first proposition shows that appropriately random inhomogeneous subgraphs of a pseudorandom bipartite graph are pseudorandom. The proof is a standard application of the Chernoff bound and we postpone it to Appendix~\ref{appendix:pseudo}.

\begin{restatable}{proposition}{bipartitepseudo}\label{prop:random_subbipartite_pseudo}
    Let $p,q,\epsilon, \delta, \eta, \gamma\in [0,1]$ be real numbers such that $pq\geq C\frac{\log n}{\delta^4 n}$ for a sufficiently large constant $C>0$ and $\epsilon, \eta\leq\min\{\delta/4,1/4\}$. Let $B$ be a bipartite graph with parts $X$ and $Y$ with $|X|=|Y|=n$ and $\phi: B\rightarrow [0,1]$ be a weight function on the edges of $B$ satisfying
    \begin{align*}
        \phi(x,y)=(1\pm \eta)q
    \end{align*}
    for every $\{x,y\}\in B$. Then the following statement holds with probability at least $1-n^{-\sqrt{C}}$: 
    \begin{enumerate}[label=(\roman*)]
         \item\label{item:random_subbipartite_deg} If $\deg_{B}(z)=(\gamma\pm\epsilon)pn$, then $\deg_{B_{\phi}}(z)=(\gamma\pm \delta)pqn$
        \item\label{item:random_subbipartite_pseudo} If $B$ is $(\epsilon,p,n)$-pseudorandom, then $B_{\phi}$ is $(\delta,pq,n)$-pseudorandom.
    \end{enumerate}
\end{restatable}

    We will also need for our proof a notion of pseudorandomness for $3$-graphs that is consistent with Definition \ref{def:bipartite_pseudorandom} in the sense that a ``bipartite projection'' of a pseudorandom $3$-graph is pseudorandom. This is obtained by the following definition. Given a $3$-graph $H$ on a set of vertices $V:=V(H)$, for a $2$-uniform matching $M\subseteq V^{(2)}$ and a set of vertices $Z\subseteq V$ we define
    \begin{align*}
        E_H(M,Z):=\left\{\{x,y,z\}\in H:\:\{x,y\}\in M,\, z\in Z\right\}
    \end{align*}
    as the set of edges in $H$ such that the first two vertices are in a pair of the matching $M$ and the last vertex is in $Z$. Correspondingly, we say that $e_H(M,Z)$ is the cardinality of the set $E_H(M,Z)$.

    \begin{definition}\label{def:hypergraph_pseudorandom}
    Let $\epsilon, p \in [0,1]$ and let $H$ be a $3$-graph on $n$ vertices. We say that $H$ is \emph{$(\epsilon,p,n)$-pseudorandom} if the following hold:
    \begin{enumerate}[label=(B\arabic*)]
        \item\label{Q1:codegree} For every $x,y\in V$, we have $\deg_H(x,y)=(1\pm\epsilon)pn$;
        \item\label{Q2:largesets} For every matching $M\subseteq V^{(2)}$ and set of vertices $Z\subseteq V$ such that $|M|$, $|Z|\geq \epsilon n$, it holds
        \begin{align*}
            e_H(M,Z)=(1\pm \epsilon)p|M||Z|;
        \end{align*}
        \item\label{Q3:smallM} For every matching $M\subseteq V^{(2)}$ and set of vertices $Z\subseteq V$ such that $|M|<\epsilon n$ and $|Z|=2|M|$, it holds
        \begin{align*}
            e_H(M,Z)\leq 4\epsilon pn|M|; 
        \end{align*}
        \item\label{Q4:smallZ} For every matching $M\subseteq V^{(2)}$ and set of vertices $Z\subseteq V$ such that $|Z|<\epsilon n$ and $|M|=2|Z|$, it holds
        \begin{align*}
            e_H(M,Z)\leq 4\epsilon pn |Z|.
        \end{align*}
    \end{enumerate}
    \end{definition}

    \begin{remark}\label{rmk:pseudo_hypergraph}
    Similarly to \ref{rmk:pseudo_bipartite}, we observe that the property of being $(\epsilon,p,n)$-pseudorandom is monotone for $3$-graphs. That is, if $\epsilon,\delta\in [0,1]$ and $\epsilon\leq \delta$, then $H$ being $(\epsilon,p,n)$-pseudorandom implies that $H$ is $(\delta,p,n)$-pseudorandom.
\end{remark}

    As in Proposition \ref{prop:random_subbipartite_pseudo}, we can show that a random subhypergraph of a pseudorandom $3$-graph is pseudorandom. Given $q\in [0,1]$ and a $3$-graph $H$, we define the \emph{random homogeneous $3$-graph $H_q$} as the $3$-graph obtained by selecting each edge $e\in H$ independently with probability $q$, i.e., for an edge $e\in H$ the probability $\PP(e\in H_q)=q$. As before, the proof also is a standard application of the concentration inequalities for binomial distributions and we also postpone it to Appendix~\ref{appendix:pseudo}.

    \begin{restatable}{proposition}{hypergraphpseudo}\label{prop:random_hypergraph_pseudo}
    Let $p,q,\epsilon, \delta\in [0,1]$ be real numbers such that $pq\geq C\frac{\log n}{\delta^4 n}$ for a sufficiently large constant $C>0$ and $\epsilon<\delta/4$. Let $H$ be an $(\epsilon,p,n)$-pseudorandom $3$-graph on $n$ vertices. Then with probability at least $1-n^{-\sqrt{C}}$ the $3$-graph $H_q$ is $(\delta,pq,n)$-pseudorandom.
    \end{restatable}

    We remark that $H_{n,p}^{(3)}=(K_n^{(3)})_p$, i.e., the binomial random $3$-graph is just the random homogeneous subhypergraph of the complete $3$-graph. Since $K_n^{(3)}$ is $(0,1,n)$-pseudorandom, the following immediate corollary follows from Proposition \ref{prop:random_hypergraph_pseudo}.

    \begin{corollary}\label{cor:binomial_hypergraph_pseudo}
        Let $p,\epsilon \in [0,1]$ be real numbers such that $p\geq \frac{C\log n}{\epsilon^4 n}$ for sufficiently large $C>0$. Then with probability at least $1-n^{-\sqrt{C}}$ the binomial random $3$-graph $H_{n,p}^{(3)}$ is $(\epsilon,p,n)$-pseudorandom.
    \end{corollary}

    We finish this section by stating the useful observation that the deletion of a small number of vertices does not destroy pseudorandomness.

    \begin{proposition}\label{prop:deletion_pseudo}
      Let $\epsilon, p \in [0,1]$ and $\eta\in [0,1/2]$ be real numbers. Let $H$ be an $(\epsilon,p,n)$-pseudorandom $3$-graph on $n$ vertices and let $D\subseteq V$ be a subset of vertices such that $|D|\leq \eta n$. Then the induced subhypergraph $\widetilde{H}:=H[V\setminus D]$ is $(\epsilon+4\eta/p, p, n-|D|)$-pseudorandom. 
    \end{proposition}

    \begin{proof}
    Note that properties \ref{Q2:largesets}, \ref{Q3:smallM} and \ref{Q4:smallZ} are inherited by induced subhypergraphs and by Remark \ref{rmk:pseudo_hypergraph} they also hold for a larger error parameter. Thus, it suffices to check that property \ref{Q1:codegree} still holds. For a fixed pair $\{x,y\} \subseteq V\setminus D$, the deletion of $|D|$ vertices affects the pair degree by at most $|D|$. Moreover the size of $D$ implies that $\widetilde{n}:=|V\setminus D|\geq (1-\eta)n$. Hence,
    \begin{align*}
        \deg_{\widetilde{H}}(x,y)&=\deg_H(x,y)\pm |D|=(1\pm(\epsilon+\eta/p))pn\\&=(1\pm(\epsilon+\eta/p))(1\pm 2\eta)p\widetilde{n}=(1\pm(\epsilon+4\eta/p))p\widetilde{n}.
    \end{align*}
    Since $\{x,y\}$ was an arbitrary pair, we have that $\widetilde{H}$ satisfies \ref{Q1:codegree}, which concludes the proof.
    \end{proof}

    \section{Bipartification}\label{sec:bipartification}
    
    One of the main methods to obtain perfect matchings in $3$-graphs in this paper is to use certain ``projections'' to reduce the problem to the problem of finding perfect matchings in bipartite graphs. The goal of this section is to describe this reduction, which we name \emph{bipartification}.

    Given any integer $n\in 3\ZZ$, we define a \emph{bipartite structure $\Xi=(X,Y,Z,M,E)$} as an ordered quintuplet where $[n]=X\cup Y\cup Z$ is an equipartition of the vertex set $[n]$ with $|X|=|Y|=|Z|=n/3$, $M\subseteq X\times Y$ is a perfect matching from $X$ to $Y$ and
    \begin{align*}
        E\subseteq M\times Z:=\big\{\{x,y,z\}:\:\{x,y\}\in M,\, z\in Z\big\}
    \end{align*}
    is a $3$-graph. Given a bipartite structure $\Xi=(X,Y,Z,M,E)$ and a $3$-graph $H$ on the same set of vertices, let $\Xi^{(3)}(H):=H\cap E$ be the \emph{bipartite projection} of $H$ with respect to $\Xi$. We call it a bipartite projection because although $\Xi^{(3)}(H)$ is a $3$-graph, one can see it as a bipartite graph $\Xi^{(2)}(H)$ by considering the set of vertices to be $M\cup Z$ and the projection map $\pi:\Xi^{(3)}(H)\rightarrow \Xi^{(2)}(H)$ given by
    \begin{align*}
        \pi(x,y,z)=\big\{\{x,y\},\, z\big\}.
    \end{align*}
    Throughout the paper we will always denote the bipartite projection of $H$ seen as a $3$-graph by $\Xi^{(3)}(H)$ and seen as a bipartite $2$-graph by $\Xi^{(2)}(H)$.

\begin{figure}[H]
  \centering
  \definecolor{bpRed}{HTML}{cd001f}
  \definecolor{bpGray}{HTML}{9A9A9A}
  \def\bpRetainedEdges{1/2,2/1,3/3,3/2}
  \begin{tikzpicture}[
    x=1cm,y=1cm,font=\small,
    bp hyperedge/.style={line width=11.3pt,
      line cap=round,line join=round},
    bp discarded/.style={bp hyperedge,draw=black,
      postaction={draw=black!5,line width=10pt}},
    bp retained/.style={bp hyperedge,draw=black!50,
      postaction={draw=black!25,line width=10pt}},
    bp pair/.style={line width=1.8pt,line cap=round},
    bp vertex/.style={circle,fill=black,inner sep=0pt,
      minimum size=4pt},
    bp label/.style={font=\scriptsize,inner sep=1pt},
    bp map/.style={->,draw=black!70,line width=.8pt,>=stealth}
  ]
    \foreach \shift/\stage in {0/1,5.6/2}{
      \begin{scope}[xshift=\shift cm]
        \foreach \idx/\py in {1/2.1,2/0,3/-2.1}{
          \coordinate (x\idx) at (0,\py);
          \coordinate (y\idx) at (1.45,\py);
          \coordinate (z\idx) at (3.5,\py);
        }
        \foreach \px/\part in {0/X,1.45/Y,3.5/Z}
          \node at (\px,2.69) {$\part$};
          \ifnum\stage=1
          \node at (1.75,3.17) {$H$};
          \draw[bp discarded] (z1) -- (z2) -- (z3);
         
          \draw[bp discarded] (x3)
            to[out=0,in=270] (y2)
            to (y1);
          \foreach \i/\j in {1/2,2/1}
            \draw[bp discarded] (x\i)
              to[out=0,in=180] (y\j) -- (z\j);
        \else
          \node at (1.75,3.17) {$\Xi^{(3)}(H)=H\cap E$};
        \fi
\foreach \i/\j in \bpRetainedEdges
          \draw[bp retained] (x\i) -- (y\i)
            to[out=0,in=180] (z\j);

        \foreach \idx/\py in {1/2.1,2/0,3/-2.1}{
          \draw[bp pair,bpRed] (x\idx) -- (y\idx);
          \foreach \part in {x,y}{
            \node[bp vertex] at (\part\idx) {};
            \node[bp label,anchor=north,yshift=-10pt] at (\part\idx)
              {$\part_{\idx}$};
          }
          \node[bp vertex] at (z\idx) {};
            \node[bp label,anchor=north,yshift=-10pt] at (z\idx)
              {$z_{\idx}$};
        }
         
      \end{scope}
    }

    \draw[bp map] (4.00,0) -- (5.08,0);
    \draw[bp map] (9.62,0) -- (10.70,0);

    \begin{scope}[xshift=11.3cm]
      \node at (1.10,3.17) {$\Xi^{(2)}(H)$};
      \node at (0,2.69) {$M$};
      \node at (2.20,2.69) {$Z$};
      \foreach \idx/\py in {1/2.1,2/0,3/-2.1}{
        \coordinate (m\idx) at (0,\py);
        \coordinate (zz\idx) at (2.20,\py);
      }
      \foreach \i/\j in \bpRetainedEdges
        \draw[black,line width=1.1pt] (m\i) -- (zz\j);
      \foreach \idx in {1,2,3}{
        \node[bp vertex,fill=bpRed,minimum size=6pt] at (m\idx) {};
        \node[bp vertex] at (zz\idx) {};
        \node[bp label,anchor=north,yshift=-10pt] at (m\idx) {$\{x_{\idx},y_{\idx}\}$};
        \node[bp label,anchor=north,yshift=-10pt] at (zz\idx) {$z_{\idx}$};
      }
    \end{scope}
  \end{tikzpicture}
  \caption{Bipartification on a 3-graph with 9 vertices. The red segments form our matching 
    $M=\{\{x_1,y_1\},\{x_2,y_2\},\{x_3,y_3\}\}$.  
    The four white hyperedges form $E$, the set of hyperedges compatible with $M$. Restricting $H$ to $E$, we obtain $\Xi^{(3)}(H)$. In the bipartite graph $\Xi^{(2)}(H)$, the edges in $M$ are contracted to vertices, and each
    compatible hyperedge $\{x_i,y_i,z_j\}$ becomes the edge $\{\{x_i,y_i\},z_j\}$.
    Observe that the perfect matching $\{\{\{x_1,y_1\},z_2\},\{\{x_2,y_2\},z_1\},\{\{x_3,y_3\},z_3\}\}$
    lifts to a perfect matching in $\Xi^{(3)}(H)$ and $H$ itself.}
  \label{fig:bipartification-three-pairs-recolored}
\end{figure}
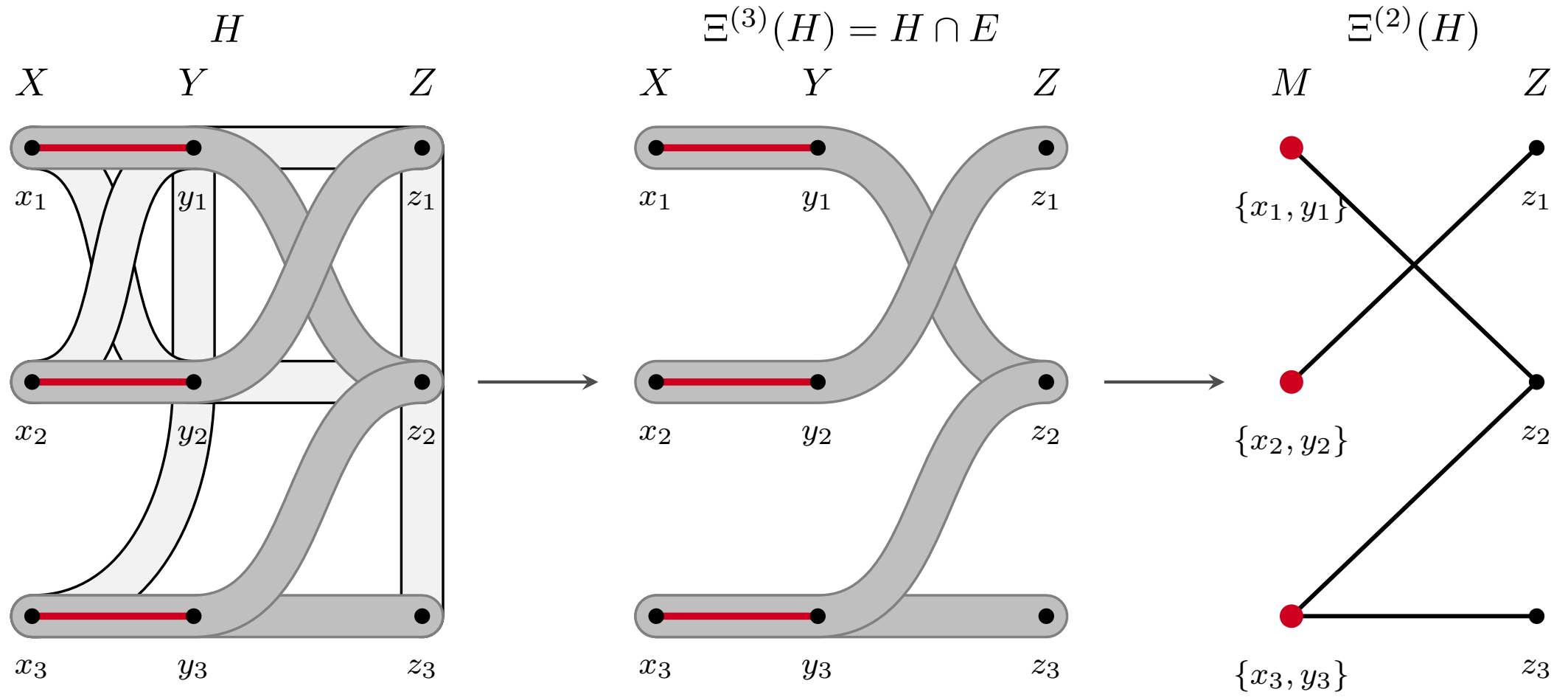

    The following definition will be used as our model of random bipartification.

    \begin{definition}[Random $(\eta, q, n)$-bipartite structure]\label{def:randombipartification}
        Given $\eta,q \in [0,1]$ and a vertex set $V$ of size $n\in 3\ZZ^+$, we define a \emph{random $(\eta,q,n)$-bipartite structure $\Xi_{\eta,q}=(X,Y,Z,M,E)$} as a quintuplet on the vertex set $V$ chosen in the following order:
        \begin{enumerate}[label=(C\arabic*)]
        \item\label{R1:tripartition} $V=X\cup Y \cup Z$ is an equipartition chosen uniformly at random;
        \item\label{R2:matching} $M\subseteq X\times Y$ is a perfect matching from $X$ to $Y$ chosen uniformly at random;
        \item\label{R3:graph} For a weight function $\phi: M\times Z \rightarrow [0,1]$ satisfying
        \begin{align*}
            \phi(x,y,z)=(1\pm \eta)q
        \end{align*}
        for $\{x,y,z\}\in M\times Z$, the $3$-graph $E\subseteq M\times Z$ is the random inhomogeneous $3$-graph $(M\times Z)_\phi$ given by
        \begin{align*}
            \PP\left(\{x,y,z\}\in E\right)=\phi(x,y,z).
        \end{align*}
        \end{enumerate}
    \end{definition}
    
    The main result in this section is that random $(\eta,q,n)$-bipartite structures project pseudorandom $3$-graphs into pseudorandom bipartite $2$-graphs with the appropriate density.

    \begin{lemma}\label{lem:bipartification}
        Let $\epsilon, \eta, p, q\in [0,1]$ be real numbers and $n\in 3\ZZ^+$ be a positive integer multiple of $3$ such that $n^{-1/5}\leq 3\epsilon\leq \frac{1}{100}$, $\eta\leq \min\{\epsilon, 1/4\}$, $p\geq \frac{C\log n}{\epsilon^2 n}$ and $q\geq \frac{C\log n}{\epsilon^4 pn}$ for a sufficiently large constant $C$. Let $H$ be a $3$-graph on $n$ vertices and $\Xi_{\eta,q}=(X,Y,Z,M,E)$ be a random $(\eta,q,n)$-bipartite structure on the same set of vertices. Then the following two statements hold with probability $1-n^{-C^{1/3}}$:
        \begin{enumerate}[label=(\roman*)]
            \item\label{bipartif1:regular} If $\deg_H(u,v)=(1\pm \epsilon)pn$ for every $u,v \in V(H)$, then
            \begin{align*}
            \deg_{\Xi_{\eta,q}^{(2)}(H)}(x)=(1\pm 100\epsilon)pqn/3
            \end{align*}
            for $x\in V(\Xi_{\eta,q}^{(2)}(H))$.
            \item\label{bipartif2:pseudo} If $H$ is an $(\epsilon,p,n)$-pseudorandom $3$-graph, then the bipartite graph $\Xi_{\eta,q}^{(2)}(H)$ is $(100\epsilon,pq,n/3)$-pseudorandom.
        \end{enumerate}
    \end{lemma}

    \subsection{The matching intersection lemma}
        
    A key ingredient to prove Lemma \ref{lem:bipartification} is the following concentration result about the intersection of an almost regular bipartite graph with a perfect matching chosen uniformly at random.

     \begin{lemma}\label{lem:switching}
        Let $\epsilon, p \in [0,1]$ be real numbers such that $0\leq \epsilon \leq \frac{1}{100}$ and $p\geq C\frac{\log n}{\epsilon^2 n}$ for a sufficiently large constant $C>0$. Let $B$ be a bipartite graph with parts $X,Y$ such that $|X|=|Y|=n$ and satisfying $$e(B)=(1\pm \epsilon)pn^2.$$ for every $z\in V(B)$. If $M$ is a perfect matching chosen uniformly at random from the set of all perfect matchings from $X$ to $Y$, then with probability at least $1-n^{-\sqrt{C}}$ it holds that $$|B\cap M|=(1\pm 6\epsilon)pn.$$
    \end{lemma}
\begin{proof}
    Let $M$ be a uniformly randomly chosen perfect matching between $X$ and $Y$. By ordering $X$ and $Y$, we can identify $M$ with the permutation $\pi_M \in S_n$. Here, $\pi_M$ is the permutation such that the edges of $M$ are given by $x_i,y_{\pi_M(i)}$.

    Let $Z = |B \cap M|$ be the random variable equal to the size of the intersection between $B$ and $M$. We plan to apply Theorem \ref{thm:mcdiarmid} to concentrate $Z$.  Since $B$ is fixed, the value of $Z$ is determined entirely by the randomly chosen permutation $\pi_M$. If two elements of $\pi_M$ (and $M$) are swapped, $Z$ changes by at most 2, since only 2 edges of $M$ are altered, so we can take $c=2$. Additionally, if $Z(\pi_M) \geq s$, this result is certified by any set of $s$ edges from $B \cap M$. This means we can take $r=1$.
    
    There is only one problem: to use this theorem effectively, we need to somehow estimate the median value of $Z$. Since every edge of $B$ belongs to $M$ with probability $1/n$, if $\mu=\E(Z)$, we have
    \begin{equation}\label{eq:mudef}\mu= \frac{e(B)}{n}= (1 \pm \eps)pn.\end{equation}
 Let $m$ be a median of $Z$. Our goal is to establish that $$|m-\mu| \leq \eps pn.$$
To do this, we prove a different concentration inequality for $Z$
 using a second moment argument. First, since the probability of two non-adjacent edges being included in $B$ is $\frac{1}{n(n-1)}$, 
$$\ex\left(Z^2 \right) \leq \frac{e(B)^2}{n(n-1)}.$$
Each edge is adjacent to at most $n$ other edges, so 
$$\frac{e(B)(e(B)-n)}{n(n-1)}\leq\ex\left(Z^2 \right).$$
Therefore, we have $$ \mu^2\frac{e(B)-n}{e(B)}\leq \ex(Z^2) \leq \mu^2\frac{n}{n-1}.$$
Note that $\frac{2n}{e(B)} \geq \max\left\{\frac{n}{e(B)},\frac{1}{n-1}\right\}$. Therefore,

$$\ex\left(Z^2-\mu^2 \right) \leq \mu^2\frac{2n}{e(B)} =2\mu.$$
We now obtain the concentration bound we wanted. By Chebyshev's inequality,
\begin{equation}\label{eq:badswitchinginequality}\Prob\left(|Z-\mu| \geq 2\sqrt{\mu}\right) \leq 1/2.\end{equation}
Using (\ref{eq:badswitchinginequality}), we can see that there is a median value $m$ for $Z$ such that $$|m-\mu| \leq 2\sqrt{\mu} \leq \eps pn.$$
The second inequality above comes from (\ref{eq:mudef}) and the fact that $p\geq C\frac{\log n}{\epsilon^2 n}$. 
We can finally apply Theorem \ref{thm:mcdiarmid} with $t=\eps pn$, $c=2$, and $r=1$ to obtain
\[
    \Prob(|Z-m|\geq \eps pn)
    \leq
    4\exp\left(
        -\frac{(\eps pn)^2}{64(m+\eps pn)}
    \right).
\]
Using (\ref{eq:badswitchinginequality}) and the fact that $p\geq C\frac{\log n}{\epsilon^2 n}$, we conclude:
\begin{align*}
    \Prob(|Z-\mu| \geq 2\eps pn) &\leq 4\exp\left(
        -\frac{(\eps pn)^2}{64(m+\eps pn)} 
    \right) \\
    &\leq 4\exp\left( -\frac{(\eps pn)^2}{65pn}\right) \\
    &\leq 4\exp \left(\frac{C}{65}\log n\right) \\
    &\leq n^{-\sqrt{C}}.
\end{align*}

\end{proof}
    
    \subsection{Proof of Lemma \ref{lem:bipartification}}

    We are now able to show the proof of the main result of this section.

    \begin{proof}[Proof of Lemma \ref{lem:bipartification}]
    We start by verifying \ref{bipartif1:regular}. Let $V(H)=[n]$. We expose the randomness from $\Xi_{\eta,q}$ in three steps. First, let $[n]=X\cup Y\cup Z$ be a random equipartition chosen uniformly at random with $|X|=|Y|=|Z|=n/3$. For any pair of vertices $\{u,v\}\subseteq [n]$, let $\deg_X(u,v)$, $\deg_{Y}(u,v)$ and $\deg_Z(u,v)$ be the degrees of the pair $\{u,v\}$ to the sets $X$, $Y$ and $Z$, respectively. Note that the random variables $\deg_X(u,v)$, $\deg_Y(u,v)$ and $\deg_Z(u,v)$ follow a hypergeometric distribution $\text{HG}(n-2, n/3, \deg_H(u,v))$ with mean 
    \begin{align*}
\mu=\E(\deg_X(u,v))=\E(\deg_Y(u,v))=\E(\deg_Z(u,v))=\frac{n\deg_H(u,v)}{3(n-2)}=(1\pm2\epsilon)pn/3.
    \end{align*}
    Therefore, by Theorem \ref{thm:hyperchernoff},
    \begin{align*}
        \PP\left(|\deg_X(u,v)-\mu|>\epsilon pn/3\right)\leq \exp\left(-\frac{\epsilon^2pn}{9(1+2\epsilon)^2}\right)\leq n^{-C/2}
    \end{align*}
    and similar bounds can be obtained for $\deg_Y(u,v)$ and $\deg_Z(u,v)$. Thus, after a union bound over all the pairs $\{u,v\}\in [n]^{(2)}$, we have with probability at least $1-n^{-C/4}$ that $\deg_X(u,v)=\deg_Y(u,v)=\deg_Z(u,v)=(1\pm3\epsilon)pn/3$. Condition on a triple $X\cup Y\cup Z$ such that this event occurs.
    
    Let $z\in Z$ be an arbitrary vertex and let $L_H(z)$ be the link of $z$ in $H$, i.e., the $2$-graph defined by
    \begin{align*}
        L_H(z):=\big\{\{x,y\} \in [n]^{(2)}:\: \{x,y,z\}\in H \big\}
    \end{align*}
    Let $B_z=L_H(z)[X,Y]$ be the bipartite graph with parts $|X|=|Y|=n/3$. The paragraph above implies that for $u\in X\cup Y$, we have $\deg_{B_z}(u)=(1\pm 3\epsilon)pn/3$, so $e(B_z) = (1 \pm 3\epsilon)pn^2/9$. Let $M$ be the perfect matching from $X$ to $Y$ chosen uniformly at random. By our hypothesis on $p$ and $\epsilon$, we obtain by Lemma \ref{lem:switching} that $|B_z\cap M|=(1\pm 18\epsilon)pn/3$ with probability $1-n^{-\sqrt{C}}$. From now on, condition on a perfect matching $M$ such that this event occurs.

    Let $\widetilde{\Xi}=(X,Y,Z,M,M\times Z)$ be the partite structure obtained by taking $E$ to be the whole $3$-graph $M\times Z$ and let $\widetilde{B}=\widetilde{\Xi}^{(2)}(H)$ be the bipartite projection of $H$ with respect to $\widetilde{\Xi}$. We claim that for every $u\in V(\widetilde{B})$ it holds that $\deg_{\widetilde{B}}(u)=(1\pm 18\epsilon)pn/3$. Indeed, if $e\in M$, then
    $\deg_{\widetilde{B}}(e)=\deg_Z(e)=(1\pm 3\epsilon)pn/3$. Otherwise, if $z\in Z$, then $\deg_{\widetilde{B}}(z)=|B_z\cap M|=(1\pm 18\epsilon)pn/3$. This finishes the proof of the claim.

    Finally, reveal the $3$-graph $E=(M\times Z)_{\phi}$. Note that by our procedure, the bipartite projection $\Xi_{\eta,q}^{(2)}(H)$ is just the random inhomogeneous graph $\widetilde{B}_\phi$. Then by Proposition \ref{prop:random_subbipartite_pseudo}.\ref{item:random_subbipartite_deg} we obtain for $x \in V(\Xi^{(2)}_{\eta,q}(H))$ that $\deg_{\Xi^{(2)}_{\eta,q}(H)}(x)=(1\pm100\epsilon)pqn/3$ with probability at least $1-n^{-\sqrt{C}}$. Since the probability that every event that we conditioned on fails is at most $n^{-\sqrt{C}}$, we obtain that \ref{bipartif1:regular} holds with probability at least $1-n^{-C^{1/3}}$.

    To verify \ref{bipartif2:pseudo}, we use that properties \ref{Q2:largesets}--\ref{Q4:smallZ} of being pseudorandom are hereditary. Indeed, we claim that if $H$ is an $(\epsilon,p,n)$-pseudorandom $3$-graph, then the bipartite projection $\widetilde{B}=\widetilde{\Xi}^{(2)}(H)$ (obtained after conditioning on the appropriate choices of $X\cup Y\cup Z$ and $M$) is $(18\epsilon,p,n/3)$-pseudorandom. Property \ref{P1:degree} follows from the fact that $\deg_{\widetilde{B}}(u)=(1\pm18\epsilon)pn/3$ for $u\in V(B)$. Properties \ref{P2:largesets}--\ref{P4:smallT} follow from \ref{Q2:largesets}--\ref{Q4:smallZ} where we take $S\subset M$ to be a subset of edges in the matching $M$ and $T\subseteq Z$ a subset of vertices in $Z$. Together with Remark \ref{rmk:pseudo_bipartite}, this concludes the proof of the claim. Since $\Xi^{(2)}_{\eta,q}(H)=\widetilde{B}_\phi$, it follows by Proposition \ref{prop:random_subbipartite_pseudo}.\ref{item:random_subbipartite_pseudo} with probability at least $1-n^{-\sqrt{C}}$ that $\Xi^{(2)}_{\eta,q}(H)$ is $(100\epsilon,pq,n/3)$-pseudorandom. The proof now follows since we condition $X$, $Y$, $Z$, $M$ and $E$ on the same events as in \ref{bipartif1:regular}.  
    \end{proof}

\section{Randomized algorithm}\label{sec:randomized}

In the last section, we describe our process to reduce the problem of finding a perfect matching in $3$-graphs to the problem of finding it in bipartite graphs. In this section, we study the related problem of how to obtain an almost $1$-factorization, i.e., a packing of disjoint perfect matchings on pseudorandom bipartite graphs. Although the problem of showing the existence of a perfect matching in pseudorandom bipartite graphs is a well-known consequence of Hall's theorem, for our application we will require that our perfect matchings are in some sense randomized. Indeed, we will show that if $B$ is a balanced pseudorandom bipartite graph, then there is a randomized algorithm that outputs an almost $1$-factorization with a certain probability distribution. 

Our approach is inspired by the work on the threshold for latin squares in \cites{kang2023thresholds, latinsquares24jain, keevash2022optimaledge}. One of the main technical steps in those works was to construct a probability distribution on $1$-factorizations for balanced dense regular bipartite graphs satisfying a so-called spreadness condition. In a certain sense, Lemma \ref{lem:randomizedalgorithm} can be seen as a variation of those results in the setting where the bipartite graph is not necessarily dense or regular.

\begin{lemma}\label{lem:randomizedalgorithm}
    Let $\epsilon,p, \gamma \in [0,1]$ be real numbers and $n\geq 1$ an integer satisfying $p\geq \frac{C\log n}{\epsilon^2 n}$ and $\epsilon\leq C^{-2}$ for a sufficiently large positive constant $C$. Let $B$ and $R$ be disjoint bipartite graphs on the vertex set $X, Y$ with $|X|=|Y|=n$ such that $B$ is $(\epsilon,p,n)$-pseudorandom and $R$ is $(\gamma\pm\epsilon)pn$-regular. Then there exists a randomized algorithm that with probability $1-n^{-C^{1/4}}$ outputs a collection of perfect matchings $\cM$ of $B\cup R$ with the following properties:
    \begin{enumerate}[label=(D\arabic*)]
        \item\label{itemR1} The edges of $R$ are covered by $\cM$, i.e., $R\subseteq \bigcup_{M\in \cM}M$.
        \item\label{itemR2} If $\widetilde{R}:=B\setminus \bigcup_{M\in\cM}M$, then for $x\in V(\widetilde{R})$, we have
        \begin{align*}
            C\epsilon pn\leq \deg_{\widetilde{R}}(x)\leq 12C\epsilon pn 
        \end{align*}
        \item\label{itemR3} Conditioned on the event that the algorithm succeeds, the probability that an edge is not covered by a matching in $\cM$ is
        \begin{align*}
            C\epsilon\leq \PP\left(e\notin\bigcup_{M\in \cM} M\right)\leq 5C\epsilon
        \end{align*}
        \item\label{itemR4} Conditioned on the event that the algorithm succeeds, for any set $A\subseteq B$ of $a$ edges, the probability that the set $A$ is contained in one of the matchings selected in $\cM$ is
        \begin{align*}
        \PP\left(\exists M\in \cM: A\subseteq M\right)\leq \left(\frac{C\log n}{pn}\right)^{a-1}
        \end{align*}
    \end{enumerate}
\end{lemma}

Before we introduce the algorithm used in Lemma \ref{lem:randomizedalgorithm}, we start with two preliminary results. The first is a standard corollary of Hall's theorem.

\begin{proposition}[\cite{diestel}, Corollary 2.1.3]\label{prop:regularfactorization}
    Let $d\geq 1$. Then every $d$-regular balanced bipartite graph $B$ has a $1$-factorization, i.e., there exists a collection $\{M_1,\ldots, M_d\}$ of disjoint perfect matchings in $B$ such that $B=\bigcup_{i=1}^d M_i$.
\end{proposition}

The second result is an application of the max-flow min-cut theorem obtained in \cite{kang2023thresholds}. We remark that even though in their proof they use a different concept of pseudorandomness, our definition implies theirs and we leave the details to Appendix \ref{appendix:kang}.

\begin{restatable}[\cite{kang2023thresholds}, Lemma 4.4]{lemma}{kangoreslemma} \label{lem:kangoreslemma}
    There exists a positive constant $\alpha_0\in  (0,1]$ such that the following holds: Let $\alpha, p, \gamma\in [0,1]$ and $n\geq 1$ such that $\alpha\leq \alpha_0$ and $\alpha pn\geq 1$. If $R$ and $B$ are disjoint bipartite graphs on the vertex set $X$, $Y$ with $|X|=|Y|=n$ such that $B$ is $(\alpha, p,n)$-pseudorandom and $R$ is $(\gamma\pm\alpha)pn$-regular, then there exists a regular spanning subgraph $G\subseteq B\cup R$ such that $R\subseteq G$.
\end{restatable}

We remark that both Proposition \ref{prop:regularfactorization} and Lemma \ref{lem:kangoreslemma} have efficient algorithmic versions using standard network flow algorithms (see, e.g., \cite{diestel}). We are going to assume these algorithmic versions in what comes next, even though they are not necessary for the understanding of our proof. We are now in a position to describe the randomized algorithm used in the proof of Lemma \ref{lem:randomizedalgorithm}.

Let $B$ and $R$ be disjoint balanced bipartite graphs where $B$ is $(\epsilon, p,n)$-pseudorandom and $R$ is $(\gamma\pm \epsilon)pn$-regular for some $\gamma \in [0,1]$. Let $C$ be a sufficiently large constant. We define the following parameters to be used in the rest of the discussion. Let $q$ be a real number and $L, k$ be integers satisfying
\begin{align}\label{eq:parameters}
    2C\epsilon\leq 2^{-L}\leq 4C\epsilon, \quad q:=\frac{1}{k\cdot 2^L},  \quad \text{and} \quad \frac{C\log n}{16}\leq pqn\leq \frac{C\log n}{8}.
\end{align}
First, as a preparation, we are going to partition the bipartite graph $B$ into $L+1$ parts $B=\bigcup_{i=1}^{L+1} B_i$, where for $1\leq i\leq L$, each part $B_i$ is further refined to $k\cdot 2^{L-i}$ parts $B_i=\bigcup_{j=1}^{k\cdot 2^{L-i}}B_{i,j}$ as follows: For $1\leq i\leq L$ and $1\leq j\leq k\cdot 2^{L-i}$, let 
\begin{align}\label{eq:q_ij}
q_{i,j}:=q=\frac{1}{k\cdot 2^L}, \quad  q_i:=\sum_{j=1}^{k\cdot 2^{L-i}}q_{i,j}=2^{-i} \quad \text{and} \quad q_{L+1}:=2^{-L}
\end{align}
be real numbers such that $\sum_{i=1}^L\sum_{j=1}^{k\cdot 2^{L-i}} q_{i,j}+q_{L+1}=1$. We then assign each edge of $B$ to a unique graph on $\bigcup_{i=1}^L\bigcup_{j=1}^{k\cdot 2^{L_i}}B_{i,j}\cup B_{L+1}$ with probability $q_{i,j}$ to be in $B_{i,j}$ and with probability $q_{L+1}$ to be in $B_{L+1}$. This in particular implies that each $B_{i,j}$ is sampled as the homogeneous random graph $B_q$ and for $1\leq i \leq L+1$, the bipartite graphs $B_{i}$ are sampled as the homogeneous random graph $B_{q_i}$. Alternatively, the partition can be seen as halving the bipartite graph $B$ at random for $L$ steps and in each step obtaining $B_i$ for $1\leq i \leq L$. We call the remaining graph $B_{L+1}$ and for $1\leq i \leq L$, each $B_i$ is further refined into $k\cdot 2^{L-i}$ parts sampled with equal probabilities. 

Let $\alpha:=\min\{\alpha_0,1/2\}$, where $\alpha_0$ is the constant given in Lemma \ref{lem:kangoreslemma}. Our first procedure is to verify that a series of bipartite graphs are pseudorandom. For $1\leq i \leq L+1$, let $B_{\geq i}=\bigcup_{i=1}^{L+1} B_i$. We verify pseudorandomness for the following three families
\begin{itemize}
    \item[\scalebox{0.5}{$\bullet$}] $B_{i,j}$ is $(\alpha, pq, n)$-pseudorandom for $1\leq i \leq L$ and $1\leq j\leq k\cdot 2^{L-i}$.
    \item[\scalebox{0.5}{$\bullet$}] $B_i$ is $(\alpha/8,p q_i, n)$-pseudorandom for $1\leq i \leq L+1$
    \item[\scalebox{0.5}{$\bullet$}] $B_{\geq i}$ is $(\alpha/8,p\cdot 2^{-(i-1)}, n)$-pseudorandom for $1\leq i \leq L+1$.
\end{itemize}
In case that one of those graphs fails pseudorandomness, we declare failure and stop immediately. After that, we start the main procedure of the algorithm by setting $R_1:=R$ and the collection of matchings $\cM$ to be the empty set. We construct our almost $1$-factorization in $L$ steps. At each step $1\leq i\leq L$, we obtain a set of matchings $\cM_i$ with the properties
\begin{enumerate}
    \item[(1)] For each $M\in \cM_i$, we have $M\subseteq R_i\cup B_i$
    \item[(2)] The edges of $R_i$ are covered by the matchings in $\cM_i$
\end{enumerate}
and a set $R_{i+1}\subseteq B_i$ of remaining edges to be covered in the next step. At the end of the process we return the collection of matchings $\cM=\bigcup_{i=1}^L \cM_i$.

Suppose that we already obtained $\cM_1, \ldots \cM_{i-1}$ and $R_1,\ldots, R_{i}$ and we want to construct $\cM_i$ and $R_{i+1}$. For each edge in $R_i$, we assign a label $j \in [k\cdot 2^{L-i}]$ uniformly at random, and let $R_{i,j}\subseteq R_i$ be the subgraph of $R_i$ with all edges labeled $j$. This partitions $R_i =\bigcup_{j=1}^{k\cdot 2^{L-i}}R_{i,j}$ into $k\cdot 2^{L-i}$ parts. For each $j\in [k\cdot 2^{L-i}]$, if the bipartite graph $R_{i,j}$ is $(\widetilde{\gamma}\pm\alpha)pqn$ for some $0\leq \widetilde{\gamma} \leq 1$, we apply Lemma \ref{lem:kangoreslemma} to $B_{i,j}$ and $R_{i,j}$ to obtain a regular subgraph $G_{i,j}\subseteq B_{i,j}\cup R_{i,j}$ containing $R_{i,j}$. Otherwise, we declare that our algorithm fails and stop immediately. In the case that there are no failures, we obtain a family of bipartite graphs $G_{i,j}$ that are all regular. By Proposition \ref{prop:regularfactorization}, we can then partition each $G_{i,j}$ into a collection of perfect matchings $\cM_{i,j}$. We set $\cM_{i}=\bigcup_{j=1}^{k\cdot 2^{L-i}}\cM_{i,j}$. Finally, we let $R_{i+1}$ be the remaining edges in $B_{i}$ that were not covered by matchings of $\cM_i$, i.e., $R_{i+1}=B_i\setminus \bigcup_{j=1}^{k\cdot 2^{L-i}}G_{i,j}$.

Note that by the end of the algorithm, all the edges of $B\cup R$ were covered except for $B_{L+1}\cup R_{L+1} \subseteq B_L\cup B_{L+1}$, which corresponds to roughly less than $2^{-(L-1)}$ proportion of the whole set of edges. We also remark that for the purpose of obtaining \ref{itemR1}, the partition into sets $B_i$ suffices, and the refinement $B_{i,j}$'s are intended to control the bounds on condition \ref{itemR3}. The following is a formal description of the algorithm.

\begin{algorithm}[H]\caption{{\bfseries Randomized algorithm}}\label{alg:randomized} 
\begin{algorithmic}[1]
\STATE Assign each edge in $B$ to a unique bipartite graph $B_{i,j}$ with probability $q_{i,j}$ or to the bipartite graph $B_{L+1}$ with probability $q_{L+1}$.
\IF{one of the $B_{i,j}$, $B_{i}$ or $B_{\geq i}$ is not pseudorandom}
     \STATE declare \textbf{failure}
\ENDIF
\STATE Set~$R_1:=R$ and~$\cM:=\emptyset$.
\FOR{$i = 1$ \textbf{to} $L$}
		\STATE Assign each edge in $R_i$ a label $j\in [k\cdot 2^{L-i}]$ uniformly at random, and let $R_{i,j}$ be the bipartite subgraphs obtained after the assignment.
        \FOR{$j=1$ \textbf{to} $k\cdot 2^{L-i}$}
            \IF{$R_{i,j}$ is not $(\widetilde{\gamma}\pm \alpha)pqn$-regular for $0\leq \widetilde{\gamma}\leq 1$}
            \STATE declare \textbf{failure}
            \ENDIF
            \STATE Apply Lemma \ref{lem:kangoreslemma} to $B_{i,j}$ and $R_{i,j}$ to obtain a subgraph $G_{i,j}\subseteq B_{i,j}\cup R_{i,j}$ containing $R_{i,j}$
            \STATE Apply Proposition \ref{prop:regularfactorization} to obtain a $1$-factorization $\cM_{i,j}$ of $G_{i,j}$
            \STATE Add $\cM_{i,j}$ to $\cM$
        \ENDFOR
        \STATE Let $R_{i+1}:=\bigcup_{j\in [k\cdot 2^{L-i}]} B_{i,j}\setminus G_{i,j}$.
\ENDFOR
\STATE Return~$\cM$.
\end{algorithmic}
\end{algorithm}

We remark that the algorithm above performs in polynomial time. Indeed, as discussed previously, Proposition \ref{prop:regularfactorization} and Lemma \ref{lem:kangoreslemma} have efficient algorithms using network flow algorithms. Moreover, properties \ref{P2:largesets}--\ref{P4:smallT} of pseudorandomness are hereditary and hence verifying that $B_i$, $B_{\geq i}$ and $B_{i,j}$ reduces to verifying that these graphs are almost regular, which can be done in polynomial time.

Although our algorithm can declare failure, in the following we show that it succeeds with high probability.

\begin{proposition}\label{prop:failure}
    The algorithm above declares failure with probability at most $n^{-C^{1/4}}$
\end{proposition}

\begin{proof}
    There are two points in the algorithm that can produce a failure. The first one is if one of the subgraphs $B_i$, $B_{i\geq i}$ or $B_{i,j}$ is not pseudorandom. The second one is if the bipartite graph $R_{i,j}$ is not $(\widetilde{\gamma}\pm \alpha)pqn$-regular for some $\widetilde{\gamma}\in [0,1]$. We analyze the graphs $B_{i,j}$ first. Let $\alpha:=\min\{\alpha_0,1/2\}$, where $\alpha_0$ is the constant given by Lemma \ref{lem:kangoreslemma}. By our choice of parameters in (\ref{eq:parameters}) and (\ref{eq:q_ij}), for $1\leq i \leq L$ and $1\leq j\leq k\cdot 2^{L-i}$, the graph $B_{i,j}$ has the same distribution as $B_q$. Since $B$ is $(\epsilon, p, n)$-pseudorandom and
    \begin{align*}
    \alpha^4pqn\geq C\alpha^4\log n/16,
    \end{align*}
    we have by Proposition \ref{prop:random_subbipartite_pseudo} that $B_{i,j}$ is $(\alpha,pq,n)$-pseudorandom  with probability $1-n^{-\alpha^2\sqrt{C}/4}$. Let $\cE$ be the event that for every $1\leq i \leq L$ and $1\leq j\leq k\cdot 2^{L-i}$ the graphs $B_{i,j}$ are pseudorandom. Then the last observation combined with a union bound gives that
    \begin{align*}
        \PP(\cE^{\mathsf{c}})\leq k\cdot 2^L \cdot n^{-\alpha^2\sqrt{C}/4}\leq n^{-C^{1/3}} 
    \end{align*}
    for sufficiently large $C$.

    For each $1\leq i\leq L+1$, let $B_{\geq i}=\bigcup_{t=i}^{L+1} B_t$ be the collection of edges in $B$ not used at steps $1$ to $i-1$. Similarly as in the last paragraph, by our choice of parameters in (\ref{eq:parameters}) and (\ref{eq:q_ij}) we have for $1\leq i\leq L+1$ that the graphs $B_i$ and $B_{\geq i}$ have the same distribution as the random homogeneous subgraph of $B$ where each edge is chosen with probability $q_i$ and $\sum_{t=i}^{L+1}q_i=2^{-(i-1)}$, respectively. Since $B$ is $(\epsilon, p, n)$-pseudorandom and
    \begin{align*}
        \frac{\alpha^4 pn}{8^4\cdot 2^{i-1}}\geq \frac{\alpha^4 pn}{8^4\cdot 2^L}> \frac{C\alpha^4\log n}{16^4}
    \end{align*}
    we have by Proposition \ref{prop:random_subbipartite_pseudo} with probability $1-n^{-\alpha^2\sqrt{C}/256}$ that $B_{\geq i}$ is $(\alpha/8,p\cdot 2^{-(i-1)},n)$-pseudorandom. Similar calculations show with probability $1-n^{-\alpha^2\sqrt{C}/256}$ that $B_i$ is $(\alpha/8,pq_i,n)$-pseudorandom. Let $\cF$ be the event that for $1\leq i\leq L+1$ the graphs $B_i$ and $B_{\geq i}$ are pseudorandom. The last observation combined with a union bound gives that $\PP(\cF^{\mathsf{c}})\leq n^{-C^{1/3}}$ for sufficiently large $C$. 

    We now turn our attention to the degrees of $R_{i,j}$. For $1\leq t \leq L$, let $\cG_t$ be the event that the algorithm has not declared failure until the end of step $t$. Suppose that we are currently at the beginning of step $i$, for $1\leq i \leq L$, conditioned on the events $\cE$, $\cF$ and $\cG_{i-1}$, i.e., all graphs $B_{i,j}$, $B_{\geq i}$ are pseudorandom and the algorithm has not declared failure yet. We claim that for some $0\leq \gamma_i\leq 1$, the graph $R_i$ is $(\gamma_i\pm\alpha/2)pq_{i-1}n$-regular. Note that the claim holds trivially from the statement hypothesis for $i=1$, so we may assume that $i\geq 2$. For $1\leq t \leq i-1$, let $G_t=\bigcup_{j\in [k\cdot 2^{L-t}]} G_{t,j}$. Note that at the beginning of step $i$, the edges of $B_{\geq i}$ were untouched and our algorithm already covered all the edges from $R\cup \bigcup_{t=1}^{i-1} B_t$ except for the edges in $R_{i}$. Therefore, it holds that
    \begin{align*}
        R_{i}=B_{i-1}\setminus G_{i-1}=\left(R\cup B\right) \setminus\left(B_{\geq i}\cup \bigcup_{t=1}^{i-1} G_t\right).
    \end{align*}
    Since each step of the algorithm so far produced a regular spanning graph $G_{t}$, the hypothesis of the lemma together with the fact that $B_{\geq i}$ is $(\alpha/8,p\cdot 2^{-(i-1)},n)$-pseudorandom and (\ref{eq:parameters}) and (\ref{eq:q_ij}) gives us that
    \begin{align}\label{eq:deggapR_i}
        \Delta(R_i)-\delta(R_i)&\leq \left(\Delta(R)-\delta(R)\right) + \left(\Delta(B)-\delta(B)\right) + \left(\Delta(B_{\geq i})-\delta(B_{\geq i})\right)\nonumber\\
        &\leq 2\epsilon pn + 2\epsilon pn +\frac{2\alpha pn}{8\cdot 2^{i-1}} < \frac{2\alpha pn}{4\cdot 2^{i-1}}=\alpha pq_{i-1}n/2.
    \end{align}
    Moreover, because $R_i\subseteq B_{i-1}$ and $B_{i-1}$ is $(\alpha/8,pq_{i-1}, n)$-pseudorandom, the relation in (\ref{eq:deggapR_i}) gives us that 
    \begin{align*}
        \deg_{R_i}(x)=\left(\gamma_i\pm \frac{\alpha}{4}\right)pq_{i-1}n,
    \end{align*}
    for every $x\in V(R_i)$, where $0\leq \gamma_i\leq 1$. This finishes the proof of the claim.

    Finally, for $j\in [k\cdot 2^{L-i}]$, we construct the graphs $R_{i,j}$ by selecting uniformly at random a label for every edge in $R_i$. That is, every $R_{i,j}$ is distributed as the random homogeneous subgraph of $R_i$, where each edge is chosen with probability $(k\cdot 2^{L-i})^{-1}$. Since $R_i$ is $(\gamma_i+\alpha/4)$-pseudorandom and
    \begin{align*}
        \frac{\alpha^4p q_{i-1} n}{k\cdot 2^{L-i}}=\frac{\alpha^4 pn}{k\cdot 2^L}=\alpha^4pqn\geq C\alpha^4\log n/16,
    \end{align*}
    we have by Proposition \ref{prop:random_subbipartite_pseudo} with probability $1-n^{-\alpha^2\sqrt{C}/4}$ that $R_{i,j}$ is $(\gamma_i \pm \alpha)pqn$-regular. A union bound shows that with probability $1-n^{-C^{1/3}}$, it holds that every $R_{i,j}$ is $(\gamma_i \pm \alpha)pqn$-regular. Therefore, because we are conditioning on the event $\cE$, for every $j\in [k\cdot 2^{L_i}]$ the bipartite graphs $B_{i,j}$ and $R_{i,j}$ satisfy the hypothesis of Lemma \ref{lem:kangoreslemma} for the chosen constant $\alpha$. That implies that $\PP(\cG_i^\mathsf{c}\mid \cE, \cF, \cG_{i-1})\leq n^{-C^{1/3}}$. Consequently, the probability that the algorithm declares failure is at most
    \begin{align*}
        \PP(\cG_L^{\mathsf{c}})\leq \sum_{i=1}^L\PP(\cG_i^{\mathsf{c}}\mid \cE, \cF, \cG_{i-1})+\PP(\cE^\mathsf{c})+\PP(\cF^{\mathsf{c}})\leq n^{-C^{1/4}}
    \end{align*}
    for sufficiently large $C$. This concludes the proof of the proposition.
\end{proof}

We are now able to finish the proof of Lemma \ref{lem:randomizedalgorithm}.

\begin{proof}[Proof of Lemma \ref{lem:randomizedalgorithm}]
    Consider the algorithm described in Algorithm \ref{alg:randomized} and let $\cG$ be the event that it succeeded. By Proposition \ref{prop:failure}, we have that $\PP(\cG)\geq 1-n^{-C^{1/4}}$. Therefore, for every event $\cA$ we have that
    \begin{align}\label{eq:correction}
        \PP(\cA)-n^{-C^{1/4}}\leq \PP(\cA)-\PP(\cG^{\mathsf{c}})\leq\PP(\cA\mid \cG)\leq \frac{\PP(\cA)}{\PP(\cG)}\leq \left(1+ n^{-C^{1/4}}\right)\PP(\cA)
    \end{align}
    Let $\cM=\bigcup_{i=1}^L M_i$ be the set of perfect matchings produced as the output of the algorithm conditioned on $\cG$, where $\cM_i$ is the set of matchings produced at step $i$. By the description of the procedure, it is clear that $R=R_1\subseteq \cM_1\subseteq \cM$. Hence, property \ref{itemR1} follows immediately.

    Note that by the end of the algorithm, we obtain as an output a subgraph $R_{L+1} \subseteq B_L$ and the edges of $B_{L+1}$ were never touched. This implies in particular that $\widetilde{R}=R_{L+1}\cup B_{L+1}$ and moreover that $B_{L+1}\subseteq \widetilde{R}\subseteq B_{L}\cup B_{L+1}$. Since the algorithm did not declare failure, then both $B_L$ and $B_{L+1}$ are $(1\pm \alpha)pq_{L+1}\cdot n$-regular. Thus, for every $x \in V(\widetilde{R})$, we have by (\ref{eq:parameters}) and (\ref{eq:q_ij}) that
    \begin{align*}
        \deg_{\widetilde{R}}(x)\leq \deg_{B_{L}}(x)+\deg_{B_{L+1}}(x)\leq 2(1+\alpha)pq_{L+1}\cdot n\leq 12C\epsilon pn.
    \end{align*}
    and that
    \begin{align*}
        \deg_{\widetilde{R}}(x)\geq \deg_{B_{L+1}}(x)\geq (1-\alpha)pq_{L+1}\cdot n \geq C\epsilon pn
    \end{align*}
    This concludes the proof of \ref{itemR2}.

    To verify \ref{itemR3}, note that the only edges that are not covered by $\cM$ are the edges in $\widetilde{R}$. Since $B_{L+1}\subseteq \widetilde{R}\subseteq B_{L}\cup B_{L+1}$, then
    \begin{align*}
        \PP(e\in B_{L+1}\mid \cG)\leq \PP\left(e\notin \bigcup_{M\in \cM}M \,\,\middle| \,\,\cG\right) \leq \PP(e\in B_{L}\cup B_{L+1}\mid \cG).
    \end{align*}
    The result now follows by (\ref{eq:correction}), the fact that $C$ is sufficiently large and because $B_L$ and $B_{L+1}$ are disjoint random homogeneous subgraphs of $B$ sampled with probability $q_{L+1}=2^{-L}$ and $L$ is chosen such that $2C\epsilon\leq 2^{-L}\leq 4C\epsilon$.

    It remains to verify \ref{itemR4}. For that we will need some preparation. Let
    \begin{align*}
        I=\{(i,j):\: 1\leq i \leq L,\, 1\leq j \leq k\cdot 2^{L-i}\} \cup\{L+1\}
    \end{align*}
    be the set of $1+\sum_{i=1}^Lk\cdot 2^{L-i}=1+k(2^L-1)$ indices corresponding to all bipartite graphs $B_{i,j}$ and the remaining $B_{L+1}$. For every edge $e\in B$, let $X_e \in I$ be a random variable chosen such that
    \begin{align*}
        \PP(X_e=(i,j))=q_{i,j} \quad \text{and} \quad \PP(X_e=L+1)=q_{L+1},
    \end{align*}
    where $q_{i,j}$ and $q_{L+1}$ are defined in (\ref{eq:q_ij}). Moreover, for every edge $e\in B$ and $1\leq i \leq L$, let $Y_{e}^{(i)} \in [k\cdot 2^{L-i}]$ be a random variable chosen independently uniformly at random. We claim that the randomized selections in Algorithm \ref{alg:randomized} can be described by the random variables $X_e$, $Y_e^{(i)}$. Indeed, note that $X_e$ just indicates that the edge $e$ belongs to the bipartite graph $B_{X_e}$. Moreover, if we are on step $i$ of the algorithm and $e\in R_i$, then $Y_e^{(i)}$ indicates the label that the edge $e$ is assigned, i.e., the edge belongs to the graph $R_{i,Y_e^{(i)}}$. Note that in the case that $e \in R_i$ or that the algorithm already produced a failure at that point, the random choice of $Y_{e}^{(i)}$ is still there, even if it is not used in the algorithm.    
    
    Let $A$ be an arbitrary set of edges in $B$. For every $1\leq i \leq L$ and $1\leq j\leq k\cdot 2^{L-i}$, we want to compute $\PP(\exists M\in \cM_{i,j}:\: A\subseteq M\mid \cG)$. Note that this event can be bounded by the event $\{A\subseteq G_{i,j} \mid \cG\}$. For an edge $e\in A$, there are two ways that it can be in $G_{i,j}$
    \begin{enumerate}
        \item[$\cA_e$:] The edge $e \in B_{i,j}$ and the application of Lemma \ref{lem:kangoreslemma} resulted in a graph $G_{i,j}$ containing it. In this case, the choice $X_e=(i,j)$ must have happened.
        \item[$\cB_e$:] The edge $e\in B_{i-1}$ and in step $(i-1)$, the edge was not covered by $G_{i-1}$, i.e., $e$ ended up in the graph $R_i$. Moreover, the edge $e$ was assigned to $R_{i,j}$ and then covered by $G_{i,j}$ after Lemma \ref{lem:kangoreslemma}. In this case, the choice $X_e \in I_{i-1}:= \bigcup_{t=1}^{k\cdot 2^{L-(i-1)}} (i-1,t)$ and $Y_e^{(i)}=j$ must have happened.
    \end{enumerate}
    Since the random variables $X_e$, $Y_e^{(i)}$ are all independent for every $1\leq i \leq L$ and $e\in B$, we obtain by (\ref{eq:parameters}), (\ref{eq:q_ij}) and (\ref{eq:correction}) that
    {\footnotesize
    \begin{align*}
        \PP(\exists M\in \cM_{i,j}: A\subseteq M\mid \cG)&\leq \PP\left(\bigcap_{e\in A}\cA_e\cup B_e\,\,\middle|\,\,\cG\right)\leq \PP\left(\bigcap_{e\in A}\{X_e=(i,j)\}\vee\left\{X_e\in I_{i-1}\wedge Y_{e}^{(i)}=j\right\}\,\,\middle|\,\, \cG\right)\\
        &=\left(1+2n^{-C^{1/4}}\right)\prod_{e\in A}\left(q+q_{i-1}\cdot \frac{1}{k\cdot 2^{L-i}}\right)\leq 2\cdot (3q)^a.
    \end{align*}
    }
    Hence, by a union bound over all $k(2^L-1)$ possible indices $(i,j)$ and (\ref{eq:parameters}) we obtain that
    \begin{align*}
        \PP(\exists M\in \cM: A\subseteq M\mid \cG)\leq 2k\cdot (2^L-1)\cdot (3q)^a\leq (3q)^{a-1}\leq \left(\frac{C\log n}{pn}\right)^{a-1}.
    \end{align*}
    This concludes the proof of the lemma.
\end{proof}

We finish the section with an application of Lemma \ref{lem:randomizedalgorithm} to obtain a randomized algorithm that outputs a perfect matching on a pseudorandom $3$-graph with a distribution that is very close to being uniform over all the edges and pairs of vertices.

\begin{lemma}\label{lem:randomized_hypergraph_matching}
Let $\epsilon,p\in [0,1]$ be real numbers and let $n\in 3\ZZ^+$ be a positive integer multiple of $3$ such that $n^{-1/5}\leq \epsilon\leq C^{-3}$ and $p\geq \frac{C\log n}{\epsilon^2 n}$ for a sufficiently large constant $C>0$. Let $H$ be an $(\epsilon,p,n)$-pseudorandom $3$-graph on $n$ vertices. Then there exists a randomized algorithm that with probability at least $1-n^{-C^{1/5}}$ outputs a perfect matching $M$ of $H$ satisfying the following properties:
\begin{enumerate}[label=(E\arabic*)]
    \item\label{itemM1} For every edge $e\in H$, it holds that $\PP(e\in M)=2(1\pm 10C\epsilon)/pn^2$.

    \item\label{itemM2} For every pair $\{u,v\}\subseteq V(H)$, the probability that there exists an edge $e\in M$ containing $\{u,v\}$ satisfies
    \begin{align*}
        \PP\left(\exists e\in M:\{u,v\}\subseteq e\right)
        =2(1\pm10 C\epsilon)/n.
    \end{align*}

    \item\label{itemM3} For every set of edges $A\subseteq H$ of size $a$, it holds that
    \begin{align*}
        \PP(A\subseteq M)\leq
        \left(\frac{9C\log n}{pn^2}\right)^{a}.
    \end{align*}
\end{enumerate}
\end{lemma}

\begin{proof}
We define a randomized algorithm as follows. First, we invoke Definition \ref{def:randombipartification} and let $\Xi_{0,1}=(X_{0,1},Y_{0,1},Z_{0,1},M_{0,1},E_{0,1})$
be a random $(0,1,n)$-bipartite structure on $V(H)$. In other words, $V(H)=X_{0,1}\cup Y_{0,1}\cup Z_{0,1}$ is a random equipartition with $|X_{0,1}|=|Y_{0,1}|=|Z_{0,1}|=n/3$, the set $M_{0,1}\subseteq X_{0,1}\times Y_{0,1}$ is a perfect matching chosen uniformly at random, and, since the parameters are chosen such that $\eta=0$ and $q=1$, we have $E_{0,1}=M_{0,1}\times Z_{0,1}$. Let $B:=\Xi_{0,1}^{(2)}(H)$ be the bipartite projection of $H$ with respect to $\Xi_{0,1}$. 

Note that, since $H$ is $(\epsilon, p, n)$-pseudorandom, Lemma \ref{lem:bipartification} says with high probability that $B$ is $(100\epsilon, p, n/3)$-pseudorandom. If the graph $B$ is not pseudorandom, then the algorithm declares failure. Otherwise, we apply Lemma \ref{lem:randomizedalgorithm} to $B$ with $R=\emptyset$. If the algorithm from Lemma \ref{lem:randomizedalgorithm} declares failure, then our algorithm also declares failure. Otherwise, let $\widetilde{\cM}=\{\widetilde{M}_1^{(2)},\ldots, \widetilde{M}_t^{(2)}\}$ be the collection of $t$ perfect matchings of $B$ obtained as an output. We choose an index $i_{\star}\in [t]$ uniformly at random and let $M:=M^{(3)}$ be the $3$-uniform matching in $\Xi_{0,1}^{(3)}(H)$ obtained by lifting the perfect matching $\widetilde{M}_{i_{\star}}$ through the projection map. Since $\widetilde{M}_{i_{\star}}$ is a perfect matching of the bipartite graph with parts $M_{0,1}$ and $Z_{0,1}$, its lift is a perfect matching of $H$. Here is a formal description of the algorithm.

\begin{algorithm}[H]\caption{{\bfseries Randomized perfect matching algorithm}}\label{alg:randomized_hypergraph_matching}
\begin{algorithmic}[1]
\STATE Let $\Xi_{0,1}=(X_{0,1},Y_{0,1},Z_{0,1},M_{0,1},E_{0,1})$ be a random $(0,1,n)$-bipartite structure on $V(H)$ and let $B:=\Xi_{0,1}^{(2)}(H)$ be the bipartite projection of $H$ with respect to $\Xi_{0,1}$.
\IF{$B$ is not $(100\epsilon,p,n/3)$-pseudorandom}
    \STATE declare \textbf{failure}
\ENDIF
\STATE Apply Lemma \ref{lem:randomizedalgorithm} to $B$ with $R=\emptyset$.
\IF{Lemma \ref{lem:randomizedalgorithm} declares failure}
    \STATE declare \textbf{failure}
\ENDIF
\STATE Let $\widetilde{\cM}=\{\widetilde{M}^{(2)}_1,\ldots,\widetilde{M}^{(2)}_t\}$ be the set of $t$ perfect matchings obtained from Lemma \ref{lem:randomizedalgorithm}.
\STATE Choose $i_{\star}\in [t]$ uniformly at random.
\STATE Let $M$ be the $3$-uniform matching in $\Xi_{0,1}^{(3)}(H)$ obtained by lifting $\widetilde{M}^{(2)}_{i_\star}$ through the projection map.
\STATE Return~$M$.
\end{algorithmic}
\end{algorithm}

Note that for sufficiently large $C$, our choice of parameters ensures that $H$ satisfies the hypothesis of Lemma \ref{lem:bipartification}. Therefore, by Lemma \ref{lem:bipartification}, with probability at least $1-n^{-C^{1/3}}$ the graph $B=\Xi_{0,1}^{(2)}(H)$ is $(100\epsilon,p,n/3)$-pseudorandom. Moreover, conditioned on this event, again for sufficiently large $C$, the bipartite graph $B$ satisfies the hypothesis of Lemma \ref{lem:randomizedalgorithm}, and consequently the lemma succeeds with probability at least $1-(n/3)^{-C^{1/4}}$. Hence, the probability that the algorithm declares failure is at most $n^{-C^{1/5}}$.

For the rest of the proof all the probabilities are conditioned on the event that the algorithm succeeds. We now verify properties \ref{itemM1}--\ref{itemM3}. Since $B$ is $(100\epsilon,p,N)$-pseudorandom, every vertex of $B$ has degree $(1\pm 100\epsilon)pn/3$. Moreover, if $\widetilde{R}:=B\setminus\bigcup_{i=1}^t \widetilde{M}^{(2)}_i$ is the remainder graph not covered by the matchings in $\widetilde{\cM}$, then, by property \ref{itemR2} of Lemma \ref{lem:randomizedalgorithm}, we have $\Delta(\widetilde{R})\leq 16\epsilon pn/3$. Since every perfect matching $\widetilde{M}_i^{(2)}$ covers every vertex of $B$ exactly once, for every vertex $v\in V(B)$ we have $\deg_B(v)=t+\deg_{\widetilde{R}}(v)$. Thus, we can bound $t$ by
\begin{align}\label{eq:number_matchings_hypergraph_algorithm}
    t=\deg_B(v)-\deg_{\widetilde{R}}(v)=(1\pm 200\epsilon)pn/3.
\end{align}

We will also use the following consequence of Lemma \ref{lem:randomizedalgorithm}. Let $f\in M_{0,1}\times Z_{0,1}$ be a fixed edge. Let $\cA_f$ be the event corresponding to the choices of the random bipartite structure $\Xi_{0,1}$ for which $f\in B$. Thus, conditioned on $\cA_f$, Property \ref{itemR3} of Lemma \ref{lem:randomizedalgorithm} gives us that 
\[
    1-5C\epsilon\leq \PP\left(f\in \bigcup_{i=1}^t \widetilde{M}^{(2)}_i \,\, \middle |\,\, \cA_f\right)\leq 1- C\epsilon.
\]
Together with \eqref{eq:number_matchings_hypergraph_algorithm}, this implies for sufficiently large $C$ that
\begin{align}\label{eq:projected_edge_selected}
    \PP(f\in \widetilde{M}^{(2)}_{i_\star} \mid \cA_f)=\frac{1}{t}\cdot \PP\left(f\in \bigcup_{i=1}^t \widetilde{M}^{(2)}_i \,\, \middle | \,\, \cA_f\right)=3(1\pm 6C\epsilon)/pn.
\end{align}

We first verify \ref{itemM1}. Fix an edge $e\in H$. In order for $e$ to be an edge of the output matching, it is necessary that $e \in \Xi_{0,1}^{(3)}(H)$. That is, the edge $e$ must contain exactly one vertex in each of $X_{0,1},Y_{0,1},Z_{0,1}$, and the two vertices of $e$ that lie in $X_{0,1}\cup Y_{0,1}$ must form an edge of $M_{0,1}$. Since the partition $X_{0,1}\cup Y_{0,1}\cup Z_{0,1}$ is chosen uniformly at random from all equipartitions, the probability that $e$ contains exactly one vertex in each of $X_{0,1},Y_{0,1},Z_{0,1}$ is $\frac{2n^2}{9(n-1)(n-2)}=2(1\pm \epsilon)/9$. Conditioned on this event, the probability that the two vertices in $X_{0,1}\cup Y_{0,1}$ are matched by $M_{0,1}$ is $3/n$. Therefore, the probability that $e$ appears as an edge of $\Xi_{0,1}^{(3)}(H)$ is $2(1\pm \epsilon)/3n$.

Let $f_e \in B$ be the edge corresponding to $e$ in the projection map $\Xi^{(3)}(H)\mapsto\Xi^{(2)}(H)$. For $e \in M$, we also need that $f_e$ is selected in the uniformly chosen matching $\widetilde{M}_{i_\star}$. By \eqref{eq:projected_edge_selected}, the probability that $f_e$ is in matching $\widetilde{M}_{i_\star}$ is $3(1\pm 6C\epsilon)/pn$. Hence, we have
\begin{align*}
    \PP(e\in M)=\PP\left(f_e\in \widetilde{M}_{i\star}\,\,\middle|\,\, e\in \Xi^{(3)}(H)\right)\cdot \PP\left(e\in \Xi^{(3)}(H)\right)=2(1\pm 10C\epsilon)/pn^2.
\end{align*}
This concludes the proof of \ref{itemM1}.

We now prove \ref{itemM2}. Fix a pair $\{u,v\}\subseteq V(H)$. There are two ways in which the output matching $M$ can contain an edge containing $\{u,v\}$.

\vspace{0.2cm}

\noindent \underline{Case 1:} The pair $\{u,v\}$ is an edge of $M_{0,1}$.

\vspace{0.2cm}

This happens precisely when one of $u,v$ lies in $X_{0,1}$, the other lies in $Y_{0,1}$, and they are matched by $M_{0,1}$. The probability of this event is given by
\begin{align}\label{eq:pair_case1}
    \frac{2n}{9(n-1)}\cdot \frac{3}{n}
    =\frac{2}{3(n-1)}=\frac{2(1\pm\epsilon)}{3n}.
\end{align}
If $\{u,v\}\in M_{0,1}$, then every perfect matching $\widetilde{M}_i$ of $B\subseteq M_{0,1}\times Z_{0,1}$ contains exactly one edge incident to the vertex $\{u,v\}\in M_{0,1}$. Therefore, the lifted matching $M$ has an edge containing~$\{u,v\}$.

\vspace{0.2cm}

\noindent \underline{Case 2:} One vertex of $\{u,v\}$ lies in $Z_{0,1}$ and the other lies in $X_{0,1}\cup Y_{0,1}$.

\vspace{0.2cm}

Let $\cE$ be the event described in this case. Since the equipartition $X_{0,1}\cup Y_{0,1}\cup Z_{0,1}$ is chosen uniformly at random, one can see that $\PP(\cE)=4n/9(n-1)$. Suppose without loss of generality that $u\in X_{0,1}$ and $v\in Z_{0,1}$. Let
\begin{align*}
    \deg_{Y_{0,1}}(u,v):=\{w\in Y_{0,1}:\{u,v,w\}\in H\}
\end{align*}
be the degree of the pair $\{u,v\}$ in $Y_{0,1}$. Note that since $Y_{0,1}$ is a uniformly chosen random set of size $n/3$, the random variable $\deg_{Y_{0,1}}(u,v)$ follows a hypergeometric distribution with parameters $\textrm{HG}(n, \deg_H(u,v),n/3)$ and mean given by $\E(\deg_{Y_{0,1}}(u,v))=\deg_H(u,v)/3$. Moreover, because $H$ is $(\epsilon,p,n)$-pseudorandom, we have $\deg_H(u,v)=(1\pm \epsilon)pn$. Thus, by Theorem \ref{thm:hyperchernoff} it holds that
\begin{align*}
\PP\left(\left|\deg_{Y_{0,1}}(u,v)-\frac{\deg_H(u,v)}{3}\right|>\frac{\epsilon pn}{3}\right)\leq 2\exp\left(-\frac{\epsilon^2pn}{9(1+\epsilon)^2}\right)\leq n^{-\sqrt{C}}.
\end{align*}
Let $\cF$ be the event that $\deg_{Y_{0,1}}(u,v)=(1\pm 3\epsilon)pn/3$. The above computation implies that $\PP(\cF^{\mathsf{c}})\leq n^{-\sqrt{C}}$. We remark that a similar computation would hold if $\{u,v\}\subseteq Y_{0,1}\cup Z_{0,1}$ except that in this case we would be interested in the quantity $\deg_{X_{0,1}}(u,v)$.

We now condition on both events $\cE$ and $\cF$. Hence, the probability that the vertex $u$ is matched in $M_{0,1}$ to some vertex $w\in Y_{0,1}$ such that $\{u,v,w\}\in H$ is given by $3\deg_Y(u,v)/n=(1\pm 3\epsilon)p$. If this happens, then the triple $\{u,v,w\}$ is an edge of $\Xi_{0,1}^{(3)}(H)$ and corresponds to a fixed edge $f\in B$. By \eqref{eq:projected_edge_selected}, the probability that this edge is selected in $\widetilde{M}_{i_\star}$ is $3(1\pm 6C\epsilon)/pn$. Thus, the probability that there is an edge $\{u,v,w\}$ that is contained in the matching $M$ in this case is given by 
\begin{align*}
    (1\pm 3\epsilon)p\cdot\frac{3(1+6C\epsilon)}{pn}=\frac{3(1\pm 7C\epsilon)}{n}
\end{align*}
Hence, the total contribution of this case is given by
\begin{align}\label{eq:pair_case2}
    \PP(\cE)\cdot \PP(\cF)\cdot\frac{3(1+7C\epsilon)}{n}=\frac{4n}{9(n-1)}\cdot (1\pm n^{-\sqrt{C}})\cdot \frac{3(1+7C\epsilon)}{n}= \frac{4(1\pm10C\epsilon)}{3n}
\end{align}

Combining \eqref{eq:pair_case1} and \eqref{eq:pair_case2}, we obtain
\begin{align*}
    \PP\left(\exists e\in M:\{x,y\}\subseteq e\right)
    =2(1\pm 10C\epsilon)/n,
\end{align*}
which concludes the proof of \ref{itemM2}.

It remains to verify \ref{itemM3}. Let $A\subseteq H$ be a fixed set of edges of $H$. Note that if $A$ is not a matching, then the desired probability is $0$. Thus, we may assume that $A$ is a matching. For the algorithm to output a matching $M$ such that $A\subseteq M$, first the random bipartite structure $\Xi_{0,1}$ must be compatible with $A$. In particular, after the random equipartition is chosen, the matching $M_{0,1}$ must agree with the induced matching of size $a$ between the vertices of $A$ lying in $X_{0,1}$ and $Y_{0,1}$. Not every equipartition contains such an induced matching, but conditioned on any equipartition for which such an induced matching exists, this happens with probability exactly $(n/3-a)!/(n/3)!$. Thus, the probability of this compatibility event is at most $(n/3-a)!/(n/3)!$.

Now suppose that $A$ is compatible with the bipartite structure $\Xi_{0,1}$. The set of edges $A$ projects to a fixed matching $F_A$ of $B$. In order for the output matching $M$ to contain $F_A$, the matching $F_A$ must be contained in one matching $M'$ of the collection $\widetilde{\cM}$ produced by Lemma \ref{lem:randomizedalgorithm}, and then $M'$ must be selected as the uniformly chosen matching $\widetilde{M}_{i_\star}$. By property \ref{itemR4} of Lemma \ref{lem:randomizedalgorithm} and \eqref{eq:number_matchings_hypergraph_algorithm}, we obtain
\begin{align*}
    \PP(A\subseteq M)
    &\leq \frac{(n/3-a)!}{(n/3)!}\cdot
    \left(\frac{C\log (n/3)}{pn/3}\right)^{a-1}
    \cdot \frac{1}{t}\leq \left(\frac{9C\log n}{pn^2}\right)^{a},
\end{align*}
where we use the bound $m!/(m-a)!=\binom{m}{a}\cdot a!\geq (m/e)^a$ for $m\geq 1$. This proves \ref{itemM3} and concludes the proof of the lemma.
\end{proof}

\section{Stage I: Reducing the density}\label{sec:stageI}

In this section, we study the first stage of our iterative absorption process. As discussed in the overview, given a current remainder graph $R$ with concentrated codegree, the goal of this stage is to reduce its density by introducing a new source of randomness coming from a random graph $H$ of density $p$. This reduction is going to be done by choosing suitable bipartite projections (see Lemma \ref{lem:decomposition_bipartite_structures}) and then applying the randomized algorithm from the previous section to obtain a collection of almost perfect matchings covering most of the edges from $H\cup R$.

We would like to make a few remarks before delving into the details. First, our bipartification process (see Lemma \ref{lem:bipartification}) requires that the number of vertices from our graph is a multiple of~$3$. Unfortunately, this is not always necessarily true. To fix that we introduce a set of vertices, which we call $S$, that are going to serve as potential vertices missing our almost perfect matching. We will see later that this set $S$ is going to be defined as the set of vertices of small degree in the original random hypergraph. A downside of this new introduction is that we will need to analyze how well distributed the almost perfect matchings are with the vertices of $S$, which will make our statements slightly more technical. Second, the whole process described in this stage is not enough to guarantee that the new remainder graph $\widetilde{R}$ has concentrated codegree. Indeed, we can only impose that $\widetilde{R}$ has bounded codegree. To fix this issue we will need to perform a second stage in our iterative absorption process that we will handle in the next section. The following is the main result of this section, which we stated for a deterministic pseudorandom hypergraph~$H$.

\begin{lemma}[Stage I]\label{lem:global_covering}
Let $\epsilon,\delta, \lambda,p,r\in[0,1]$ be real numbers and let $n$ be a sufficiently large integer such that
\[
    n^{-1/5}\leq \epsilon\leq \delta\leq C^{-3},\qquad
    C\frac{\log n}{\epsilon^2 n}\leq r \leq p,\qquad
    C\epsilon^2\log n\leq \delta^6, \quad \text{and} \quad  \ Cn^{-1}\leq \lambda
\]
Let $S\subseteq[n]$ be a set with $|S|=n/100$. Let $H$ and $R$ be edge-disjoint $3$-graphs on $[n]$ such that $H$ is $(\epsilon,p,n)$-pseudorandom and $\deg_R(u,v)=(1\pm\delta)rn$ for every pair $\{u,v\}\subseteq[n]$. Moreover, suppose that the degrees of both $H$ and $R$ satisfy 
\begin{align*}
    \deg_H(u)=\frac{pn^2}{2}\pm \lambda n^2 \qquad \text{and}\qquad \deg_R(u)=\frac{rn^2}{2}\pm \lambda n^2
\end{align*}
for every vertex $u\in [n]$. Then there exists a collection of edge disjoint matchings $\cK$ in $H\cup R$ such that the following properties hold:
\begin{enumerate}[label=(F\arabic*)]
    \item\label{global1} The collection $\cK$ covers $R$, i.e., $R\subseteq \bigcup_{K\in\cK}K\subseteq H\cup R$.

    \item\label{global2} For every $K\in\cK$, we have $[n]\setminus V(K)\subseteq S$.

    \item\label{global3} For every $u\in S$, the number of matchings $K\in \cK$ missing the vertex $u$ is at most $Cpn$.

    \item\label{global4} Let $\widetilde{R}:=H\setminus\bigcup_{K\in\cK}K$. Then there exists a positive constant $c_{\star}$ bounded both above and below by a function of $C$ such that
    \begin{align*}
        \deg_{\widetilde{R}}(u)=c_{\star}\delta pn^2 \pm 5\lambda n^2 \qquad\text{and} \qquad \Delta_2(\widetilde{R})\leq C^3\delta pn,
    \end{align*}
    for $u\in [n]$
\end{enumerate}
Moreover, the same statement holds if we take $R=\emptyset.$
\end{lemma}

We start by proving a very useful decomposition into bipartite structures. We point out that a very similar decomposition was already introduced in \cites{ferber20181, ferber2017packing}.

\begin{lemma}\label{lem:decomposition_bipartite_structures}
Let $\eta\in[0,1]$ and let $n$ be an integer such that $\eta\geq C/n$ for sufficiently large $C$. Let $m\in 3\ZZ^+$ be an integer multiple of $3$ such that $n-m\in\{0,1,2\}$ and let $S\subseteq [n]$ be a set with $|S|=n/100$. Let $t\in\ZZ^+$ be a positive integer satisfying $t\geq C\eta^{-2}\log n$, let $\ell=tn$ and set $q:=3/2t$. Then there exists a randomized construction of bipartite structures $\Xi_i=(X_i,Y_i,Z_i,M_i,E_i)$, for $1\leq i\leq \ell$, such that each $\Xi_i$ is distributed as a random $(\eta,q,m)$-bipartite structure and with probability at least $1-n^{-\sqrt{C}}$, the following properties hold:
\begin{enumerate}[label=(G\arabic*)]
    \item\label{bipdec1} The sets $E_1,\ldots,E_\ell$ form a partition of $[n]^{(3)}$, i.e., $[n]^{(3)}=E_1\cup\cdots\cup E_\ell$.

    \item\label{bipdec2} If $V_i:=X_i\cup Y_i\cup Z_i$, then $|V_i|=m$ and $[n]\setminus V_i\subseteq S$ for every $1\leq i\leq \ell$.

    \item\label{bipdec3} For every $u\in S$, the number of bipartite structures $\Xi_i$ such that $u\notin V_i$ is at most $Ct/2$.

    \item\label{bipdec4} For every pair $\{u,v\}\subseteq[n]$, the number of bipartite structures $\Xi_i$ such that $\{u,v\} \in M_i$ is $(1\pm\eta)2\ell/3m$
\end{enumerate}
\end{lemma}

\begin{proof}
Let $\xi:=n-m\in\{0,1,2\}$. We start by constructing auxiliary bipartite structures $\widehat{\Xi}_i=(X_i,Y_i,Z_i,M_i,\widehat{E}_i)$, for $1\leq i \leq \ell$,
as follows. For each $1\leq i\leq \ell$, choose a set $S_i\subseteq S$ uniformly at random among all subsets of $S$ of size $\xi$, and let $V_i:=[n]\setminus S_i$. Then choose an equipartition $V_i=X_i\cup Y_i\cup Z_i$ uniformly at random, and choose a perfect matching $M_i\subseteq X_i\times Y_i$ uniformly at random. Finally, set $\widehat{E}_i:=M_i\times Z_i$. In particular, each $\widehat{\Xi}_i$ is a random $(0,1,m)$-bipartite structure on the vertex set $V_i$.

For every triple $e\in [n]^{(3)}$, let
\begin{align*}
    J_e:=\{i\in[\ell]:e\in \widehat{E}_i\}
\end{align*}
be the set of indices $i\in [\ell]$ of bipartite structures $\widehat{\Xi}_i$ containing the edge $e$ in the set $\widehat{E}_i$. Let $\cE$ be the event that for every triple $e\in [n]^{(3)}$ it holds that $|J_e|=(1\pm \eta/3)2t/3$. We will first show that event $\cE$ holds with high probability. Fix a triple $e\in[n]^{(3)}$. Since $\xi\leq 2$, $|S|=n/100$ and $\eta\geq C/n$, we have
\begin{align*}
    \PP(e\subseteq V_i)
    =\frac{\binom{|S|-|S\cap e|}{\xi}}{\binom{|S|}{\xi}}
    =1\pm \frac{\eta}{20},
\end{align*}
for $1\leq i \leq \ell$. Conditioned on the event $e\subseteq V_i$, the probability that $e\in M_i\times Z_i$ is
\begin{align*}
    \frac{2m^2}{9(m-1)(m-2)}\cdot \frac{3}{m}
    =\frac{2m}{3(m-1)(m-2)}
    =\left(1\pm \frac{\eta}{20}\right)\frac{2}{3n},
\end{align*}
where in the first equality we used that the probability that the three vertices of $e$ lie in three different parts of $X_i\cup Y_i\cup Z_i$ is $\frac{2m^2}{9(m-1)(m-2)}$, and then the two vertices lying in $X_i\cup Y_i$ are matched by $M_i$ with probability $3/m$. Therefore, since $\ell=tn$, we obtain
\begin{align}\label{eq:expectation-Ne}
    \E(|J_e|)=\ell\cdot \PP(e\in \widehat{E}_i)=\left(1\pm \frac{3\eta}{20}\right)\frac{2t}{3}.
\end{align}
Moreover, the random variable $|J_e|$ follows a binomial distribution, and thus by a Chernoff bound (Theorem \ref{chernoff}) and our choice of $t$, it holds that
\begin{align*}
    \PP\left(\big||J_e|-\E(|J_e|)\big|>\frac{\eta}{10}\E(|J_e|)\right)
    \leq 2\exp\left(-\frac{\eta^2\E(|J_e|)}{300}\right)
    \leq 2\exp\left(-\frac{\eta^2t}{500}\right)
    \leq n^{-C/600}.
\end{align*}
Taking a union bound over all triples $e\in[n]^{(3)}$, we obtain that $\PP(\cE)\geq 1-n^{-C/1000}$.

We now define the final edge sets $E_1,\ldots,E_\ell$. We split the construction into two cases. If the event $\cE$ occurs, then for every triple $e\in[n]^{(3)}$ we choose independently and uniformly at random an index $j(e)\in J_e$ and we put $e$ into $E_{j(e)}$. Equivalently,
\begin{align*}
    E_i:=\{e\in[n]^{(3)}:j(e)=i\}.
\end{align*}
for every $1\leq i \leq \ell$. If $\cE$ fails, then for each $1\leq i \leq\ell$ we instead choose $E_i\subseteq \widehat{E}_i$ as a random homogeneous subhypergraph of $\widehat{E}_i$ with edge probability $q$. Finally, let $\Xi_i=(X_i,Y_i,Z_i,M_i,E_i)$ be the desired bipartite structures.

We first note that the construction implies that each $\Xi_i$ is distributed as a random $(\eta,q,m)$-bipartite structure. Indeed, conditioned on the choices of the auxiliary structures $\widehat{\Xi}_1,\ldots,\widehat{\Xi}_\ell$, the choices of the edges of $E_i$ are independent. Moreover, on the event $\cE$, every edge $e\in \widetilde E_i$ is chosen with probability 
\begin{align}\label{eq:inverse-Ne-revised}
    \frac{1}{|J_e|}=(1\pm\eta)\frac{3}{2t}=(1\pm \eta)q,
\end{align}
while in the backup construction every edge is chosen with probability $q$.

We now verify properties \ref{bipdec1}--\ref{bipdec4} in order. We first note that Property \ref{bipdec1} holds if we condition on the event $\cE$. Indeed, in this case, by construction, each triple $e\in[n]^{(3)}$ is assigned to exactly one index $J_e$ and the sets $E_1,\ldots,E_\ell$ form a partition of $[n]^{(3)}$. Property \ref{bipdec2} follows directly from the definition of $V_i$. For every $1\leq i\leq \ell$, we have $|V_i|=m$ and $[n]\setminus V_i=S_i\subseteq S$.

We now verify property \ref{bipdec3}. Fix $u\in S$ and let
\begin{align*}
    L_u:=|\{i\in[\ell]:u\notin V_i\}|=|\{i\in[\ell]:u\in S_i\}|
\end{align*}
be the number of bipartite structures missing the vertex $u$. Note that since the sets $S_i$ are chosen uniformly, we have $L_u\sim \Bin(\ell,\xi/|S|)$. Because $\xi\leq 2$ and $|S|=n/100$, we have $\E(L_u)\leq 200t$. Moreover,
\begin{align*}
    \PP(L_u> Ct/2)
    &\leq \binom{\ell}{Ct/2}\left(\frac{\xi}{|S|}\right)^{Ct/2}
    \leq \left(\frac{2e\ell \xi}{Ct|S|}\right)^{Ct/2}
    \leq \left(\frac{400e}{C}\right)^{Ct/2}
    \leq n^{-C},
\end{align*}
where in the last inequality we used that $t\geq C\log n$ and that $C$ is sufficiently large. Thus, by a union bound over all $u\in S$, with probability at least $1-n^{-C/2}$, property \ref{bipdec3} holds for every $u\in S$.

Finally, we verify property \ref{bipdec4}. Fix a pair $\{u,v\}\subseteq[n]$ and let
\begin{align*}
    N_{u,v}:=|\{i\in[\ell]:\{u,v\}\in M_i\}|
\end{align*}
be the number of bipartite structures $\Xi_i$ such that the matching $M_i$ contains $\{u,v\}$. We estimate the expectation of $N_{u,v}$. Since $\xi\leq 2$ and $\eta\geq C/n$, we have
\begin{align*}
    \PP(\{u,v\}\subseteq V_i)=\frac{\binom{|S|-|S\cap \{u,v\}|}{\xi}}{\binom{|S|}{\xi}}
    =1\pm \frac{\eta}{20}.
\end{align*}
Conditioned on the event $\{u,v\}\subseteq V_i$, by symmetry every pair of vertices in $V_i$ has the same probability of being an edge of $M_i$. Since $M_i$ has $m/3$ edges, we obtain
\begin{align*}
    \PP(\{u,v\}\in M_i\mid \{u,v\}\subseteq V_i)
    =\frac{m/3}{\binom{m}{2}}
    =\frac{2}{3(m-1)}
    =\left(1\pm \frac{\eta}{20}\right)\frac{2}{3n}.
\end{align*}
Consequently, it holds that $\E(N_{u,v})=2(1\pm 3\eta/20)/3t$. Since $N_{u,v}$ follows a binomial distribution, by  Chernoff bound (Theorem \ref{chernoff}) it holds that
\begin{align*}
    \PP\left(|N_{u,v}-\E(N_{u,v})|>\frac{\eta}{10}\E(N_{u,v})\right)
    \leq 2\exp\left(-\frac{\eta^2\E(N_{u,v})}{300}\right)\leq 2\exp\left(-\frac{\eta^2t}{500}\right)
    \leq n^{-C/600}.
\end{align*}
Taking a union bound over all pairs $\{u,v\}\subseteq[n]$, we obtain that with probability at least $1-n^{-C/1000}$, property \ref{bipdec4} holds for every pair $\{u,v\}\subseteq[n]$.

Putting together the estimates for properties \ref{bipdec1}, \ref{bipdec3} and \ref{bipdec4}, and recalling that property \ref{bipdec2} holds deterministically, we conclude by a union bound that with probability at least $1-n^{-\sqrt{C}}$ all properties \ref{bipdec1}--\ref{bipdec4} hold simultaneously. This concludes the proof of the lemma.
\end{proof}

We now turn our attention to the proof of Lemma \ref{lem:global_covering}. 

\begin{proof}[Proof of Lemma \ref{lem:global_covering}]
    We prove the result by a randomized construction. Let $m \in 3\ZZ^+$ be the unique integer multiple of $3$ such that $n-m\in \{0,1,2\}$. Let $t, \ell\in \ZZ^+$ and $q\in[0,1]$ be such that 
    \begin{align}\label{eq:parametersII}
        t:= C\delta^{-2}\log n, \qquad  \ell:=tn \qquad \text{and} \qquad q:=3/2t.
    \end{align} 
    We first apply Lemma \ref{lem:decomposition_bipartite_structures} with $\eta=\delta$ to obtain, for $1\leq i\leq \ell$, that $\Xi_i=(X_i,Y_i,Z_i,M_i,E_i)$ is distributed as a random $(\delta,q,m)$-bipartite structure. Let $\cE_{\textrm{dec}}$ be the event that properties \ref{bipdec1}--\ref{bipdec4} hold for these structures. By Lemma \ref{lem:decomposition_bipartite_structures}, we have that $\PP(\cE_{\textrm{dec}})\geq 1-n^{-\sqrt{C}}$.

     For each $1\leq i\leq \ell$, let $B_i^{(2)}:=\Xi_i^{(2)}(H)$ and $R_i^{(2)}:=\Xi_i^{(2)}(R)$ be the bipartite projections of $H$ and $R$ on these bipartite structures. We now verify that, with high probability, these bipartite graphs satisfy the hypothesis of Lemma \ref{lem:randomizedalgorithm}. For every $1\leq i\leq \ell$, let $V_i:=X_i\cup Y_i\cup Z_i$ and $S_i:=[n]\setminus V_i\subseteq S$. Since $|S_i|=n-m\leq 2$, Proposition \ref{prop:deletion_pseudo} implies that the induced hypergraph $H[V_i]$ is $(\epsilon+C/pn,p,m)$-pseudorandom and $R[V_i]$ is such that $\deg_{R[V_i]}(u,v)=(1\pm(\delta+C/pn))rm$ for every $\{u,v\} \subseteq [n]$. By our hypothesis on the parameters, we have $\delta\geq \epsilon>C/pn$ and consequently, $H[V_i]$ is $(2\delta,p,m)$-pseudorandom and $\deg_{R[V_i]}(u,v)=(1\pm 2\delta)rn$ for every pair $\{u,v\}\subseteq V_i$.

    Since each $\Xi_i$ is distributed as a random $(\delta,q,m)$-bipartite structure, we should be able to apply Lemma \ref{lem:bipartification} for $H[V_i]$ and $R[V_i]$. Indeed, note that by~\eqref{eq:parametersII} and the hypothesis on $p$ and $C\epsilon^2 \log n\leq \delta^6$, we have
    \begin{align}\label{eq:condition_pq}
        p\geq C\frac{\log n}{\epsilon^2 n}\geq C\frac{\log n}{(2\delta)^2 n} \quad \text{and}\quad pq=\frac{3p}{2t}\geq \frac{\delta^2}{\epsilon^2 n}\geq C\frac{\log n}{(2\delta)^4 n}.
    \end{align}
    Similar calculations give that $r\geq C\frac{\log n}{(2\delta)^2 n}$ and $rq\geq C\frac{\log n}{(2\delta)^4 n}$. Thus, Lemma~\ref{lem:bipartification} implies that for $1\leq i \leq \ell$, with probability at least $1-m^{-C^{1/3}}$, the bipartite graphs $B_i^{(2)}$ and $R_i^{(2)}$ satisfy
    \begin{align*}
        \text{$B_i^{(2)}$ is $(200\delta,pq,m/3)$-pseudorandom} \qquad \text{and} \qquad \text{$R_i^{(2)}$ is $(1\pm 200\delta)rqm/3$-regular.}
    \end{align*}
    We remark that the same event holds for $R_i^{(2)}$ if $R=\emptyset$. In this case, we do not need to apply Lemma~\ref{lem:bipartification} to $R$. The fact that $R$ is the empty set gives immediately that $R_i^{(2)}=\emptyset$ and the condition holds for $r=0$. The remainder of the proof works whether $R$ is dense or empty without changing a sentence. Let $\cE_{\textrm{con}}$ be the event that the above conclusion holds for $1\leq i \leq \ell$. Then a union bound shows that $\PP(\cE_{\textrm{con}})\geq 1-n^{-C^{1/4}}$. Note that conditioned on the event $\cE_{\textrm{con}}$, the graphs $B_i^{(2)}$ and $R_i^{(2)}$ satisfy the hypothesis of \ref{lem:randomizedalgorithm}. Indeed, since $0\leq r/p\leq 1$, we have
    \begin{align*}
        (1\pm 200\delta)rqm/3=(r/p\pm200\delta)pqm/3.
    \end{align*}
    Hence, $R_i^{(2)}$ is $(\gamma\pm 200\delta)pqm/3$ regular for $\gamma=r/p$. Moreover, by~\eqref{eq:condition_pq} the density $pq$ satisfies the required condition for sufficiently large $C$. Therefore, we can apply Lemma \ref{lem:randomizedalgorithm} to $B_i^{(2)}$ and $R_i^{(2)}$ with density $pq$, error parameter $200\delta$ and $\gamma=r/p$.

    For each $1\leq i\leq \ell$, apply Lemma \ref{lem:randomizedalgorithm} to $B_i^{(2)}$ and $R_i^{(2)}$. If the algorithm succeeds, let $\cK_i^{(2)}$ be the collection of perfect matchings of $B_i^{(2)}\cup R_i^{(2)}$ obtained as the output, and let $\cK_i^{(3)}$ be the collection of $3$-uniform matchings obtained by lifting the matchings in $\cK_i^{(2)}$ through the projection map. Let
    \begin{align*}
        \cK:=\bigcup_{i=1}^{\ell}\cK_i^{(3)}.
    \end{align*}
    Let $\cF$ be the event that the algorithm succeeds for $1\leq i \leq \ell$. By Lemma \ref{lem:randomizedalgorithm}, the probability that one of the $\ell$ applications of Algorithm \ref{alg:randomized} fails is at most $n^{-C^{1/5}}$. Therefore, by a union bound, it holds that $\PP(\cF)\geq 1-n^{-C^{1/5}/2}$. From now on we assume that the events $\cE_\textrm{dec}, \cE_{\textrm{con}}$ and $\cF$ hold.

    We first verify \ref{global1}. By property \ref{itemR1} of Lemma \ref{lem:randomizedalgorithm}, the matchings in $\cK_i^{(2)}$ cover all edges of $R_i^{(2)}$. Lifting back to the $3$-graph, this implies that the matchings in $\cK_i^{(3)}$ cover all edges of $\Xi_i^{(3)}(R)=R\cap E_i$. Since, by property \ref{bipdec1}, the sets $E_1,\ldots,E_{\ell}$ form a partition of $[n]^{(3)}$, we obtain that $R\subseteq \bigcup_{K\in\cK}K$. On the other hand, each matching in $\cK_i$ is contained in $\Xi_i^{(3)}(H\cup R)\subseteq H\cup R$. Therefore, $R\subseteq \bigcup_{K\in\cK}K\subseteq H\cup R$, which proves \ref{global1}.

    We now verify \ref{global2}. Note that every matching in $\cK_i^{(2)}$ is a perfect matching on $\Xi^{(2)}(H\cup R)$. Therefore, its lift is a $3$-uniform matching covering exactly the set of vertices $V_i=X_i\cup Y_i\cup Z_i$. By property \ref{bipdec2}, we have $[n]\setminus V_i\subseteq S$. Hence, every matching $K\in\cK$ satisfies $[n]\setminus V(K)\subseteq S$, proving \ref{global2}.

    To prove \ref{global3}, first note that if $u\in V_i$, then every matching in $\cK_i^{(3)}$ covers the vertex $u$. Thus, a matching in $\cK_i^{(3)}$ can miss $u$ only if $u\notin V_i$. We now bound the number of matchings in each collection $\cK_i^{(3)}$. Let $\widetilde{R}_i^{(2)}:=B_i^{(2)}\setminus \bigcup_{K\in\cK_i^{(2)}}K$ be the bipartite remainder from the $i$-th application of Lemma \ref{lem:randomizedalgorithm}. Since every perfect matching in $\cK_i^{(2)}$ covers every vertex of $B_i^{(2)}\cup R_i^{(2)}$, for any vertex $x\in V(B_i^{(2)})$ we have
    \begin{align*}
        |\cK_i^{(3)}|=\deg_{B_i^{(2)}}(x)+\deg_{R_i^{(2)}}(x)-\deg_{\widetilde{R}_i^{(2)}}(x)\leq \deg_{B_i^{(2)}}(x)+\deg_{R_i^{(2)}}(x)\leq pqm,
    \end{align*}
    where we used $r<p$ in the last inequality. By property \ref{bipdec3}, for every $u\in S$ there are at most $Ct/2$ indices $i$ for which $u\notin V_i$. Combining these two observations with the fact that $m=n-O(1)$ and $q=3/2t$, we obtain that, for $u\in S$, the number of matchings $K\in \cK$ missing the vertex $u$ is at most 
    \begin{align*}
        \sum_{i:\: u\notin V_i}|\cK_i|\leq \frac{Ct}{2} \cdot pqm\leq Cpn.
    \end{align*}
    This concludes the proof of \ref{global3}.
    
    It remains to verify \ref{global4}. For each $1\leq i \leq \ell$, let $\widetilde{R}_i=\Xi^{(3)}_i(H)\setminus \bigcup_{K\in \cK_i} K$ be the remainder graph after we remove the lifted matchings $\cK_i$ from the corresponding part $\Xi^{(3)}_i(H)$ of $H$ and let
    \begin{align*}
        \widetilde{R}:=H\setminus\bigcup_{K\in\cK}K=\bigcup_{i=1}^\ell \widetilde{R}_i
    \end{align*}
    be the remainder hypergraph obtained after our process. We start by proving the codegree bound. Fix a pair $\{u,v\}\subseteq[n]$. We split the edges of $\widetilde{R}_i$ containing $\{u,v\}$ into two types, depending on whether $\{u,v\}$ is an edge of the auxiliary matching $M_i$ in the bipartite structure $\Xi_i$.

    First, consider the contribution when $\{u,v\}\in M_i$. In this case, every leftover edge of $\Xi_i^{(3)}(H)$ containing $\{u,v\}$ corresponds to an edge of the bipartite remainder $\widetilde{R}_i^{(2)}:=B_i^{(2)}\setminus \bigcup_{K^{(2)}\in \cK_i^{(2)}}K^{(2)}$ incident to the vertex $\{u,v\}\in M_i$. By property \ref{itemR2} of Lemma \ref{lem:randomizedalgorithm}, the degree of any vertex $x\in V(\widetilde{R}_i^{(2)})$ is 
    \begin{align*}
        \deg_{\widetilde{R}_i^{(2)}}(x)\leq 2400C\delta pqm/3.
    \end{align*}
    Moreover, by property \ref{bipdec4}, the number of indices $i$ with $\{u,v\}\in M_i$ is at most $Ct$. Therefore, by (\ref{eq:parametersII}), the total contribution of this type is at most
    \begin{align*}
        Ct\cdot 800C\delta pqm\leq C^3\delta pn/2.
    \end{align*}

    It remains to consider the indices $i$ such that $\{u,v\}\notin M_i$. For each such $i$, there is at most one edge in $E_i$ containing the pair $\{u,v\}$. Indeed, if such an edge exists, then one of $u,v$ lies in $Z_i$ and the other lies in a matched pair of $M_i$, which uniquely determines the third vertex of the edge. Since the sets $E_i$ form a partition of $[n]^{(3)}$, the total number of such possible edges over all $i$ is at most $\deg_H(u,v)\leq 2pn$.

    Let $W$ be the random variable counting the number of edges containing $\{u,v\}$ of the second type. Note that for each $1\leq i \leq \ell$, the event that there exists an edge $e\in E_i$ containing $\{u,v\}$ depends on a unique application of Lemma \ref{lem:randomizedalgorithm}, which is independent of all the other applications. Moreover, by property \ref{itemR3} of Lemma \ref{lem:randomizedalgorithm}, each of these possible edges remains uncovered with probability at most $1000C\delta$. Hence, the number of leftover edges of this type is a sum of independent Bernoulli variables and its mean is at most $1000C\delta\cdot  2 pn\leq C^2\delta pn$ . Therefore, by the Chernoff bound (Theorem \ref{chernoff}),
    \begin{align*}
        \PP\left(W>C^3\delta pn/2\right)
        \leq \exp(-C^3\delta pn/2)\leq n^{-C}.
    \end{align*}
    Let $\cG$ be the event that the above inequality holds for every pair $\{u,v\}\subseteq [n]$. A union bound shows that $\PP(\cG)\geq 1-n^{-\sqrt{C}}$. Consequently, if $\cG$ holds, we have by summing the contribution of both types that
    \begin{align*}
        \Delta_2(\widetilde{R})\leq C^3\delta pn/2+C^3\delta pn/2=C^3\delta pn.
    \end{align*}

    We now prove the degree estimate in \ref{global4}. We first estimate the average degree of $\widetilde{R}$.  Note that Property \ref{itemR2} of Lemma \ref{lem:randomizedalgorithm} gives us that $200C\delta pqm/3\leq \deg_{\widetilde{R}_i^{(2)}}(x)\leq 2400C\delta pqm/3$
    for every $1\leq i \leq \ell$ and vertex $x\in V(\widetilde{R}_i^{(2)})$. Hence, by lifting the edges of $\widetilde{R}^{(2)}_i$ to the $3$-graph $H$ and summing over all $1\leq i\leq \ell$, we obtain by (\ref{eq:parametersII}) that there exist positive constants $C_1, C_2>0$ such that
    \begin{align*}
        C_1\delta pn^3 \leq e(\widetilde{R})\leq C_2\delta pn^3.
    \end{align*}
    Let $c_{\star}:=3e(\widetilde{R})/{\delta pn^3}$.
    Then by the previous computation we have that $c_\star>C_2>0$ is a positive constant. Moreover, the average degree of $\widetilde{R}$ is exactly $c_{\star}\delta pn^2$.

    It remains to show that all vertex degrees concentrate. Let $s:=|\cK|$ and, for $u\in[n]$, let
    \begin{align*}
        a_u:=|\{K\in\cK:u\notin V(K)\}|
    \end{align*}
    be the number of matchings missing the vertex $u$. Since every matching in $\cK$ that contains $u$ covers exactly one edge incident to $u$, and since $\cK$ covers all the edges of $R$ (by \ref{global1}), we have that
    \begin{align}\label{eq:global_degree_identity}
        \deg_{\widetilde{R}}(u)=\deg_{H\cup R}(u)-s+a_u.
    \end{align}
    Now, by the hypothesis in the statement we have that $\deg_H(u)=pn^2/2\pm \lambda n^2$ and $\deg_R(u)=rn^2/2 \pm \lambda n^2$ for every $u\in [n]$. Hence, we have that for every pair of vertices $u,v\in[n]$,
    \begin{align}\label{eq:degree_gap_HR}
        |\deg_{H\cup R}(u)-\deg_{H\cup R}(v)|\leq  4\lambda n^2.
    \end{align}
    Moreover, by \ref{global3}, we have that $a_u=0$ if $u\notin S$ and $a_u\leq Cpn$ for $u\in S$. Combining this with \eqref{eq:global_degree_identity} and \eqref{eq:degree_gap_HR}, we obtain
    \begin{align*}
        |\deg_{\widetilde{R}}(u)-\deg_{\widetilde{R}}(v)|
        \leq 4\lambda n^2+Cpn\leq 5\lambda n^2.
    \end{align*}
    Therefore, every degree of $\widetilde{R}$ differs from the average degree by at most $5\lambda n^2$. This proves \ref{global4} and concludes the proof of the lemma.
\end{proof}

\section{Stage II: Concentrating the degree}\label{sec:stageII}

We now focus on the second stage of our absorption process. The goal here as described before is different than the one in Section \ref{sec:stageI}. In the last section, we use an auxiliary pseudorandom $3$-graph to in some sense `reduce' the density of our remainder set of edges. However, by doing that we obtain a new remainder that only has bounded codegree. In order to reiterate Lemma \ref{lem:global_covering}, we would need the codegrees to be concentrated. This is the goal of Lemma \ref{lem:stageII}: At the cost of slightly increasing the density of our current remainder, by using a random sprinkling of a pseudorandom $3$-graph we are able to obtain a codegree concentration that can be reused later. Again, because our technical tools require that the number of vertices from our graph is a multiple of $3$, we will need to introduce a set $S$ of potential vertices missing our almost perfect matchings. The following is the statement of the main result of this section.

\begin{lemma}[Stage II]\label{lem:stageII}
Let $\epsilon,\delta,p, q, r,\lambda\in[0,1]$ be real numbers, and let $k\in\ZZ^+$ be a positive integer satisfying
\begin{align*}
    \log n\leq k\leq n^{0.01},\qquad \epsilon\leq\frac{1}{k^2}, \qquad \delta\geq \frac{C^3}{k^2}, \qquad p\geq 0.01,\\ \frac{k^7\log n}{n} \leq r\leq \frac{1}{k^6}, \qquad q\geq C^4kr\log^3n, \qquad C\frac{\log n}{n}\leq \lambda\leq \frac{r}{3k^2}, 
\end{align*}
for a sufficiently large constant $C$ and integer $n$. Let $S\subseteq[n]$ be a set with $|S|=n/100$. Let $H$ and $R$ be edge-disjoint $3$-graphs on $[n]$ such that $H$ is $(\epsilon,p,n)$-pseudorandom and $\Delta_2(R)\leq Crn$. Moreover, suppose that both $H$ and $R$ satisfy the degree condition
\begin{align*}
    \deg_H(u)=\frac{pn^2}{2}\pm\lambda n^2 \qquad \text{and} \qquad \deg_R(u)=\frac{rn^2}{2}\pm \lambda n^2
\end{align*}
for every $u\in[n]$. Then, with probability at least $1-n^{-\sqrt{C}}$, the random homogeneous $3$-graph $H_q\subseteq H$ satisfies the following. There exists a collection of edge disjoint matchings $\cK$ in $H_q\cup R$ with the properties:
\begin{enumerate}[label=(H\arabic*)]
    \item\label{stageII1} The collection $\cK$ covers $R$, i.e., $R\subseteq \bigcup_{K\in\cK}K\subseteq H_q\cup R.$

    \item\label{stageII2} For every $K\in\cK$, we have $[n]\setminus V(K)\subseteq S$.

    \item\label{stageII3} For every $u\in S$, the number of matchings $K\in\cK$ missing the vertex $u$ is at most $Ckrn$.

    \item\label{stageII4} Let $\widetilde{R}:=H_q\setminus\bigcup_{K\in\cK}K$. Then for every pair $\{u,v\}\subseteq[n]$, it holds
    \begin{align*}
        \deg_{\widetilde{R}}(u)=\frac{\widetilde{r}n^2}{2}\pm 4\lambda n^2 \qquad \text{and} \qquad \deg_{\widetilde{R}}(u,v)=(1\pm \delta)\widetilde{r}n,
    \end{align*}
    where $\widetilde{r}:=pq+r-2|\cK|/n^2>0$. Moreover, it holds that $pq-2rk<\widetilde{r}<pq$.
\end{enumerate}
\end{lemma}

The proof of Lemma \ref{lem:stageII} is quite technical and because of that we will first give a description of the main steps and later analyze them in more detail. Let $H$ and $R$ be $3$-graphs satisfying the hypothesis of the lemma. Let $C$ be a constant larger than all the previous constants in the preceding Lemmas and Propositions. Let $\epsilon, \delta, p, q, r, \lambda \in [0,1]$ be parameters satisfying
\begin{align}\label{eq:parametersIII}
   \log n\leq k\leq n^{0.01},\qquad \epsilon\leq\frac{1}{k^2}, \qquad \delta\geq \frac{C^3}{k^2}, \qquad p\geq 0.01, \\ \frac{k^7\log n}{n} \leq r\leq \frac{1}{k^6}, \qquad q\geq C^4kr\log^3n, \quad \text{and} \quad C\frac{\log n}{n}\leq \lambda\leq \frac{r}{3k^2}. \nonumber
\end{align}
We begin by showing that there exists a partition of $R=E_1\cup \ldots \cup E_s$ into $s=(1+o(1))\Delta(R)=(1+o(1))rn^2/2$ matchings with the property that at most $s/k$ of those matchings satisfy $|E_i|<n/3-n/k$ (see Proposition \ref{prop:partition_remainder_into_matchings}). This is a consequence of an averaging argument together with the most recent bounds on the chromatic index theorem for hypergraphs with bounded degrees \cites{kang2023new, molloy2000near}. The result allows us to split the matchings into two types: For every $1\leq i \leq s$, we say that $E_i$ is \emph{full} if $|E_i|\geq n/3-n/k$, otherwise if $|E_i|<n/3-n/k$ we say that $E_i$ is \emph{meager}.

For every $1\leq i \leq s$, we partition each $E_i=\bigcup_{j=1}^k E_{i,j}$ into $k$ smaller matchings randomly where each $E_{i,j}$ receives an edge independently with probability $1/k$. Set $\xi=n\pmod 3$. For every $1\leq j\leq k$, let $S_{i,j}\subseteq S\setminus V(E_{i,j})$ be a set chosen uniformly at random of size $\xi$. We define the set of vertices
\begin{align*}
    U_{i,j}=[n]\setminus (V(E_{i,j})\cup S_{i,j}).
\end{align*}
Note by our choice of $\xi$ that $|U_{i,j}|$ is a multiple of $3$. The next step of the proof (see Proposition~\ref{prop:random_split_preserves_pseudorandomness}) shows that with high probability the induced hypergraphs $H[U_{i,j}]$ are all pseudorandom. In particular, if $E_i$ is full, then a consequence of Proposition \ref{prop:random_split_preserves_pseudorandomness} is that $H[U_{i,j}]$ is pseudorandom with very good error parameters. This is a consequence of standard concentration inequalities. Thus, after this step, every matching $E_{i,j}$ has a large pseudorandom reservoir $H[U_{i,j}]$ in which it can be completed to an almost perfect matching. For technical reasons, we also keep track at this point, for every vertex $u\in S$, in how many of the $U_{i,j}$ this vertex is missing.

The third step of the proof is to run a randomized algorithm (see Lemma \ref{lem:randomized_completion}) on $H$ that outputs with high probability a hypergraph $F=\bigcup_{i,j}F_{i,j}$ such that
\begin{enumerate}
    \item $F_{i,j}\subseteq H[U_{i,j}]$ is a matching.
    \item For every $\{u,v\}$, it holds that $\deg_F(u,v)=(1\pm \delta)(k-1)rn$. 
\end{enumerate}
Note that the output graph $F$ almost achieves our goal. Indeed, by taking for every $\{i,j\}\in [s]\times[k]$ the matching $K_{i,j}=E_{i,j}\cup F_{i,j}$ and setting $\cK=\{K_{i,j}\}$, we obtain a family of almost perfect matchings that covers $R$. Moreover, since $F$ has concentrated codegree, the graph $H\setminus F$ also has concentrated codegree. This is the most technical part of the proof and essentially uses all the machinery developed in the previous sections. We would like to emphasize that the error $(1\pm \delta)$ can only be achieved by the fact that most of the matchings $E_{i,j}$ originate from a full matching $E_i$ and that the corresponding $3$-graphs $H[U_{i,j}]$ are pseudorandom with very good error parameters. 

Even though Lemma \ref{lem:randomized_completion} outputs the desired graph $F$ complementing the matchings $E_{i,j}$, it has the caveat that it requires the density of $H$ to be nearly constant. This creates a problem, since the potential remainder $H\setminus F$ has density much larger than the previous remainder. The fourth step in the proof is to transfer the above deterministic completion inside $H$ to a completion inside a much sparser random subgraph $H_q$. For this, we use that the randomized completion algorithm in Lemma \ref{lem:randomized_completion} gives a probability distribution on the family $\mathcal{F}$ of all possible successful graphs $F$ satisfying a spread condition. In particular, it is going to be shown that the family $\mathcal F$ is $(C^2kr\log n)$-spread for a sufficiently large constant $C$. Then, by an application of the Kahn--Kalai conjecture, or more precisely by its fractional-threshold form, the random subgraph $H_q$ contains some $F\in\mathcal \cF$ with high probability provided that $q\geq C^4kr\log^3 n$. The proof is then completed in the last step by taking the remainder $\widetilde{R}=H_q\setminus F$. We now proceed to explain every step in more detail.

\subsection{Step I: Partitioning $R$ into matchings}

As described in the outline, we start the proof by decomposing $R$ into $(1+o(1))\Delta(R)$ matchings. Before we prove such a decomposition, we recall the following result on the chromatic index of bounded codegree hypergraphs by Kang, K\"{u}hn, Methuku and Osthus \cite{kang2023new}.

\begin{theorem}[\cite{kang2023new}, Corollary 1.6]\label{thm:kangnibble}
Let $D_1, D_2, N>0$, and let $0<\gamma,\mu<1$, and let $0<\eta<1/96$. Then there exists
$N_0=N_0(\gamma,\eta,\mu)$ such that the following holds for all $N\geq N_0$
and $D_1\geq \exp(\log^\mu N)$. Let $H$ be a $3$-graph on $N$ vertices with maximum degree at most $D_1$
and codegree at most $D_2\leq D_1^{1-\gamma}$. Then the chromatic index of $H$ is at most $D_1+D_1\left(D_1/D_2\right)^{-1/3-\eta}$.
\end{theorem}

A simple averaging argument combined with the last theorem gives the following.

\begin{proposition}\label{prop:partition_remainder_into_matchings}
There exists a partition $R=E_1\cup\cdots\cup E_s$ into matchings such that $s=\left(1\pm1/k^2\right)rn^2/2$. Moreover, all but at most $s/k$ of the matchings $E_i$ satisfy $|E_i|\geq n/3-n/k$.
\end{proposition}

\begin{proof}
Recall that the hypothesis on $R$ gives us that
\begin{align*}
    \Delta(R)\leq rn^2/2+\lambda n^2 \qquad \text{and} \qquad \Delta_2(R)\leq Crn.
\end{align*}
Note that by our hypothesis on $r$ and $\lambda$, the $3$-graph $R$ satisfies the hypothesis of Theorem \ref{thm:kangnibble}. Thus, the theorem gives a partition $R=E_1\cup\cdots\cup E_s$
into $s$ matchings such that
\begin{align}\label{eq:number-colors-bound}
    s\leq
    \Delta(R)+\Delta(R)\left(\frac{\Delta(R)}{\Delta_2(R)}\right)^{-1/3}\leq \left(1+\frac{2\lambda}{r}+Cn^{-1/3}\right)\frac{rn^2}{2}\leq (1+1/k^2)\frac{rn^2}{2},
\end{align}
where we use that $\lambda\leq r/3k^2$, $ k^3\leq n^{0.03}<n^{-1/3}$. The lower bound on $s$ comes from the fact that 
\begin{align}\label{eq:lower_s}
    s\geq \Delta(R)\geq \left(1-\frac{2\lambda}{r}\right)\frac{rn^2}{2}\geq (1-1/k^2)\frac{rn^2}{2}.
\end{align}
It remains to show that almost all of the matchings are large. Let
\begin{align*}
    I:=\left\{i\in[s]: |E_i|<\frac{n}{3}-\frac{n}{k}\right\}
\end{align*}
be the set of indices corresponding to matchings that are small. Then, by the hypothesis on $R$, it holds that
\begin{align*}
    |I|\cdot \left(\frac{n}{3}-\frac{n}{k}\right)+ (s-|I|)\cdot \frac{n}{3}\geq e(R)\geq \frac{\delta(R)\cdot n}{3}\geq \frac{rn^3}{6}-\frac{\lambda n^3}{3}.
\end{align*}
This together with (\ref{eq:number-colors-bound}), (\ref{eq:lower_s}), $\lambda \leq r/3k^2$ and $k\geq \log n$ implies that 
\begin{align*}
     |I|\leq \frac{k}{3}\cdot \left(s-\frac{rn^2}{2}+\lambda n^2\right)\leq \frac{s}{k}.
\end{align*}
This concludes the proof of the proposition.
\end{proof}

Let $R=E_1\cup \ldots \cup E_s$ be as given in the proposition above. We split the partition of matchings into the two types already mentioned: We say that an index $i\in [s]$ is \emph{full} if $|E_i|\geq n/3-n/k$. Otherwise, if $|E_i|<n/3- n/k$, we say that $i\in [s]$ is \emph{meager}. Proposition \ref{prop:partition_remainder_into_matchings} implies that there are at most $s/k$ meager indices.

\subsection{Step II: Further refining the matchings}

We now further split the collection of matchings $\{E_i\}_{1\leq i \leq s}$. For each $1\leq i \leq s$, we partition $E_i=\bigcup_{j=1}^k E_{i,j}$ into $k$ matchings at random such that for every $e\in E_i$ we have
\begin{align*}
    \PP(e\in E_{i,j})=\frac{1}{k}.
\end{align*}
Let $\xi \in \{0,1,2\}$ be such that $\xi=n \pmod 3$. For every $1\leq j\leq k$, let
\begin{align*}
    S_{i,j}\subseteq S\setminus V(E_{i,j})
\end{align*}
be a random set chosen uniformly at random of size $\xi$. Moreover, let 
\begin{align*}
    U_{i,j}=[n]\setminus (V(E_{i,j})\cup S_{i,j})
\end{align*}
be the set of vertices excluding the vertices from the matchings $E_{i,j}$ and the vertices in the random set $S_{i,j}$. The next lemma records some properties of the sets $E_{i,j}, U_{i,j}$ and $S_{i,j}$. 

\begin{proposition}\label{prop:random_split_preserves_pseudorandomness}
Let $E_{i,j}$, $U_{i,j}$ and $S_{i,j}$ be defined as above and let the parameters satisfy (\ref{eq:parametersIII}).
Then, with probability at least $1-n^{- C^{1/3}}$, the following properties hold:
\begin{enumerate}[label=(I\arabic*)]
    \item\label{split5} If the index $i$ is full, then the sets $U_{i,j}$ satisfy for $1\leq j\leq k$
    \begin{align*}
        |U_{i,j}|=\left(1-\frac{1}{k}\pm \frac{10}{k^2}\right)n.
    \end{align*}

    \item\label{split1} If the index $i$ is full, then the induced $3$-graph $H[U_{i,j}]$ is $(C/k^2, p, |U_{i,j}|)$-pseudorandom for every $1\leq j\leq k$.

    \item\label{split2} If the index $i$ is meager, then the induced $3$-graph $H[U_{i,j}]$ is $(C/k, p, |U_{i,j}|)$-pseudorandom for every $1\leq j \leq k$.

    \item\label{split3} For every vertex $u\in S$, the number of $S_{i,j}$ containing $u$ is at most $Csk/n$.
    
    \item\label{split4} For every pair $\{u,v\}\subseteq[n]$, the number of $U_{i,j}$ containing the pair $\{u,v\}$ is
    \[
        I_{u,v}:=|\{(i,j)\in[s]\times[k]:\{u,v\}\subseteq U_{i,j}\}|
        =
        \left(1\pm \frac{10}{k^2}\right)(sk-2s).
    \]
\end{enumerate}
\end{proposition}

\begin{proof}
First, for each $1\leq i \leq s$, we reveal the random partitions $E_i=\bigcup_{j=1}^k E_{i,j}$. For every $(i,j)\in [s]\times [k]$, let 
\begin{align*}
W_{i}:=[n]\setminus V(E_i) \qquad \text{and} \qquad W_{i,j}:=[n]\setminus V(E_{i,j}).
\end{align*} 
We will show that, with probability at least $1-n^{-\sqrt C}$, the following properties hold simultaneously:
\begin{enumerate}[label=(\roman*)]
    \item\label{item_spliti} $|E_{i,j}|\leq n/k$ for every $(i,j)\in[s]\times [k]$;
    \item\label{item_splitii} If the index $i$ is full, then $H[W_{i,j}]$ is $(300/k^2,p,|W_{i,j}|)$-pseudorandom
    for every $j\in[k]$;
\end{enumerate}
After that, we will reveal the sets $S_{i,j}$ and verify the properties \ref{split1}--\ref{split4}.

We start with \ref{item_spliti}. For fixed $i$ and $j$, the random variable $|E_{i,j}|$ follows a binomial distribution $\Bin(|E_i|, 1/k)$. Since $|E_i|\leq n/3$ and $k\leq n^{0.01}$, we have
\begin{align*}
    \PP\left(|E_{i,j}|>\frac{n}{k}\right)\leq \binom{|E_i|}{n/k}\left(\frac{1}{k}\right)^{n/k}\leq \left(\frac{e}{3}\right)^{n/k}\leq \exp(-\sqrt{n}).
\end{align*}
A union bound over all pairs $(i,j)\in[s]\times[k]$ gives that \ref{item_spliti} holds with probability at least $1-n^{-C}$.

We now prove \ref{item_splitii}. Fix a pair $(i,j)\in [s]\times [k]$ with $i$ being full. We first estimate the size of the set $W_{i,j}$ in this case. Note that since $i$ is full, we have $|E_i|\geq n/3-n/k\geq n/6$. Thus, because $|E_{i,j}|$ follows $\Bin(|E_i|, 1/k)$ with mean $|E_i|/k$, by Chernoff's bound (Theorem \ref{chernoff}), we obtain
\begin{align*}
    \PP\left(\big||E_{i,j}|-\E|E_{i,j}|\big|>\frac{n}{k^2}\right)\leq 2\exp\left(-\frac{n^2}{3k^3|E_i|}\right)\leq 2\exp(-\sqrt{n}),
\end{align*}
where we use the assumption that $k\leq n^{0.01}$. Let $\cE_1$ be the event that for every index $i$ that is full, we have that $|E_{i,j}|=|E_i|/k \pm n/k^2$. Then, by taking a union bound over all possible pairs $(i,j)$, the calculation above shows that $\PP(\cE_1)\geq 1-n^{-2C}$. In particular, conditioned on $\cE_1$, note that
\begin{align}\label{eq:size-U-prime}
|W_{i,j}|=n-3|E_{i,j}|=n-\frac{3|E_i|}{k}\pm \frac{3n}{k^2}=\left(1-\frac{1}{k}\pm \frac{6}{k^2}\right)n
\end{align}
since $|E_i|=n/3\pm n/k$.

We now try to control the codegrees over $W_{i,j}$. Fix a pair $\{u,v\}\subseteq W_{i,j}$, and let $N_H(u,v):=\{w\in[n]:\{u,v,w\}\in H\}$ be the neighborhood of $\{u,v\}$ in $H$. For every edge $e\in E_i$, we define the intersection of this edge with the neighborhood of $\{u,v\}$ by
\begin{align*}
    a_e:=|e\cap N_H(u,v)|.
\end{align*}
Note that $0\leq a_e\leq 3$ for every $e\in E_i$. Let $X_{u,v}^{i,j}$ be the random variable counting the number of neighbors of $\{u,v\}$ that are in $V(E_{i,j})$, i.e.,
\begin{align*}
    X_{u,v}^{i,j}:=|N_H(u,v)\cap V(E_{i,j})|
      =\sum_{e\in E_i}a_e\1_{\{e\in E_{i,j}\}}.
\end{align*}
The mean of $X_{u,v}^{i,j}$ is given by $\E(X_{u,v}^{i,j})=|N_H(u,v)\cap V(E_i)|/k$. Moreover, since $X_{u,v}^{i,j}$ is a sum of independent random variables, each on the interval $[0,3]$, then by Hoeffding's inequality (Theorem \ref{thm:Hoeffding}), we obtain that
\begin{align*}
    \PP\left(\left|X-\E X\right|>\frac{pn}{k^2}\right)
    \leq
    2\exp\left(-\frac{p^2n^2}{9k^4|E_i|}\right)
    \leq \exp(-\sqrt{n}),
\end{align*}
where we used the assumption $p\geq 0.01$, $|E_i|\leq n/3$ and $k\leq n^{0.01}$. Let $\cE_2$ be the event that for every pair $(i,j)\in [s]\times [k]$ and every pair $\{u,v\}\subseteq [n]$ we have
\begin{align}\label{eq:codegree-random-split-good}
    |N_H(x,y)\cap V(E_{i,j})|
    =
    \frac{|N_H(x,y)\cap V(E_i)|}{k}\pm \frac{pn}{k^2}.
\end{align}
By taking a union bound, we have that $\PP(\cE_2)\geq 1-n^{-2C}$.

We now condition on the events $\cE_1$ and $\cE_2$ and verify the pseudorandom properties \ref{Q1:codegree}--\ref{Q4:smallZ} for $H[W_{i,j}]$. We start with \ref{Q1:codegree}. Note by (\ref{eq:codegree-random-split-good}) that 
\begin{align*}
    \deg_{H[W_{i,j}]}(u,v)
    =\deg_H(u,v)-X_{u,v}^{i,j}=\deg_H(u,v)-\frac{|N_H(u,v)\cap V(E_i)|}{k}\pm \frac{pn}{k^2}
\end{align*}
Note that because $i$ is full, we have
\begin{align*}
    \deg_H(u,v)\geq |N_H(u,v)\cap V(E_i)|\geq \deg_H(u,v)-|W_i|\geq \deg_H(u,v)-\frac{n}{k}.
\end{align*}
Moreover, note that $\deg_H(u,v)=(1\pm \epsilon)pn$ since $H$ is $(\epsilon,p,n)$-pseudorandom. These two observations together with the assumptions $\epsilon\leq 1/k^2$, $p\geq 0.01$ and \eqref{eq:size-U-prime} imply that
\begin{align*}
    \deg_{H[W_{i,j}]}(u,v)
    =\left(1-\frac{1}{k}\pm \frac{200}{k^2}\right)pn= \left(1\pm \frac{300}{k^2}\right)p|W_{i,j}|.
\end{align*}
This concludes the proof of \ref{Q1:codegree}. Properties \ref{Q2:largesets}--\ref{Q4:smallZ} follow from the fact that they are hereditary, and from $\epsilon<300/k^2$ and Remark \ref{rmk:pseudo_hypergraph}. Since $\cE_1 \cap \cE_2$ holds with probability at least $1-n^{-C}$, this implies that with at least this probability the $3$-graph $H[W_{i,j}]$ is $(300/k^2,p,|W_{i,j}|)$-pseudorandom for every $i$ that is full. Therefore, with probability $1-n^{-\sqrt{C}}$ both \ref{item_spliti} and \ref{item_splitii} hold.

We now condition on the event that \ref{item_spliti} and \ref{item_splitii} hold, and reveal the sets $S_{i,j}$. We first verify \ref{split5} and \ref{split1}. Note that $U_{i,j}=W_{i,j}\setminus S_{i,j}$. Since each $S_{i,j}$ has size at most two, we have $|U_{i,j}|=|W_{i,j}|\pm 2$ and property \ref{split5} follows from (\ref{eq:size-U-prime}) and the fact that $k\leq n^{0.01}$. Moreover, Proposition \ref{prop:deletion_pseudo} implies that deleting $S_{i,j}$ increases the error parameter in the pseudorandomness condition by at most $8/pn\leq 1/k^2$,
where we use that $p\geq 0.01$ and $k\leq n^{0.01}$. Therefore, the $3$-graph $H[U_{i,j}]$ is $(C/k^2,p,|U_{i,j}|)$-pseudorandom. This proves \ref{split1}.

To verify \ref{split2} we also use Proposition \ref{prop:deletion_pseudo}. Fix a meager index $i$ and $j\in[k]$. By \ref{item_spliti}, the set $V(E_{i,j})$ has size at most $3n/k$. Therefore $U_{i,j}=[n]\setminus V(E_{i,j})\cup S_{i,j}$ is obtained from $[n]$ by deleting at most $3n/k+2\leq 4n/k$ vertices. Since $H$ is $(\epsilon,p,n)$-pseudorandom, by Proposition \ref{prop:deletion_pseudo}, the induced $3$-graph $H[U_{i,j}]$ is $(\epsilon+16/pk,p,|U_{i,j}|)$-pseudorandom.
The result now follows since $p\geq 0.01$ and $\epsilon\leq 1/k^2$.

Now we turn our attention to \ref{split3}. Fix $u\in S$. For each pair $(i,j)$, we have by \ref{item_spliti} and $k\geq \log n$ that
\begin{align*}
|S\setminus V(E_{i,j})|=|S|-|V(E_{i,j})|\geq \frac{n}{100}-\frac{3n}{k}\geq \frac{n}{200}.   
\end{align*}
Therefore, the probability that $u$ is selected in $S_{i,j}$ is at most
\[
    \PP(u\in S_{i,j})\leq \frac{2}{|S\setminus V(E_{i,j})|}
    \leq \frac{400}{n}.
\]
The choices of the sets $S_{i,j}$ are independent over all pairs $(i,j)$. Hence, we have that
\begin{align*}
    \PP\left(|\{(i,j):u\in S_{i,j}\}|>\frac{Csk}{n}\right)\leq \binom{sk}{\frac{Csk}{n}}\left(\frac{400}{n}\right)^{Csk/n}\leq \left(\frac{400e}{C}\right)^{Csk/n}\leq n^{-C},
\end{align*}
where we use the fact that $s\geq (1-1/k^2)rn^2/2$, $k\geq \log n$ and $r\geq k^7\log n/n$. A union bound over all $u\in S$ gives that \ref{split3} holds with probability $1-n^{-\sqrt{C}}$. Let $\cE_S$ be the event corresponding to it.

Finally, we analyze \ref{split4}. We first ignore the random sets $S_{i,j}$ and count the number of pairs $(i,j)$ such that $\{u,v\}\subseteq W_{i,j}$. For a fixed pair $\{u,v\}\subseteq[n]$, let $I_{u,v}^{\neq}$ be the set of indices $i\in[s]$ such that $u$ and $v$ are both covered by $E_i$, but they lie in two distinct edges of $E_i$. For each $i\in I_{u,v}^{\neq}$, let $Y_i^{u,v}$ be the indicator of the event that these two edges are assigned to the same part in the random partition of $E_i$. Set
\[
    Y^{u,v}:=\sum_{i\in I_{u,v}^{\neq}}Y_i^{u,v}.
\]
The random variables $Y_i^{u,v}$ are independent $0/1$-random variables and each has expectation $1/k$. Since $|I_{u,v}^{\neq}|\leq s$, we have that $\E(Y^{u,v})\leq s/k$.
Hence, it holds that
\begin{align*}
    \PP\left(Y^{u,v}>\frac{3s}{k}\right)\leq \binom{|I_{u,v}^{\neq}|}{3s/k}\left(\frac{1}{k}\right)^{3s/k}\leq \left(\frac{e\cdot |I_{u,v}^{\neq}|}{3s}\right)^{3s/k}\leq \exp(-\sqrt{n}),
\end{align*}
where we use the fact that $s\geq (1-1/k^2)rn^2/2$, $k\leq n^{0.01}$ and $r\geq k^7\log n/n$. Taking a union bound over all pairs $\{u,v\}\subseteq[n]$, with probability at least $1-n^{-C}$ we have that $Y^{u,v}\leq 3s/k$ for every pair $\{u,v\}\subseteq[n]$. Let $\cE_3$ be this event.

Assume that $\cE_3$ holds. Let
\begin{align*}
I_{u,v}^0:=\{(i,j)\in [s]\times [k]:\: \{u,v\} \subseteq W_{i,j}\}.
\end{align*}
We claim that
\begin{align}\label{eq:I0-count-pair}
    |I^0_{u,v}|
    =sk-\deg_R(u)-\deg_R(v)+\deg_R(u,v)+Y^{u,v}.
\end{align}
Indeed, note that $\{u,v\}\subseteq W_{i,j}$ is the same as $u, v\notin V(E_{i,j})$. For a fixed value of $i$, note that there is exactly one value of $j$ such that $u\in V(E_{i,j})$ if $u \in E_i$, and no values otherwise. Similarly, there is exactly one value of $j$ such that $v\in V(E_{i,j})$ if $v\in E_i$ and no values otherwise. If $u$ and $v$ lie in the same edge of $E_i$, then the same value of $j$ has been counted twice in the sum $\deg_R(u)+\deg_R(v)$, and we must subtract it back once, which is done by the term $\deg_R(u,v)$. If $u$ and $v$ lie in distinct edges of $E_i$ but these two edges are assigned to the same part, then again the same value of $j$ has been counted twice in $\deg_R(u)+\deg_R(v)$, and we must subtract it back once, which is done by the term $Y^{u,v}$. This proves the claim.

By the assumptions on $R$, Proposition \ref{prop:partition_remainder_into_matchings}, and the fact that $\lambda\leq r/3k^2$, we have
\begin{align*}
    \deg_R(x)=\left(1\pm \frac{2\lambda}{r}\right)\frac{rn^2}{2}=\left(1\pm \frac{3}{k^2}\right)s,
\end{align*}
for every $x\in [n]$. Moreover, since $k\leq n^{0.01}$, we also have $\deg_R(u,v)\leq Crn\leq 3Cs/n\leq s/k$. Putting everything together, we obtain from (\ref{eq:I0-count-pair}) and from $k\geq \log n$ that 
\begin{align}\label{eq:I0-final-count}
    |I^0_{u,v}|
    =sk-2s\pm \frac{5s}{k}
    =\left(1\pm \frac{5}{k^2}\right)(sk-2s)
\end{align}
for every pair $\{u,v\}\subseteq[n]$.

It remains to account for the random sets $S_{i,j}$. Assume that $\cE_3$ and $\cE_S$ hold. Recall that
\begin{align*}
    I_{u,v}:=\{(i,j)\in[s]\times[k]:\{u,v\}\subseteq U_{i,j}\}.
\end{align*}
The only pairs $(i,j)$ which belong to $I^0_{u,v}$ but not to $I_{u,v}$ are those for which $u\in S_{i,j}$ or $v\in S_{i,j}$. Therefore, if $Z_u:=|\{(i,j):\: u \in S_{i,j}\}|$, then
\begin{align*}
    |I^0_{u,v}|-Z_u-Z_v\leq |I_{u,v}|\leq |I^0_{u,v}|.
\end{align*}
By \ref{split3} and the fact that $k\leq n^{0.01}$, we have that $Z_u\leq Csk/n\leq s/k$. Hence, combining this with \eqref{eq:I0-final-count}, we obtain
\[
    |I_{u,v}|=\left(1\pm \frac{10}{k^2}\right)(sk-2s)
\]
for every pair $\{u,v\}\subseteq[n]$. This completes the proof of \ref{split4}. Since all the properties hold conditioned on the events $\cE_1,\cE_2,\cE_3$ and $\cE_S$, we obtain that all the properties hold with probability at least $1-n^{-C^{1/3}}$.
\end{proof}

\begin{remark}\label{rmk:p_large}
    We would like to remark that Proposition \ref{prop:random_split_preserves_pseudorandomness} is the only point in the proof where it is required that $p$ is nearly constant. In particular, this requirement is necessary to prove \ref{split1} (More specifically, the error in relation (\ref{eq:codegree-random-split-good})). Although a similar statement could be obtained with an appropriate $p=o(1)$, for simplicity we decided to just use the lower bound $p\geq 0.01$. 
\end{remark}

\subsection{Step III: A randomized completion algorithm}

By Proposition \ref{prop:random_split_preserves_pseudorandomness} with high probability, the sets $\{E_{i,j}\}$, $\{U_{i,j}\}$ and $\{S_{i,j}\}$ satisfy properties \ref{split5}--\ref{split4}. Hence, there exists a choice of such sets satisfying the properties. From now on, fix one of those choices. The next step in the proof is to show the existence of a certain type of graphs $F$ that have the property that they complete the matchings $E_{i,j}$ to almost perfect matchings. We start with a definition. Let $L$ be an integer satisfying $L=\sum_{i,j} \frac{|U_{i,j}|}{3}$. Note that this is indeed an integer, since by our construction, the quantity $|U_{i,j}|$ is always a multiple of $3$.

\begin{definition}\label{def:cF}
    We say that a $3$-graph $F$ is \emph{good} if $F$ has $L$ edges and there exists an edge disjoint partition $F=\bigcup_{i,j} F_{i,j}$ for $(i,j)\in [s]\times [k]$ satisfying the following properties:
    \begin{enumerate}[label=(J\arabic*)]
        \item\label{cF1} $F_{i,j}$ is a perfect matching with vertex set $U_{i,j}$
        \item\label{cF2} For every $\{u,v\} \in [n]$, it holds that 
        \begin{align*}
        \deg_{F}(u,v)=\left(1\pm \frac{C^3}{2k^2}\right)\left(\frac{2sk-2s}{n}\right).
        \end{align*}
    \end{enumerate}
    Let $\cF$ be the family of all good $3$-graphs $F$ with vertex set $[n]$
\end{definition}

The main result in this step is a randomized algorithm that with high probability outputs a graph $F\in \cF$ from $H$ with a spread probability distribution.

\begin{lemma}\label{lem:randomized_completion}
Let $H$ be an $(\epsilon,p,n)$-pseudorandom $3$-graph on $[n]$. Let $\{E_{i,j}\}$, $\{U_{i,j}\}$ and $\{S_{i,j}\}$ be as fixed above and let all parameters satisfy (\ref{eq:parametersIII}). Then there exists a randomized algorithm that, with probability at least
$1-n^{-C^{1/6}}$, outputs a good $3$-graph $F\in \cF$. Moreover, conditioned on the event that the algorithm succeeds, for every set of edges $A\subseteq H$ of size $a$, we have
\begin{align*}
\PP(A\subseteq F)\leq\left(C^2kr\log n\right)^a.
\end{align*}
\end{lemma}

Before we give a proof of Lemma \ref{lem:randomized_completion}, we first describe the randomized algorithm used in the lemma. For the sake of convenience, we order the pairs $(i,j) \in [s]\times[k]$ in lexicographical order and relabel them by $t:=(i,j)\in [s]\times [k]$. We will abuse the notation and often treat $t$ as an element of $[sk]$, instead of $[s]\times[k]$. This, in particular, allows us to talk about the step right before $t$ as $t-1$. 

At each step $t$, suppose that the algorithm has already constructed
the matchings $F_{\ell}$ for all pairs $\ell<t$, and let $F^{<t}:=\bigcup_{\ell<t}F_{\ell}$ be the $3$-graph built so far. The goal of the next step is to find a perfect matching $F_t$ in $U_t$. We consider the induced $3$-graph
\begin{align*}
    \widetilde{H}_{t}:=\left(H\setminus F^{<t}\right)[U_t].
\end{align*}
That is, $\widetilde{H}_t$ is the $3$-graph $H[U_t]$ minus the edges of the $3$-graph $F^{<t}$. If the first coordinate of $t$ is a full index, then by Proposition \ref{prop:random_split_preserves_pseudorandomness} the $3$-graph $H[U_t]$ is $(C/k^2,p,|U_t|)$-pseudorandom. Based on that, we check whether $\widetilde{H}_t$ is $(2C/k^2,p, |U_t|)$-pseudorandom. In case that $\widetilde{H}_t$ is not pseudorandom, the algorithm declares failure. Similarly, if the first coordinate of $t$ is a meager index, then we check whether $H_{i,j}$ is $(2C/k,p,|U_t|)$-pseudorandom. Again, in case that this does not hold, the algorithm declares failure.

After the algorithm passes the pseudorandom test, we apply Lemma
\ref{lem:randomized_hypergraph_matching} to $\widetilde{H}_t$ to obtain a perfect matching $F_{t}\subseteq \widetilde{H}_{t}\subseteq H[U_t]$ with vertex set $U_{t}$. In case that the application of Lemma \ref{lem:randomized_hypergraph_matching} fails, then again our algorithm declares failure. Otherwise, it proceeds to the next step. By the end of all steps, if the algorithm never declared failure, then it outputs a $3$-graph $F:=\bigcup_{t\in [s]\times [k]}F_{t}$. Here is a formal description of the procedure.

\begin{algorithm}[H]
\caption{Randomized completion algorithm}
\label{alg:randomized_completion}
\begin{algorithmic}[1]
\STATE Set $F:=\emptyset$.
\FOR{$t=1$ \textbf{to} $sk$}
        \STATE Set $\widetilde{H}_t:=(H\setminus F)[U_t]$.
        \IF{$\widetilde{H}_{t}$ is not appropriately pseudorandom}
            \STATE Declare \textbf{failure}.
        \ENDIF
        \STATE Apply Lemma \ref{lem:randomized_hypergraph_matching} to $\widetilde{H}_{t}$ to obtain a perfect matching $F_{t}$ of $\widetilde{H}_{t}$.
        \IF{Lemma \ref{lem:randomized_hypergraph_matching} fails}
            \STATE Declare \textbf{failure}.
        \ENDIF
        \STATE Let $F:=F\cup F_{i,j}$.
\ENDFOR
\STATE Return $F$.
\end{algorithmic}
\end{algorithm}

We first prove that the algorithm succeeds with high probability.

\begin{proposition}\label{claim:completion_succeeds}
Algorithm \ref{alg:randomized_completion} succeeds with probability at least $1-n^{-C^{1/6}}$.
\end{proposition}

\begin{proof}
There are two possibilities of failure in Algorithm \ref{alg:randomized_completion}: Either the current $3$-graph $\widetilde{H}_t$ fails the pseudorandom test or an application of Lemma \ref{lem:randomized_hypergraph_matching} fails. We begin with a deterministic observation. Let $\cE^{<t}$ be the event that at the beginning of some step $t\in [s]\times [k]$, the algorithm has not declared failure and we have $\Delta_2(F^{<t})\leq Csk/n$.
We claim that, conditioned on $\cE^{<t}$, the $3$-graph $\widetilde{H}_t$ never fails the pseudorandomness test. That is, if the first coordinate of $t$ is a full index, then $\widetilde{H}_t$ is $(2C/k^2,p,|U_t|)$-pseudorandom. Otherwise, if the first coordinate of $t$ is a meager index, then $\widetilde{H}_t$ is $(2C/k,p,|U_t|)$-pseudorandom.

We start by analyzing the full case first. Note by Proposition \ref{prop:random_split_preserves_pseudorandomness} that $H[U_t]$ is $(C/k^2,p,|U_t|)$-pseudorandom. We just need to show that the deletion of the edges of $F^{<t}$ will not destroy the pseudorandom property. For every pair $\{u,v\}\subseteq U_t$, the deletion of $F^{<t}$ changes the codegree by at most $\Delta_2(F^{<t})$. By using that $s\leq (1+1/k^2)rn^2/2$, $r\leq 1/k^6$, $p\geq 0.01$ and \ref{split5} we have
\begin{align*}
    \Delta_2(F^{<t})\leq \frac{Csk}{n}\leq \frac{C^2pn}{k^5}\leq\frac{C}{k^2}p|U_t|  
\end{align*}
Therefore,
\[
    \deg_{\widetilde{H}_t}(u,v)=
    \deg_{H[U_{t}]}(u,v)\pm \frac{C}{k^2}p|U_t|=\left(1\pm \frac{2C}{k^2}\right)p|U_t|.
\]
Hence, property \ref{Q1:codegree} holds for error $2C/k^2$.

We now verify the remaining pseudorandomness properties. For property \ref{Q2:largesets}, let $M\subseteq U_{t}^{(2)}$ be a matching and let $Z\subseteq U_{t}$. Since every edge of $F^{<t}$ containing an element of $M$ contains exactly one pair of $M$, we have
\begin{align*}
    e_{F^{<t}}(M,Z)
    \leq |M|\cdot \Delta_2(F^{<t})
    \leq \frac{Csk}{n}|M|.
\end{align*}
If $|M|,|Z|\geq C |U_{t}|/k^2$, then again by the fact $s\leq (1+1/k^2)rn^2/2$, $r\leq 1/k^6$, $k\geq \log n$, $p\geq 0.01$ and $\ref{split5}$ we have
\begin{align*}
    e_{F^{<t}}(M,Z)\leq  \frac{Csk}{n}|M|
    \leq \frac{C^2pn}{k^5}|M|\leq 
    \frac{C}{k^2}p|M||Z|.
\end{align*}
Thus property \ref{Q2:largesets} is preserved after deleting $F^{<t}$. Properties \ref{Q3:smallM} and \ref{Q4:smallZ} follow automatically from Remark \ref{rmk:pseudo_hypergraph} and the fact that $\widetilde{H}_t$ has fewer edges than $H[U_t]$. We remark that the case that the first coordinate is a meager index is analogous to the previous one. Indeed, in no place of the proof did we use that the index is full except for the target error parameter.  Since in the meager case  the allowed error parameter is $C/k$, which is larger than $C/k^2$, the proof follows the same lines. This proves the deterministic observation.

It remains to prove that the events $\cE^{<t}$ hold throughout the algorithm with high probability. Let $\tau$ be the first step where the algorithm has either declared failure, or that $\Delta_2(F^{<\tau})>Csk/n$. That is, either in the beginning of step $\tau$ we are not at event $\cE^{<\tau}$ and consequently, we cannot guarantee that we are passing the pseudorandomness test; or we declare failure in the application of Lemma \ref{lem:randomized_hypergraph_matching} at this step. In case that such a step does not exist, we set $\tau=\infty$. Fix a pair $\{u,v\}\subseteq[n]$. For each step $t<\tau$, let $X_t^{u,v}$ be the indicator of the event that the matching chosen at step $t$ has an edge containing the pair $\{u,v\}$. We set $X_t^{u,v}=0$ for every $t\geq \tau$

Note that if the algorithm has not declared failure until the step $t$, then $\deg_{F^{<t}}(u,v)=\sum_{\ell\leq t} X_{\ell}^{u,v}$. We will bound this quantity by using Freedman's inequality. Let $\{\cA_t\}_{t\in[s]\times [k]}$ be a filtration collecting the history of all probabilistic choices of the algorithm over time, i.e., $\cA_t$ is the history of the algorithm until the end of the step $t$. We claim that for every $t\in [s]\times [k]$, it holds that
\begin{align}\label{eq:prob_Xuv_I}
    \PP(X_t^{u,v}=1\mid \cA_{t-1})\leq \frac{3}{n}.
\end{align}
Indeed, if $\{u,v\}\not\subseteq U_t$ or if $t\geq \tau$, then the probability above is zero. Otherwise, the $3$-graph $\widetilde{H}_t$ is pseudorandom and the application of Lemma \ref{lem:randomized_hypergraph_matching} succeeded. We split into two cases. If the first coordinate of $t$ is a full index, then $\widetilde{H}_t$ is $(2C/k^2, p, |U_t|)$-pseudorandom. Hence, by property \ref{itemM2} from Lemma \ref{lem:randomized_hypergraph_matching}, property \ref{split5} of Proposition \ref{prop:random_split_preserves_pseudorandomness} and the assumptions from (\ref{eq:parametersIII}), we have that
\begin{align*}
\PP(X_t^{u,v}=1 \mid \cA_{t-1})=\left(1\pm\frac{20C^2}{k^2}\right)\frac{2}{|U_t|}\leq \frac{3}{n}.
\end{align*}
If the first coordinate of $t$ is a meager index, then a similar calculation holds and we obtain (\ref{eq:prob_Xuv_I}). 

Let $p_t^{u,v}=\E(X_t^{u,v}\mid \cA_{t-1})$ be the random variable measuring the expected value of $X_t^{u,v}$ given the history $\cA_{t-1}$. Note that by (\ref{eq:prob_Xuv_I}), the random variable always satisfies $p_{t}^{u,v}\leq 3/n$. Let $\{Y_t^{u,v}\}_{t\in [s]\times[k]}$ be the martingale defined by
\begin{align*}
    Y_t^{u,v}:=\sum_{\ell\leq t}\left(X_{\ell}^{u,v}-\E(X_\ell^{u,v}\mid \cA_{\ell-1})\right)=\sum_{\ell\leq t}(X_{\ell}^{u,v}-p_{\ell}^{u,v}).
\end{align*}
Then $\{Y_t^{u,v}\}_{t\in [s]\times[k]}$ is a martingale with increments bounded by $1$. Moreover, since $X_{\ell}^{u,v}$ is an indicator random variable and (\ref{eq:prob_Xuv_I}) holds, we have that the predictable quadratic variation is always bounded by
\begin{align*}
    W_t^{u,v}=\sum_{\ell\leq t}\E\left(\left(X_{\ell}^{u,v}-p_{\ell}^{u,v}\right)^2 \,\,\middle |\,\,\cA_{\ell-1}\right)\leq \sum_{\ell \leq t} p_{\ell}^{u,v}\leq \frac{3sk}{n}
\end{align*}

Therefore, by the two-sided Freedman's inequality (Theorem \ref{thm:freedman}), we obtain
\begin{align*}
    \PP\left(\exists t\geq 0 :\: |Y_t^{u,v}|>\frac{s}{kn}\right)
    \leq
    2\exp\left(-\frac{s^2/2k^2n^2}{3sk/n+s/3kn}\right)\leq 2\exp\left(-\frac{s}{8k^3n}\right)
    \leq n^{-C},
\end{align*}
where we use $r\geq k^7\log n/n$, $s=(1\pm 1/k^2)rn^2/2$, and $k\geq \log n$. Let $\cE_{\deg}$ be the event that for every $\{u,v\}\subseteq [n]$ and $t\in [s]\times[k]$, we have
\begin{align}\label{eq:cE_deg}
     |\sum_{\ell\leq t}X_\ell^{u,v}-\sum_{\ell\leq t}p_{\ell}^{u,v}|<\frac{s}{kn}.
\end{align}
Then a union bound together with the estimate given by Freedman's inequality gives that $\PP(\cE_{\deg})\geq 1-n^{-\sqrt{C}}$. Moreover, let $\cE_{\ref{lem:randomized_hypergraph_matching}}$ be the event that whenever the application of Lemma \ref{lem:randomized_hypergraph_matching} was called, the application was successful. We remark that $\cE_{\ref{lem:randomized_hypergraph_matching}}$ is not the same as the event that Algorithm \ref{alg:randomized_completion} has not declared failure, because we still consider outcomes in $\cE_{\ref{lem:randomized_hypergraph_matching}}$ where the algorithm fails due to the pseudorandom test. Then by a union bound over all possible numbers of steps and the fact that one application of Lemma \ref{lem:randomized_hypergraph_matching} succeeds with probability $1-n^{-C^{1/5}}$, we obtain that $\PP(\cE_{\ref{lem:randomized_hypergraph_matching}})=1-n^{-C^{1/5}/2}$. 

Now suppose that $\omega \in \cE_{\deg}\cap \cE_{\ref{lem:randomized_hypergraph_matching}}$ is an outcome in the intersection of both events. We claim that $\omega \in \cE^{<t}$ for every $t\in [s]\times[k]$, and in particular that the algorithm finishes without declaring failure. We proceed by induction on $t$. The base case holds trivially. Suppose that for every $\ell <t$ we have that $\omega\in \cE^{<\ell}$. This implies that the algorithm has not declared failure at the beginning of step $t-1$ and $\Delta(F^{<t-1})\leq Csk/n$. Thus, by our deterministic observation, on step $t-1$ the hypergraph $\tilde{H}_{t-1}$ passes the pseudorandom test. Since $\omega\in \cE_{\ref{lem:randomized_hypergraph_matching}}$, we also have that the application of Lemma \ref{lem:randomized_hypergraph_matching} succeeded. Therefore, the algorithm has not declared failure at the beginning of step $t$. It remains to show the codegree bound for $\Delta(F^{<t})$. For this, note that by (\ref{eq:prob_Xuv_I}) and (\ref{eq:cE_deg}) it holds that
\begin{align*}
    \deg_{F^{<t}}(u,v)=\sum_{\ell\leq t}X_{\ell}^{u,v}(\omega)<\sum_{\ell \leq t}p_{\ell}^{u,v}(\omega)+\frac{s}{nk}<\frac{Csk}{n}.
\end{align*}
This concludes the induction and the proof of the claim. Hence, the probability that Algorithm \ref{alg:randomized_completion} succeeds is at least $\PP(\cE_{\deg}\cap \cE_{\ref{lem:randomized_hypergraph_matching}})\geq 1-n^{-C^{1/6}}$, which proves the proposition.
\end{proof}

We now proceed to prove the lemma.

\begin{proof}[Proof of Lemma \ref{lem:randomized_completion}]
Let $\cE_{\deg}$ be the event defined by (\ref{eq:cE_deg}) and let $\cE_{\ref{lem:randomized_hypergraph_matching}}$ be the event that the application of Lemma \ref{lem:randomized_hypergraph_matching} never declared failure. We showed in Proposition \ref{claim:completion_succeeds} that the probability of being in the intersection of these two events is given by $\PP(\cE_{\deg}\cap \cE_{\ref{lem:randomized_hypergraph_matching}})\geq 1-n^{-C^{1/6}}$ and if $\omega \in \cE_{\deg} \cap \cE_{\ref{lem:randomized_hypergraph_matching}}$, then $\omega$ is an outcome where the algorithm never declared failure. We claim that for $\omega\in \cE_{\deg}\cap \cE_{\ref{lem:randomized_hypergraph_matching}}$, the $3$-graph $F$ obtained as the output of the algorithm is good, i.e., it satisfies properties \ref{cF1} and \ref{cF2}.

Since for each step $t\in [s]\times [k]$, the matching $F_{t}$ is a perfect matching on $U_{t}$, property \ref{cF1} follows immediately. Thus, the only property in Definition \ref{def:cF} that remains to verify is the codegree estimate in \ref{cF2}. Fix a pair $\{u,v\}\subseteq[n]$. By the fact that $\omega\in \cE_{\deg}$, that $\omega$ is an outcome in which the algorithm never declared failure ($\tau=\infty$) and (\ref{eq:cE_deg}), we have
\begin{align}\label{eq:deg_F_uv}
    \deg_F(u,v)=\sum_{t\in [s]\times[k]}X_{t}^{u,v}=\sum_{t\in [s]\times[k]}p_t^{u,v}\pm \frac{s}{kn}.
\end{align}
Recall that $I_{u,v}:=\{t\in[s]\times[k]:\{u,v\}\subseteq U_{t}\}$ is the number of $U_t$ containing the pair $\{u,v\}$. To estimate the sum $\sum_{t}p_t^{u,v}$, note that $p_{t}^{u,v}=0$ if and only if $t\notin I_{u,v}$. Moreover, if the first coordinate of $t$ is a full index, then by property \ref{itemM2} from Lemma \ref{lem:randomized_hypergraph_matching}, property \ref{split5} of Proposition \ref{prop:random_split_preserves_pseudorandomness}, the fact that $\widetilde{H}_t$ is $(2C/k^2,p, |U_t|)$-pseudorandom and the assumptions from~(\ref{eq:parametersIII}) 
\begin{align*}
    p_t^{u,v}(\omega)=\PP(X_t^{u,v}=1 \mid \cA_{t-1})=\left(1\pm\frac{20C^2}{k^2}\right)\frac{2}{|U_t|}=\left(1\pm \frac{30C^2}{k^2}\right)\frac{2k}{(k-1)n}.
\end{align*}
Similarly, if the first coordinate of $t$ is a meager index, then we have
\begin{align*}
p_t^{u,v}(\omega)=\PP(X_t^{u,v}=1 \mid \cA_{t-1})=\left(1\pm\frac{20C^2}{k}\right)\frac{2}{|U_t|}=\left(1\pm \frac{30C^2}{k}\right)\frac{2k}{(k-1)n}.
\end{align*}
Since, by Proposition \ref{prop:partition_remainder_into_matchings}, there are at most $s/k$ meager indices, we obtain that there are at most $s$ steps $t\in [s]\times [k]$ such that the first coordinate is a meager index. Let $I_{u,v}^{\text{mea}}$ be the set of such steps $t$. Hence, we can bound the sum of $p_t^{u,v}$ by
\begin{align*}
    \sum_{t\in [s]\times[k]} p_t^{u,v}&=\left(1\pm \frac{30C^2}{k^2}\right)\cdot \frac{2k}{(k-1)n}\cdot |I_{u,v}\setminus I_{u,v}^{\text{mea}}|+ \left(1\pm\frac{30C^2}{k}\right)\cdot\frac{2k}{(k-1)n}\cdot |I_{u,v}^{\text{mea}}|\\
    &=\left(|I_{u,v}|\pm\frac{30C^2}{k}s \pm \frac{30C^2}{k^2}|I_{u,v}|\right)\cdot \frac{2k}{(k-1)n}
\end{align*}
By property \ref{split4} of Proposition \ref{prop:random_split_preserves_pseudorandomness}, we obtain that $|I_{u,v}|=(1\pm 10/k^2)(sk-2s)$ and consequently
\begin{align*}
    \sum_{t\in [s]\times[k]} p_t^{u,v}=\left(1\pm\frac{C^3}{4k^2}\right)\frac{2sk-2s}{n},
\end{align*}
for sufficiently large $C$. Hence, property \ref{cF2} now follows from (\ref{eq:deg_F_uv}) and error parameter $C^3/2k^2$. This concludes the proof of the claim and the first part of Lemma \ref{lem:randomized_completion}, since $\PP(\cE_{\deg}\cap \cE_{\ref{lem:randomized_hypergraph_matching}})\geq 1- n^{-C^{1/6}}$.

It remains to prove the spread estimate. Let $\cE_{\text{com}}$ be the event that the algorithm succeeds and outputs a good graph $F\in \cF$. By what was just proved above, $\PP(\cE_{\text{com}})\geq 1-n^{-C^{1/6}}$. Let $A\subseteq H$ be a fixed set of edges with $|A|=a$. We first estimate the probability that both events $\{A\subseteq F\}$ and $\cE_{\text{com}}$ hold. For every function $\psi:A\rightarrow [s]\times[k]$, set $A_{t}:=\psi^{-1}(t)$ for $t\in [s]\times [k]$. Note that if $A\subseteq F$, then there exists a partition $A=\bigcup_{t\in[s]\times[k]}A_t$ such that $A_t\subseteq F_t$. Therefore, there is a function $\psi$ corresponding to such a partition.

Fix a function $\psi$ and a partition $A=\bigcup_{t}A_{t}$. Note that conditioned on the event that the algorithm has not declared failure until step $t$ and has passed the pseudorandom test, we have by property \ref{itemM3} of Lemma \ref{lem:randomized_hypergraph_matching}, by \ref{split5} and by the assumptions in \eqref{eq:parametersIII} that the probability that $A_{t}\subseteq F_t$ is bounded by
\begin{align*}
\PP(A_t\subseteq F_t\mid \text{has not declared failure})\leq \left(\frac{9C\log n}{p|U_t|^2}\right)^{|A_t|}\leq \left(\frac{10C\log n}{pn^2}\right)^{|A_t|}.
\end{align*}
We remark that this bound is independent of the error parameter of the pseudorandomness of $\widetilde{H}_t$ as long as the error parameter is sufficiently small, which is achieved by the fact that $k\geq \log n$ (independently of whether the first coordinate is full or meager). Thus, by a union bound over all functions $\psi:A\rightarrow[s]\times[k]$ and by the bound above, we have
\begin{align*}
    \PP(\{A\subseteq F\} \cap \cE_{\text{com}})
    \leq
    \sum_{\psi:A\rightarrow[s]\times[k]}
    \prod_{t\in [s]\times[k]}
    \left(\frac{10C\log n}{pn^2}\right)^{|A_t|} \leq 
    \left(\frac{10Csk\log n}{pn^2}\right)^a.
\end{align*}
Since $s=(1\pm 1/k^2)rn^2/2$, $p\geq 0.01$ and the fact that conditioning on $\cE_{\text{com}}$ only changes this upper bound by a factor of $(1-n^{-C^{1/6}})^{-1}\leq2$, we obtain
\begin{align*}
    \PP(A\subseteq F\mid \cE_{\text{com}})
    \leq
    \left(C^2kr\log n\right)^a,
\end{align*}
for sufficiently large $C$. This completes the proof of the lemma.
\end{proof}

\subsection{Step IV: Transferring to a sparse random subhypergraph}

As discussed in the overview of Section \ref{sec:stageII} and also in Remark~\ref{rmk:p_large}, in order for us to obtain a copy of $F\in \cF$, the density of the pseudorandom $3$-graph $H$ cannot be too small. This could be potentially a problem for our strategy, since we do not want to increase the density of the remainder graph by too much. To fix that, in this step, we transfer the result of Lemma~\ref{lem:randomized_completion} to sparse random subhypergraphs of $H$. The main ingredient is the fractional version of the Kahn-Kalai conjecture due to Frankston, Kahn, Narayanan and Park~\cite{frankston2021thresholds}. 

We first remind the reader of the notation and the statement of the expectation-threshold theorem. Let $X$ be a finite set and, for $p\in[0,1]$, let $\mu_p$ denote the
product probability measure on $2^X$ in which each element of $X$ is included independently with probability $p$, i.e., $\mu_p(B)=p^{|B|}(1-p)^{|X\setminus B|}$. Let $\cB\subseteq 2^X$ be a non-empty increasing family of $X$. Since $\cB$ is increasing, it follows that $\mu_p(\cB):=\sum_{B\in \cB} \mu_p(B)$ is strictly increasing in $p$. Let $p_c(\cB)$ denote the unique value of $p$ such that $\mu_p(\cB)=1/2$.

We say that $\cB$ is \emph{weakly $p$-small} if there exists a function $g:2^X\longrightarrow \mathbb{R}^+$ such that
\begin{align}\label{eq:weakly_p_small}
    \sum_{S\subseteq B} g(S)\geq 1, \quad \forall B\in \cB \qquad \text{and} \qquad \sum_{S\subseteq X} g(S)p^{|S|}\leq \frac{1}{2}.
\end{align}
The \emph{fractional expectation-threshold} of $\cB$ is defined as the maximum value of $p$ such that $\cB$ is weakly $p$-small,
\begin{align*}
    q_f(\cB)
    :=
    \max\{p\in[0,1]:\cB\text{ is weakly $p$-small}\}.
\end{align*}
Finally, let $\ell(\cB)$ be the size of a largest minimal element of $\cB$. The fractional Kahn-Kalai conjecture states that $p_c(\cB)$ is closely bounded by the fractional threshold $q_f(\cB)$.

\begin{theorem}[\cite{frankston2021thresholds}]
\label{thm:FKNP}
There exists a universal constant $K>0$ such that, for every finite set $X$ and every non-empty increasing family $\mathcal B\subseteq 2^X$,
\begin{align*}
    p_c(\cB)
    \leq
    K q_f(\cB)\log \ell(\cB).
\end{align*}
\end{theorem}

The main reason why Theorem~\ref{thm:FKNP} is useful for us comes from the following standard result in the literature, which we provide a proof of for completeness. We say that a probability measure $\nu:\cB\rightarrow [0,1]$ is \emph{$p$-spread} if for $p\in[0,1]$, it holds that a random $B\subseteq \cB$ chosen according to $\nu$ satisfies $\PP(S\subseteq B)\leq p^{|S|}$ for every $S\subseteq X$.

\begin{proposition}\label{prop:spread_KK}
    Let $p\in[0,1]$. Let $\cB \subseteq 2^X$ be a non-empty increasing family and $\nu: \cB\rightarrow [0,1]$ be a $p$-spread measure. Then $q_f(\cB)\leq p$. 
\end{proposition}

\begin{proof}
    Let $g:2^X\rightarrow \RR_{\geq 0}$ be any covering function satisfying $\sum_{S\subseteq B}g(S)\geq 1$ for every $B\in \cB$. Let $B\in\cB$ be a random set chosen according to the probability measure $\nu$. Since $\nu$ is $p$-spread, we have
    \begin{align*}
     1\leq \E\left(\sum_{S\subseteq B} g(S)\right)=\sum_{S\subseteq X} g(S)\PP(S\subseteq B)\leq \sum_{S\subseteq X} g(S)p^{|S|}.
    \end{align*}
    This contradicts condition~\eqref{eq:weakly_p_small}. Therefore, the family $\cB$ is not weakly $p$-small and consequently $q_f(\cB)\leq p$.
\end{proof}

We are now able to use both Theorem~\ref{thm:FKNP} and Proposition~\ref{prop:spread_KK} to prove our sparse transferring lemma.

\begin{lemma}\label{lem:sparse_transfer}
    Let $H$ be an $(\epsilon,p,n)$-pseudorandom $3$-graph on $H$. Let $\{E_{i.j}\}$, $\{U_{i,j}\}$ and $\{S_{i,j}\}$ be as chosen in Step III, and let all parameters satisfy~\eqref{eq:parametersIII}. Let $\cF$ be the family of graphs defined in Definition~\ref{def:cF}. Then, with probability at least $1-n^{-C}$, the random homogeneous subgraph $H_q\subseteq H$ contains a good $3$-graph $F\in \cF$.
\end{lemma}

\begin{proof}
    Let $\cF_H:=\{F\in\cF:F\subseteq H\}$ be the family of good graphs $F\in \cF$ that are subgraphs of $H$ and let $\cB\subseteq 2^{H}$ be the increasing family generated by
    $\cF_H$, i.e.,
    \begin{align*}
        \cB:=\{G\subseteq H:\: \text{there exists }F\in\cF_H
        \text{ such that }F\subseteq G\}.
    \end{align*}
    By Lemma~\ref{lem:randomized_completion}, the family $\cF_H$ is
    non-empty. Moreover, conditioned on the event that the randomized
    algorithm in Lemma~\ref{lem:randomized_completion} succeeds, its
    output induces a $\rho$-spread probability measure $\nu$ on $\cF_H$, where $\rho:=C^2kr\log n$. We extend the measure to the whole set $\cB$, by assigning measure zero to every member of $\cB\setminus\cF_H$. It is easy to check that such an extension does not change the fact that $\nu$ is $\rho$-spread.

    In order to apply Theorem~\ref{thm:FKNP}, we first estimate the parameters $q_f(\cB)$ and $\ell(\cB)$. By our construction of $\nu$ and Proposition~\ref{prop:spread_KK}, we obtain that $q_f(B)\leq \rho$. To estimate $\ell(\cB)$, we note, by construction, that the minimal elements of $\cB$ are precisely the elements of $\cF_H$. Since each $F\in \cF$ has size $L=\sum_{i,j}|U_{i,j}|/3\leq n^3$, we have $\log \ell(B)\leq \log L\leq 3\log n$. Therefore, by Theorem~\ref{thm:FKNP}, it follows that
    \begin{align*}
        p_c(\cB)\leq Kq_f(\cB)\log\ell(\cB)\leq 3KC^2kr\log^2 n.
    \end{align*}
    for some universal constant $K>0$. Translating from the definition of $p_c(\cB)$, this means that for $q_0:=3KC^2kr\log^2 n$, it holds that the probability that $H_{q_0}$ does not contain an $F\in \cF$ is
    \begin{align}\label{eq:prob_12}
        \PP\left(\bigwedge_{F\in \cF}  F\not\subseteq H_{q_0}\right)=\frac{1}{2}.
    \end{align}

    The probability obtained in~\eqref{eq:prob_12} is still not large enough for our needs. To improve it, we use the following standard boosting technique. Let $b:=C\log n$ and let $H^{(1)},\ldots,H^{(b)}$ be independent copies of the random homogeneous subgraph $H_{q_0}\subseteq H$. By \eqref{eq:prob_12} and
    the fact that all copies are independent, we obtain that
    \begin{align}\label{eq:boost-transfer}
        \PP\left(\bigwedge_{j\in [b]}\bigwedge_{F\in \cF} F\not\subseteq H^{(j)}\right)= \prod_{j\in [b]}\PP\left(\bigwedge_{F\in \cF}F\not\subseteq H^{(j)}\right)=2^{-b}
        \leq n^{-C}.
    \end{align}
    A standard coupling argument shows that the union $\bigcup_{j=1}^bH^{(j)}$ is distributed as the random homogeneous subgraph $H_{q_\star}\subseteq H_q$, where
    \begin{align*}
        q_\star=1-(1-q_0)^b\leq bq_0\leq 3KC^3kr\log^3 n\leq C^4kr\log ^3n\le q
    \end{align*}
    holds for sufficiently large $C$. Therefore, we can couple $\bigcup_{j=1}^b H^{(j)}$ and $H_q$ such that $\bigcup_{j=1}^b H^{(j)}\subseteq H_q$. Hence, by~\eqref{eq:boost-transfer},
    \begin{align*}
        \PP\left(\bigwedge_{F\in \cF}F\not\subseteq H_q\right)\leq n^{-C},
    \end{align*}
    which concludes the proof of the proposition.
\end{proof}

\subsection{Step V: Putting it all together}

We now use the previous steps to finish the proof of Lemma~\ref{lem:stageII}.

\begin{proof}[Proof of Lemma~\ref{lem:stageII}]
    We start the proof by collecting the results obtained in the previous step. Let $R=E_1\cup\cdots\cup E_s$ be the partition of the original remainder graph $R$ obtained from Proposition~\ref{prop:partition_remainder_into_matchings}, where $s=(1\pm1/k^2)rn^2/2$. By Proposition~\ref{prop:random_split_preserves_pseudorandomness}, there exists a choice of the sets $\{E_{i,j}\}$, $\{U_{i,j}\}$ and
    $\{S_{i,j}\}$ satisfying properties~\ref{split5}--\ref{split4}.
    We fix such a choice for the remainder of the proof. Let $\cF$ be the family of graphs defined in Definition~\ref{def:cF}.

    Now we turn our attention to the choice of the random homogeneous subgraph $H_q\subseteq H$. Note by~\eqref{eq:parametersIII} that $p$, $q$, $\epsilon$ and $\delta$ satisfy that
    \begin{align*}
        pq\geq C^3kr\log ^3 n\geq \frac{C^3k^8\log^4 n}{n}\geq C\frac{\log n}{\delta n^4} \quad \text{and} \quad \epsilon\leq \frac{1}{k^2}\leq \frac{\delta}{8}.
    \end{align*}
    Thus, by Proposition~\ref{prop:random_hypergraph_pseudo}, we obtain from the fact that $H$ is $(\epsilon,p,n)$-pseudorandom that with probability at least $n^{1-\sqrt{C}}$ $H_q$ is $(\delta/2, pq, n)$-pseudorandom. Let $\cE_{\textrm{co}}$ be such an event. We also concentrate the degree of $H_q$. Note that by the hypothesis, the $3$-graph $H$ has degree $\deg_H(u)=pn^2/2\pm \lambda n^2$. Since $\deg_{H_q}(u)$ is distributed as the binomial $\Bin(\deg_H(u),q)$, we have that $\E(\deg_{H_q}(u))=pqn^2/2\pm \lambda qn^2$. Thus, by Theorem~\ref{chernoff} (Chernoff) and the fact that $\lambda\geq C\log n /n$, it holds that
    \begin{align*}
        \PP\left(\left|\deg_{H_q}(u)-\E(\deg_{H_q}(u))\right|>\lambda n^2\right)\leq 2\exp\left(-\frac{\lambda^2 n^4}{3(pqn^2+\lambda qn^2)}\right)\leq 2\exp\left(-\frac{\lambda^2n^2}{6}\right)\leq n^{-C}.
    \end{align*}
    Let $\cE_{\textrm{deg}}$ be the event that for every vertex $u \in [n]$, it holds that $\deg_{H_q}(u)=pqn^2/2\pm 2\lambda n^2$. Then, by a union bound and the inequality stated above, we have $\PP(\cE_{\textrm{deg}})\geq 1-n^{-C/2}$. Finally, let $\cE_{\textrm{good}}$ be the event that $H_q$ has a good $3$-graph $F\in \cF$ as a subgraph. By Lemma~\ref{lem:sparse_transfer} we obtain that $\PP(\cE_{\textrm{good}})\geq 1-n^{-C}$. For the rest of the proof we will fix an $H_q\subseteq H$ that is an outcome in $\cE_{\textrm{co}}\cap \cE_{\textrm{deg}}\cap \cE_{\textrm{good}}$. By the computations above, such an outcome happens with probability at least $1-n^{-1/3}$.

    Let $F=\bigcup_{(i,j)\in [s]\times[k]}F_{i,j} \subseteq H_q$ be the copy of a good $3$-graph inside of $H_q$. For every $(i,j)\in [s]\times[k]$, we define $K_{i,j}:=E_{i,j}\cup F_{i,j}$ to be the almost perfect matching at step $(i,j)$ and define
    \begin{align*}
        \cK:=\{K_{i,j}:(i,j)\in[s]\times[k]\}
    \end{align*}
    to be the family of all those matchings. Finally, let $\widetilde{R}=H_q\setminus \bigcup_{K\in \cK}K=H_q\setminus F$.

    We now verify properties~\ref{stageII1}--\ref{stageII4}. Condition~\ref{stageII1} follows immediately from the fact that the sets $E_{i,j}$ form a partition of $R$, while $\bigcup_{K\in \cK} K=R\cup F\subseteq R\cup H_q$. Condition~\ref{stageII2} follows from property~\ref{cF1} and the definition of $U_{i,j}$. Indeed, it holds for every $(i,j)\in [s]\times[k]$ that $[n]\setminus  V(K_{i,j})=S_{i,j}\subseteq S$.

    To prove~\ref{stageII3}, fix $u\in S$. By the discussion above, the matching $K_{i,j}$ misses a vertex $u$ if and only if
    $u\in S_{i,j}$. Therefore, by property~\ref{split3},
    \begin{align}\label{eq:a_u_final}
        a_u:=\left|\{K\in\cK:u\notin V(K)\}\right|=
        \left|\{(i,j)\in[s]\times[k]:u\in S_{i,j}\}\right|\leq \frac{Csk}{n}\leq Ckrn,
    \end{align}
    where in the last inequality we used that $s=(1\pm1/k^2)rn^2/2$ and $k\geq \log n$. This concludes the proof of~\ref{stageII3}.

    It remains to verify~\ref{stageII4}. Recall that $\widetilde{R}=H_q\setminus\bigcup_{K\in\cK}K=H_q\setminus F$. Set
    \begin{align}\label{eq:rtilde-final}
        \widetilde{r}:=pq+r-\frac{2|\cK|}{n^2}=pq+r-\frac{2sk}{n^2}.
    \end{align}
     Note by~\eqref{eq:parametersIII} and the fact that $s=(1\pm 1/k^2)rn^2/2$ it follows that $\widetilde{r}>0$ and that $pq-2rk<\widetilde{r}<pq$. We will see later that~$\widetilde{r}$ corresponds to the density of the remainder $3$-graph $\widetilde{R}$.

    We first estimate the vertex degrees of $\widetilde{R}$. For
    $u\in[n]$, let $a_u$ be the number of matchings in $\cK$ that miss the vertex $u$ defined in~\eqref{eq:a_u_final}. Note that if $u\in V(K_{i,j})$, then either $u\in U_{i,j}$ or $u\in V(E_{i,j})$. In the first case, the edge of $K_{i,j}$ containing $u$ is also an edge of $F_{i,j}\subseteq F$. In the second case, the edge of $K_{i,j}$ containing $u$ is an edge of $E_{i,j}\subseteq R$. Since $|\cK|=sk$, this implies that the degree of $u$ in $F$ is given by
    \begin{align}\label{eq:degree-F-final}
        \deg_F(u)=sk-\deg_R(u)-a_u.
    \end{align}
    Moreover, recall that $H_q$ was chosen such that $\deg_{H_q}(u)=pqn^2/2+2\lambda n^2$. Therefore, by the hypothesis on $R$,~\eqref{eq:a_u_final},~\eqref{eq:rtilde-final} and the fact that $r\leq 1/k^6$, $k\geq \log n$ and $\lambda\geq C\log n/n$, we have that
    \begin{align*}
        \deg_{\widetilde R}(u)&=\deg_{H_q}(u)-\deg_F(u)=\deg_{H_q}(u)-sk+\deg_R(u)+Z_u\\
        &=
        \frac{\widetilde r n^2}{2}\pm3\lambda n^2 \pm Ckrn\leq  \frac{\widetilde r n^2}{2}\pm 4\lambda n^2.
    \end{align*}
    This concludes the proof of the degree estimate in~\ref{stageII4}.

    Finally, we prove the codegree estimate. Fix a pair
    $\{u,v\}\subseteq[n]$. Since our outcome $H_q$ is in the event $\cE_{\textrm{co}}$, we have that $H_q$ is $(\delta/2,pq,n)$-pseudorandom. Therefore, by property~\ref{cF2} and the fact that $\delta\ge C^3/k^2$
    \begin{align*}
        \deg_{\widetilde R}(u,v)=\deg_{H_q}(u,v)-\deg_F(u,v)=\left(1\pm \frac{\delta}{2}\right)pqn-\left(1\pm\frac{\delta}{2}\right)\left(\frac{2sk}{n^2}-\frac{2s}{n^2}\right)n.
    \end{align*}
    The result now follows by noting that $s=(1\pm1/k^2)rn^2/2$ implies that $2s/n^2=(1\pm\delta/10)r$, i.e., $\deg_{\widetilde R}(u)=(1\pm \delta)\widetilde{r}n^2/2$. This concludes the proof of~\ref{stageII4} and the lemma.
\end{proof}

\section{The Iterative Absorption Procedure} \label{sec:iter}

Our goal in this section is to iterate the two stages developed in the previous sections in order to find edge-disjoint maximum matchings covering almost all of $H^{(3)}_{n,p}$. Before we enter into the details of the statement and proof, we should take a step back and remind ourselves of the general picture.

Let $H\sim H_{n,p}^{(3)}$ and $x_{\star}$ be the vertex of maximum degree in $H$. To optimally cover all the edges of $H$ into matchings, we need to ensure that every single matching used in the covering contains the vertex $x_{\star}$. As explained in the overview, our main strategy would be to iterate Stages I and II several times to obtain a sparse remainder graph $R$ such that $x_{\star}$ is still the vertex of maximum degree by  a considerable margin. The proof then could be finished by a procedure introduced in~\cites{hefetz2014optimal, knox2015edge}: One would reserve a set of edges $E$ incident to $x_{\star}$; apply the current bounds on the chromatic index of hypergraphs to cover almost all the edges of $R\setminus E$; and then use a greedy absorption procedure to cover the remaining edges by matchings using edges of $E$.

Two important numerical factors play a role here: The gap between the degrees of $x_{\star}$ and the remaining vertices and how the sparsity of $R$. The next theorem limits the first one of the factors. The proof follows the same idea as the proof for graphs in~\cite[Chapter~3]{bollobas1998random} and, for this reason, we postpone it to Appendix~\ref{appendix:max_deg}.

\begin{theorem}\label{thm:deggap}
   For $p:=p(n)$ satisfying $\log^4 n/n^{2} \leq p \leq 0.99$, a random $3$-graph $H\sim H_{n,p}^{(3)}$ with high probability satisfies the following. If $x_\star \in V(H)$ is the vertex of maximum degree, then
   \begin{align*}
       \deg(x_\star)-\deg(y) \geq \frac{\sqrt{pn^{2}}}{\log^2 n}
   \end{align*}
   for every $y \in V(H) \setminus \{x^\star\}$.
\end{theorem}

The second factor is limited by the nature of the proof strategy implemented in Stages~I and~II. To obtain the results in those stages, we relied on suitable projections of $3$-graphs into pseudorandom bipartite graphs. The pseudorandomness comes with the natural barrier that the codegree of the $3$-graphs used in those stages must be concentrated, and hence of density at least~$1/n$. The main result in this section gives us a $3$-graph $R$ that essentially achieves both numerical barriers.

\begin{lemma}\label{lem:mainiterative}
    Let $\alpha\in(0,1)$ and $p:=p(n)$ be real numbers satisfying $n^{-1+\alpha}\leq p \leq 0.99$. Let $H\sim H_{n,p}^{(3)}$ and $x_{\star}$ be the vertex of maximum degree in $H$. Then, with high probability, the $3$-graph $H$ contains a collection of edge disjoint matchings $\cK$ and a remainder $3$-graph $R=H\setminus \bigcup_{K\in \cK} K$ with the following properties
    \begin{enumerate}[label=(K\arabic*)]
        \item \label{K1} The $3$-graph $R$ satisfies $\Delta(R) \leq n^{1+\alpha}$ and $\Delta_2(R)\leq n^{\alpha}$.
        \item \label{K2} For every $K\in \cK$, it holds that $x_{\star}\in V(K)$.
        \item\label{K3}  For every vertex $y\neq x_{\star}$ it holds that
        \begin{align*}
            \deg_R(x_{\star})-\deg_R(y) \geq \frac{\sqrt{pn^2}}{\log^2 n}.
        \end{align*}
    \end{enumerate}
\end{lemma}

Now we provide an overview of the proof of Lemma~\ref{lem:mainiterative}. The proof will be based on a multiple exposure argument together with several applications of Stages I and II (Lemmas~\ref{lem:global_covering} and~\ref{lem:stageII}). For an integer $t$ that is a constant depending only on $\alpha$, we will expose the $3$-graph $H\sim H_{n,p}^{(3)}$ in $2t$ rounds. In the $i$-th round, we expose the hypergraph $G^{i}$ by including each edge not chosen in an earlier round with probability $q_i$. We choose the probabilities $\{q_i\}_{1\leq i \leq 2t}$ such that
\begin{align*}
    1-p=\prod_{i=1}^{2t}(1-q_i).
\end{align*}
Hence, by a standard argument $H$ can be coupled with $\bigcup_{i=1}^{2t}G^i$. Moreover, for each $1\leq i \leq 2t$, with high probability the $3$-graphs $G^i$ have the properties of binomial random hypergraphs (see Proposition~\ref{prop:coupling}).

We construct the collection of matchings $\cK$ iteratively, by using auxiliary sequences of hypergraphs $\{R^i\}_{1\leq i \leq 2t+1}$ and $\{H^i\}_{1\leq i \leq 2t}$. Depending on the parity of the round, we will either perform Stage I (for the odd rounds) or Stage II (for the even rounds). In the $i$-th round, no matter if $i$ is even or odd, we are going to use as input for our Lemmas~\ref{lem:global_covering} and~\ref{lem:stageII}, the hypergraphs $R^i$ and $H^i$. As an output, in each round we are going to obtain a collection of matchings $\cK_i$ and a remainder $\widetilde{R}^i$ to be set later as $R^{i+1}$. The hypergraphs $H^i$ are going to be chosen as a function of the sprinkled sequence $\{G^i\}_{1\leq i \leq 2t}$ (see more in Subsection~\ref{subsec:iterative}). At the end of the $2t$ rounds, we obtain the collection $\cK=\bigcup_{i=1}^{2t}\cK_i$ and a remainder graph $R:=R^{2t+1}=H\setminus\bigcup_{K\in \cK} K$. 

There are many technical considerations to be discussed in the procedure above that we will postpone for later. We mention two at this moment that might be worthwhile. The first is that in order for the hypothesis of Lemma~\ref{lem:stageII} to be satisfied, the remainder graph used as input needs to be relatively sparse. This and the fact that $R_1$ is taken to be the empty set will imply that a special treatment is required for the very first round of the iteration.

The second point we want to mention is that to ensure that all the matchings $K\in \cK$ contain the vertex $x_{\star}$, we need to choose the exceptional sets $\{S_i\}_{1\leq i\leq 2t}$ in the statements of Stages I and II in such a way that $x_{\star}$ is not contained in any of them. However, without having the full hypergraph $H$ exposed in a current round, that could be a potential challenge. Fortunately, the probabilities $\{q_i\}_{1\leq i\leq 2t}$ are going to be chosen in such a way that $q_1$ is very close to $p$. As a consequence, the first exposure hypergraph $G^1$ with high probability already has $x_{\star}$ as the vertex of maximum degree. This fact will allow us to construct a fixed set $S$ that would work for the entire process (see Proposition~\ref{prop:max_G1}). We now proceed with the details of the proof.

\subsection{Exposing $H^3_{n,p}$ in multiple rounds}\label{subsec:multiple_exposure}

Let $k$ and $t$ be integers chosen such that
\begin{align}\label{eq:parameters_k_t}
k:=\lfloor n^{\alpha/1000} \rfloor \qquad \text{and} \qquad n^{-1+\alpha/3}\leq pk^{-t/4} \leq n^{-1+\alpha/2}.
\end{align}
Moreover, let $\{q_i\}_{1\leq i \leq 2t}$ be real numbers in $(0,1)$ chosen such that
\begin{align}\label{eq:parametersIV}
  1-p=\prod_{i=1}^{2t}(1-q_i), \quad  q_2=k^{-8} q_1  \quad \text{and} \quad  q_i=\begin{cases}
    2q_{i-1}, &\quad \text{if $i\geq 3$ is odd}\\
    k^{-1/4} q_{i-1}. &\quad \text{if $i\geq 3$ is even}
    \end{cases}.
\end{align}
Before we describe the exposure of our $3$-graph, we first show that our choice of probabilities in~\eqref{eq:parametersIV} is well-defined.

\begin{proposition}\label{prop:q_well_defined}
    There exists a choice of $\{q_i\}_{1\leq i\leq 2t}$ satisfying~\eqref{eq:parametersIV}. Moreover, for any such choice it holds that
    \begin{align*}
        (1-4k^{-8})p\leq q_1\leq p, \qquad \text{and} \qquad  n^{-1+\alpha/4}\leq q_{2t}\leq n^{-1+\alpha/2}.
    \end{align*}
\end{proposition}

\begin{proof}
Note that by solving the recursion given in~\eqref{eq:parametersIV}, we obtain that
\begin{align*}
q_i=k^{-8}\cdot k^{-\frac{1}{4}\cdot \lfloor(i-2)/2\rfloor}\cdot 2^{\lfloor (i-1)/2\rfloor}\cdot q_1
\end{align*}
for $2\leq i \leq 2t$. That is, all the probabilities can be expressed as linear functions of $q_1$. Thus, it suffices to show that there is a suitable choice of $q_1$ satisfying~\eqref{eq:parametersIV}.

Let $f_1,\ldots, f_{2t}$ be polynomials defined by
\begin{align*}
    f_1(x):=x \quad \text{and} \quad f_i(x)=k^{-8}\cdot k^{-\frac{1}{4}\cdot\lfloor(i-2)/2\rfloor}\cdot 2^{\lfloor (i-1)/2\rfloor}\cdot x.
\end{align*}
Consider the polynomial $F(x):=\prod_{i=1}^{2t} (1-f_i(x))$. Since $F(0)=1$, $F(1)=0$ and $F$ is continuous on $[0,1]$, we obtain that there exists a value $q_1\in [0,1]$ such that $F(q_1)=1-p$. Moreover, by~\eqref{eq:parameters_k_t}, we have that $k$ is sufficiently large. Consequently, we obtain that $q_i=f_i(q_1)\in (0,1)$ are valid probabilities for $2\leq i \leq 2t$. This proves the existence of a choice satisfying~\eqref{eq:parametersIV}.

For the bounds claimed in the proposition, first note that
\begin{align*}
    1-p=\prod_{i=1}^{2t}(1-q_i)\leq 1-q_1
\end{align*}
immediately implies that $p\geq q_1$. On the other hand, by the fact that $\prod_{i=1}^n(1-x_i)\geq 1-\sum_{i=1}^n x_i$ for $x_i\in [0,1]$ and by~\eqref{eq:parametersIV} we obtain that
\begin{align}\label{eq:p_q_i}
    1-p=\prod_{i=1}^{2t}(1-q_i)\geq 1-\sum_{i=1}^{2t}q_i.
\end{align}
By the upper bound just obtained on $q_1$, the closed formula obtained for $q_i$ and the fact that $k$ is sufficiently large, we have
\begin{align*}
    \sum_{i=2}^{2t}q_i\leq k^{-8}q_i\left(3\sum_{j=0}^{t-1}k^{-j/4}2^{j}\right)\leq 4k^{-8}q_i\leq 4k^{-8}p.
\end{align*}
Putting this together with~\eqref{eq:p_q_i} we obtain that $q_1\geq (1-4k^{-8})p$.

Finally, using the upper bound $q_1\leq p$ together with~\eqref{eq:parameters_k_t}, we obtain that
\begin{align*}
    q_{2t}
    =k^{-8}\cdot k^{-(t-1)/4}\cdot 2^{t-1}\cdot q_1
    \leq k^{-31/4}\cdot 2^{t-1}\cdot pk^{-t/4}\leq pk^{-t/4}\leq n^{-1+\alpha/2},
\end{align*}
where the second to last inequality follows from the fact that $s$ depends only on $\alpha$ and $k$ tends to infinity with $n$. The lower bound follows similarly by considering the inequality $q_1\geq (1-4k^{-8})p$ and the fact that $k= \lfloor n^{\alpha/1000}\rfloor$ and $t\geq 2$ (see~\eqref{eq:parameters_k_t}). Indeed,
\begin{align*}
    q_{2t}=k^{-8}\cdot k^{-(t-1)/4}\cdot 2^{t-1}\cdot q_1\geq 2(1-4k^{-8}) \cdot k^{-31/4}\cdot pk^{-t/4}\geq k^{-31/4}\cdot n^{-1+\alpha/3}\geq n^{-1+\alpha/4} 
\end{align*}
This concludes the proof of the proposition.
\end{proof}

The exposure process consists of $2t$ random rounds. During these rounds, we construct random sequences $\{A^i\}_{1\leq i \leq 2t}$ and $\{G^i\}_{1\leq i\leq 2t}$ as follows. Set $A^1=K_n^{(3)}$ to be the complete $3$-graph on the vertex set $[n]$. We define the remaining elements in $\{A^i\}_{2\leq i \leq 2t}$ and $\{G^i\}_{1\leq i\leq 2t}$ sequentially. For $1\leq i \leq 2t$, let $G^i:=(A^i)_{q_i}$, i.e., $G_i$ is the random subgraph of $A^i$, where each edge was chosen independently with probability $q_i$. For $2\leq i \leq 2t$, let $A^i:=A^{i-1}\setminus G^{i-1}$ be the remaining graphs of edges not chosen so far in our random process. Moreover, let 
\begin{align}\label{eq:b_i}
    b_1=1 \qquad \text{and} \qquad b_i=\prod_{j=1}^{i-1}(1-q_j)
\end{align}
for $2\leq i \leq 2t$.

It is a standard well-known fact that all the $3$-graphs defined in this process are distributed as binomial random hypergraphs, although not necessarily independent. In particular, by our choice of parameters in~\eqref{eq:parametersIV}, we have that
\begin{align}\label{eq:coupling}
    H:=\bigcup_{i=1}^{2t}G^i\sim H_{n,p}^{(3)}.
\end{align}
That is, the union of all the sprinkled graphs has the same distribution as an instance of $H_{n,p}^{(3)}$. From now on, we will work with this model instead of exposing the hypergraph only once. 
Another consequence is the following result that shows that, with high probability, all $\{A^i\}_{1\leq i \leq 2t}$ and $\{G^i\}_{1\leq i\leq 2t}$ are pseudorandom.

\begin{proposition}\label{prop:coupling}
    Let $\{A^i\}_{1\leq i \leq 2t}$ and $\{G^i\}_{1\leq i\leq 2t}$ be the random $3$-graphs defined above using the parameters in~\eqref{eq:parameters_k_t}--\eqref{eq:b_i}. Then, with high probability, the following holds simultaneously
    \begin{enumerate}[label=(L\arabic*)]
        \item\label{couple:A_i} For every $1\leq i \leq 2t$, the $3$-graph $A^i$ is $(k^{-50},b_i,n)$-pseudorandom and
        \begin{align*}
            \deg_{A^i}(u)=\frac{b_i n^2}{2}\pm n\log n.
        \end{align*}
        for every $u\in [n]$.
        \item\label{couple:G_i} For every $1\leq i \leq 2t$, the $3$-graph $G^i$ is $(k^{-50},b_iq_i,n)$-pseudorandom and
        \begin{align*}
            \deg_{G^i}(u)=\frac{b_iq_i n^2}{2}\pm n\log n.
        \end{align*}
        for every $u\in [n]$.
    \end{enumerate}
\end{proposition}

\begin{proof}
Note that by our random procedure, both $A^i$ and $G^i$ can be coupled as random hypergraphs. More specifically, we have
\begin{align*}
    A^i\sim H_{n,b_i}^{(3)} \qquad \text{and} \qquad G^i\sim H_{n,b_iq_i}^{(3)}
\end{align*}
for $1\leq i \leq 2t$. 

We now verify that both sequences $\{A^i\}_{1\leq i \leq 2t}$ and $\{G^i\}_{1\leq i\leq 2t}$ of random $3$-graphs are pseudorandom with high probability. To verify the first sequence, note by our hypothesis on $p$ and that $k\leq n^{\alpha/1000}$ that 
\begin{align}\label{eq:lowerb_i}
    b_i=\prod_{j=1}^{i-1}(1-q_j)\geq \prod_{j=1}^{2t}(1-q_j)=1-p\geq 0.01\geq \frac{Ck^{200}\log n}{n},
\end{align}
for any sufficiently large constant $C$. Let $\cE_A$ be the event that for every $1\leq i \leq 2t$, the $3$-graph $A^i$ is $(k^{-50},b_i,n)$-pseudorandom. Hence, by Corollary~\ref{cor:binomial_hypergraph_pseudo} applied to $\epsilon=k^{-50}$ and a union bound, we obtain that $\PP(\cE_A)\geq 1-n^{-C^{1/3}}$ for sufficiently large $C$.

Similarly, let $\cE_G$ be the event that for every $1\leq i \leq 2t$, the $3$-graph $G^i$ is $(k^{-50},b_iq_i,n)$-pseudorandom. Then by the choices of probabilities in~\eqref{eq:parameters_k_t},~\eqref{eq:parametersIV} and Proposition~\ref{prop:q_well_defined} we can verify that
\begin{align}\label{eq:b_iq_i}
   b_iq_i\geq (1-p)q_{2t}\geq 0.01\cdot n^{-1+\alpha/4}\geq \frac{Ck^{200}\log n}{n}.
\end{align}
Hence, by Corollary~\ref{cor:binomial_hypergraph_pseudo} applied to $\epsilon=k^{-50}$ and a union bound, we obtain that $\PP(\cE_G)\geq 1-n^{-C^{1/3}}$ for sufficiently large $C$.

It remains to verify the degree estimates. However, this follows straightforwardly from the concentration of binomial distributions. Let $J\sim H_{n,q}^{(3)}$ where $q\geq C\log n/n$ for some large constant $C$. Note that by~\eqref{eq:lowerb_i} and~\eqref{eq:b_iq_i}, this assumption applies to both $A^i$ and $G^i$. Then the quantity $\deg_J(u)$ follows a binomial distribution $\Bin\left(\binom{n}{2}, q\right)$ with mean given by $\E(\deg_J(u))=qn^2/2\pm n$. Hence, by Theorem~\ref{chernoff},
\begin{align*}
    \PP\left(
        \left|\deg_J(u)-\E(\deg_J(u))\right|
        >\frac{n\log n}{2}
    \right)
    \leq
    2\exp\left(
        -\frac{n^2\log^2 n}
        {12\,\E(\deg_J(u))}
    \right)
    \leq n^{-C}.
\end{align*}
Therefore, by taking a union bound over all vertices and all the $4t$ random $3$-graphs in the two sequences, with probability at least $1-n^{\sqrt{C}}$ we have simultaneously
\begin{align*}
    \deg_{A^i}(u)=\frac{b_i n^2}{2}\pm n\log n \quad \text{and} \quad 
    \deg_{G^i}(u)
    &=\frac{b_iq_i n^2}{2}\pm n\log n
\end{align*}
for every $u\in[n]$ and every $1\leq i\leq 2t$. Taking the intersection of this event with $\cE_A$ and $\cE_G$ concludes the proof of the proposition.
\end{proof}

We finish this subsection by showing that the first exposure hypergraph $G^1$ already has the same vertex of maximum degree as $H$. Consequently, one can define, at the moment that $G^1$ was exposed, a fixed set $S$ not containing the vertex of maximum degree during all the stages of our iterative process.

\begin{proposition}\label{prop:max_G1}
    Let $\{A^i\}_{1\leq i \leq 2t}$ and $\{G^i\}_{1\leq i\leq 2t}$ be the sequence of $3$-graphs defined above with parameters satisfying~\eqref{eq:parameters_k_t} and~\eqref{eq:parametersIV}. Let $H$ be the $3$-graph defined as in~\eqref{eq:coupling}. Then with high probability the following holds. If $x_{\star}$ is the vertex of maximum degree in $H$ and $S$ is the set of the $n/100$ vertices of $G^1$ of smallest degree, then
    \begin{enumerate}[label=(M\arabic*)]
        \item\label{M:samemax} The vertex $x_\star$ is also the vertex of maximum degree in $G^1$
        \item\label{M:deggap} For every $y\in S$, it holds that
        \begin{align*}
            \deg_H(x_{\star})-\deg_H(y)\geq \frac{\sqrt{pn^2\log n}}{100}
        \end{align*}
        for some positive constant $c>0$.
    \end{enumerate}
\end{proposition}

\begin{proof}
Set $G_L:=\bigcup_{i=2}^{2t}=H\setminus G^1$ to be the $3$-graph exposed after the first round in the remaining $2t-1$ rounds. Our standard coupling shows that $G_L\sim H_{n,q_L}^{(3)}$ where $q_L=p-q_1$. Therefore, by Proposition~\ref{prop:q_well_defined}, we have that
\begin{align*}
    q_L=p-q_1\leq 4k^{-8}p.
\end{align*}
As a consequence, since the quantity $\deg_{G_L}(u)$ follows a binomial distribution $\Bin\left(\binom{n}{2}, q_L\right)$ with mean $\E(\deg_{G_L}(u))=q_L n^2/2\pm n$, we obtain by Theorem~\ref{chernoff} and the previous inequality on $q_L$ that
\begin{align*}
    \PP\left(\left|\deg_{G_L}(u)-\E(\deg_{G_L}(u))\right|>C\sqrt{k^{-8}pn^2\log n}\right)\leq \exp\left(-\frac{C^2\log n}{2}\right)\leq n^{-C}.
\end{align*}
Therefore, by a union bound we obtain with high probability that
\begin{align}\label{eq:deg_GL}
\deg_{G_L}(u)=\frac{q_L n^2}{2}\pm 2C\sqrt{k^{-8}pn^2\log n}.
\end{align}
for every $u\in [n]$. Let $\cE_L$ be the event described by~\eqref{eq:deg_GL}. 

Let $x_1$ be the vertex of maximum degree in $G^1$. By Proposition~\ref{prop:q_well_defined}, we have $q\geq (1-4k^{-8})p$. Since $p\geq n^{-1+\alpha}\geq \log^4 n/n^2$, we have with high probability by Theorem~\ref{thm:deggap} that
\begin{align}\label{eq:deggap_G1}
    \deg_{G^1}(x_1)-\deg_{G^1}(y)\geq \frac{\sqrt{q_1 n^2}}{\log^2 n}.
\end{align}
Moreover, it is a standard fact about degree concentration that the vertex $x_1$ satisfies
\begin{align}\label{eq:maxdeg_G1_bound}
    \deg_{G^1}(x_1)\geq\frac{q_1n^2}{2}+\frac{\sqrt{pn^2\log n}}{20}.
\end{align}
for $p\leq 0.99$. For completeness, and because it is an essential part of the proof of Theorem~\ref{thm:deggap}, we present a proof of~\eqref{eq:maxdeg_G1_bound} in the Appendix~\ref{appendix:max_deg} (see Lemma~\ref{lem:maxdeg1}). Let $\cE_{\textrm{max}}$ be the event such that both~\eqref{eq:deggap_G1} and~\eqref{eq:maxdeg_G1_bound} hold.

Finally, let $Z\subseteq [n]$ be the set of vertices in $G^1$ defined by
\begin{align*}
    Z:=\left\{u\in [n]:\: \deg_{G^1}(u)>\frac{q_1n^2}{2}+\frac{\sqrt{pn^2\log n}}{100}\right\}.
\end{align*}
That is, the vertices of $G^1$ of really high degree. Then, similarly to before, the quantity $\deg_{G^1}(u)$ follows a binomial distribution with probability $q_1$ and mean $\E(\deg_{G^1}(u))=q_1n^2/2\pm n$. Thus, by using the Chernoff bound (Theorem~\ref{chernoff}) and that $q_1<p$, we obtain that
\begin{align*}
\PP\left(\left|\deg_{G^1}(u)-\E(\deg_{G^1}(u))\right|>\frac{\sqrt{pn^2\log n}}{50}\right)\leq 2\exp\left(-\frac{\log n}{300}\right)\leq n^{-1/500}.
\end{align*}
Therefore, the inequality above gives immediately that $\PP(u\in Z)\leq n^{-1/500}$. Consequently, an application of Markov's inequality gives us that
\begin{align*}
    \PP\left(|Z|\geq \frac{99n}{100}\right)\leq \frac{\E|Z|}{99n/100}< 2n^{-1/500}.
\end{align*}
Hence, with high probability the set $S$ of the $n/100$ vertices of smallest degree of $G^1$ satisfies
\begin{align}\label{eq:setS}
    \deg_{G^1}(u)\leq \frac{q_1n^2}{2}+\frac{\sqrt{pn^2\log n}}{50}
\end{align}
for every $u\in S$. Let $\cE_S$ be the event describing~\eqref{eq:setS}.

We claim that conditioned on the event $\cE_L \cap \cE_{\textrm{max}}\cap  \cE_S$, the hypergraphs described in the statement satisfy~\ref{M:samemax} and~\ref{M:deggap}. To prove~\ref{M:samemax}, note by~\eqref{eq:parameters_k_t},~\eqref{eq:deg_GL},~\eqref{eq:deggap_G1} and Proposition~\ref{prop:q_well_defined} that 
\begin{align*}
    \deg_H(x_1)-\deg_H(u)&=(\deg_{G^1}(x_1)-\deg_{G^1}(u))+(\deg_{G_L}(x_1)-\deg_{G_L}(u))\\
    &\geq \frac{\sqrt{q_1 n^2}}{\log^2 n}-4C\sqrt{k^{-8}pn^2\log n}>0
\end{align*}
for every $u\neq x_1$. Hence, the vertex $x_1$ must be exactly $x_{\star}$.

It remains to prove~\ref{M:deggap}. Since $x_1=x_{\star}$ by~\ref{M:samemax}, we obtain by~\eqref{eq:maxdeg_G1_bound} and~\eqref{eq:setS} that
\begin{align*}
    \deg_{G^1}(x_{\star})-\deg_{G^1}(u)>\frac{\sqrt{pn^2\log n}}{50}
\end{align*}
for every $u\in S$. Hence, putting this together with~\eqref{eq:deg_GL}, we obtain that
\begin{align*}
    \deg_{H}(x_{\star})-\deg_H(u)&>\frac{\sqrt{pn^2\log n}}{50}+\deg_{G_L}(x_{\star})-\deg_{G_L}(u)\\&\geq \frac{\sqrt{pn^2\log n}}{50}-4C\sqrt{k^{-8}pn^2\log n}\geq \frac{\sqrt{pn^2\log n}}{100}.
\end{align*}
Since $\cE_L \cap \cE_{\textrm{max}}\cap  \cE_S$ holds whp, this concludes the proof of the proposition.
\end{proof}

\subsection{Formal iterative procedure}\label{subsec:iterative}

In this subsection, we describe formally the iterative procedure to prove Lemma~\ref{lem:mainiterative}. We start by defining all the variables to be used later in the procedure. Let $k$, $t$, $\{q_i\}_{1\leq i\leq 2t}$ be the parameters defined in~\eqref{eq:parameters_k_t} and~\eqref{eq:parametersIV} and let $\{A^i\}_{1\leq i \leq 2t}$, $\{G^i\}_{1\leq i \leq 2t}$ be the sequence of $3$-graphs obtained in the previous subsection. Let $S$ be the set of vertices defined in Proposition~\ref{prop:max_G1}. We define sequences of real numbers $\{\epsilon_i\}_{1\leq i \leq 2t}$, $\{\delta_i\}_{1\leq i \leq 2t}$ and $\{\lambda_i\}_{1\leq i \leq 2t}$ such that
\begin{equation}\label{eq:paramters_eps_del_lam}
\begin{alignedat}{3}
    \epsilon_1 & = k^{-50}, 
    &\qquad \epsilon_i & = k^{-15}, 
    &\qquad & \forall\, 2\leq i\leq 2t,\\
    \delta_1 & = k^{-10}, 
    &\qquad \delta_i & = k^{-3/2}, 
    &\qquad & \forall\, 2\leq i\leq 2t,\\
    \lambda_1 & = \frac{\log^2 n}{n}, 
    &\qquad \lambda_i & = 5^{i-1}\lambda_1, 
    &\qquad & \forall\, 2\leq i\leq 2t.
\end{alignedat}
\end{equation}
Finally, we set the sequence of hypergraphs $\{H^i\}_{1\leq i \leq 2t}$ and sets $\{S_i\}_{1\leq i \leq 2t}$ such that
\begin{align}\label{eq:sequence_Hi}
    H^i=\begin{cases}
        G^i, &\quad \text{if $i$ is odd},\\
        A^i, &\quad \text{if $i$ is even},
    \end{cases}
    \quad \text{and} \quad S_i=S,
\end{align}
for every $1\leq i \leq 2t$.

The procedure that we are going to describe next is a randomized one and it has a probability of failure. This is a consequence of the fact that the pseudorandomness of the sequence $\{H^i\}_{1\leq i \leq 2t}$ and Stage II (see Lemma~\ref{lem:stageII}) are not deterministic statements. The iteration, as described in the overview, also produces a sequence of remainder $3$-graphs $\{R^i\}_{1\leq i \leq 2t+1}$, where we initialize by setting $R^1=\emptyset$. For a step $1\leq j \leq 2t$, suppose that we already obtained the remainder graph $R^j$ and that the procedure has not declared failure. If $j$ is odd, then we are applying Stage I (Lemma~\ref{lem:global_covering}) with the hypergraphs $H^j$ and $R^j$ and the set $S_j$ as inputs and parameters $\epsilon_j$, $\delta_j$ and $\lambda_j$. If there are no parameters $p_j$ and $r_j$ satisfying the hypothesis of the lemma, then we declare failure. Otherwise, we obtain a remainder graph $R^{j+1}:=\widetilde{R}^j$ and a collection of matchings $\cK_j$ satisfying properties~\ref{global1}--\ref{global4}.

In the case that $j$ is even, then we apply Stage II (Lemma~\ref{lem:stageII}) with $H^j$ and $R^j$ and the set $S_j$ as inputs and parameters $\epsilon_j, \delta_j, \lambda_j, q_j$ and $k$ as defined previously. If there are no parameters $p_j$ and $r_j$ satisfying the hypothesis of the lemma or if the lemma does not succeed, then we declare failure. Otherwise, we obtain a remainder graph $R^{j+1}:=\widetilde{R}^j$ and a collection of matchings $\cK_j$ satisfying properties~\ref{stageII1}--\ref{stageII4}. An important feature of this step is that the random subgraph in the statement is $H^i_{q_i}=(A^i)_{q_i}=G^i$, i.e., the randomness coming from the statement in Stage II is encoded in the construction of $G^i$. Finally, if the process did not declare failure after the $2t$ steps, then we return $R:=R^{2t+1}$ and $\cK=\bigcup_{i=1}^{2t}\cK_i$. The following is an algorithmic description of the procedure.

\begin{algorithm}[H]
\caption{Iterative absorption procedure}
\label{alg:iterative}
\begin{algorithmic}[1]
\STATE Set $R^1:=\emptyset$ and $\cK=\emptyset$.
\FOR{$j=1$ \textbf{to} $2t$}
        \IF{$j$ is odd}
            \STATE Apply Lemma~\ref{lem:global_covering} to $H^j$, $R^j$ and $S_j$
            \IF{$H^j$, $R^j$ and $S_j$ fail hypothesis}
                \STATE Declare \textbf{failure}
            \ENDIF
            \STATE Let $R^{j+1}:=\widetilde{R}^j$ and $\cK:=\cK\cup \cK_j$.
        \ENDIF
        \IF{$j$ is even}
            \STATE Apply Lemma~\ref{lem:stageII} to $H^j$, $R^j$ and $S_j$
            \IF{$H^j$, $R^j$ and $S_j$ fail hypothesis or Lemma~\ref{lem:stageII} fails}
                \STATE Declare \textbf{failure}
            \ENDIF
            \STATE Let $R^{j+1}:=\widetilde{R}^j$ and $\cK:=\cK\cup \cK_j$.
        \ENDIF.
\ENDFOR
\STATE Return $R^{2t+1}$ and $\cK$.
\end{algorithmic}
\end{algorithm}

We now show that with high probability the procedure above succeeds. 

\begin{proposition}\label{prop:iteration_succeeds}
   The iterative procedure in Algorithm~\ref{alg:iterative} does not declare failure with high probability.
\end{proposition}

\begin{proof}
Let $\{b_i\}_{1\leq i \leq 2t}$ be the sequence defined in~\eqref{eq:b_i}. Set $\{p_i\}_{1\leq i \leq 2t}$ to be the sequence of real numbers defined by
\begin{align}\label{eq:parameter_p_i}
    p_i=\begin{cases}
        b_i, &\quad \text{if $i$ is even},\\
        b_iq_i, &\quad \text{if $i$ is odd}.
    \end{cases}
\end{align}
A consequence of Proposition~\ref{prop:coupling}, our choice of parameters in~\eqref{eq:paramters_eps_del_lam} and~\eqref{eq:sequence_Hi} is that with high probability $H^i$ is $(\epsilon_i,p_i,n)$-pseudorandom for every $1\leq i \leq 2t$ and
\begin{align*}
    \deg_{H^i}(u)=\frac{p_in^2}{2}\pm n\log n=\frac{p_in^2}{2}\pm \lambda_i n^2,
\end{align*} 
for every $u\in [n]$. Let $\cE_H$ be this event.

Let $\cE_2$ be the event that our random procedure fails because one of the applications of Lemma~\ref{lem:stageII} failed. That is, even though $H^i$, $R^i$ and $S_i$ satisfied the hypothesis of the lemma, the lemma itself returned a failure. A union bound over all possible applications of Lemma~\ref{lem:stageII} gives that $\PP(\cE_2)\leq t \cdot n^{-\sqrt{C}}$. Since $C$ is sufficiently large and $t$ is a constant depending only on $\alpha$ defined by~\eqref{eq:parameters_k_t}, we obtain that $\cE_2^\mathsf{c}$ holds with high probability. As already mentioned before, the source of randomness playing a role in the event $\cE_2$ comes exclusively from our random multiple round exposure of the sequences $\{A^i\}_{1\leq i \leq 2t}$ and $\{G^i\}_{1\leq i \leq 2t}$.

We claim that Algorithm~\ref{alg:iterative} succeeds for every outcome in $\cE_H\cap \cE_2^{\mathsf{c}}$. We prove it by induction on the number of steps. For each step $1\leq i \leq 2t$, we suppose by induction hypothesis that there exists a $3$-graph $R^i$ satisfying the hypothesis of the step with parameter $r_i$. Our goal is to show that the remainder $3$-graph $R^{i+1}$, obtained at the end of the step, satisfies the hypothesis of the next iteration. As already discussed, we initialize the entire induction by taking $R_1=\emptyset$ and $r_1=0$. For each step we are going to apply Lemmas~\ref{lem:global_covering} and~\ref{lem:stageII} with sufficiently large constants $\{C_i\}_{1\leq i \leq 2t}$ satisfying
\begin{align*}
    C_1\ll C_2\ll\cdots \ll C_{2t}.
\end{align*}
Note that since $2t$ is a constant depending on $\alpha$, the constants $C_i$ might also depend on $\alpha$. Since we are on the event $\cE_H$ and $S_i=S$ for every $1\leq i \leq 2t$, there is no need to check that $H^i$ and $S_i$ satisfy the hypothesis of the lemmas.

We start the induction by doing the case $i=1$ separately. In this case, note that $R^1=\emptyset$. Thus, we may use the last assertion of Lemma~\ref{lem:global_covering}. Indeed, one can verify that the remaining numerical hypothesis in the statement holds from~\eqref{eq:parameters_k_t} and~\eqref{eq:paramters_eps_del_lam}. Thus, Lemma~\ref{lem:global_covering} outputs a $3$-graph $R^{2}:=\widetilde{R}^1$ such that
\begin{align*}
    \deg_{R^{2}}(u)=c^{(1)}_{\star}\delta_1p_1 n^2\pm 5\lambda_1n^2 \quad \text{and} \quad \Delta_2(R^2)\leq C_1^3\delta_1p_1n,
\end{align*}
for every $u\in [n]$, where $c^{(1)}_\star$ is a constant bounded both above and below by a function of $C_1$. Therefore, for sufficiently large $C_2$, $r_2=2c^{(1)}_{\star}\delta_1p_1$ and our definition of $\lambda_2$ it holds that
\begin{align*}
\deg_{R^{2}}(u)=\frac{r_2n^2}{2}\pm \lambda_2 n^2 \quad \text{and} \quad \Delta_2(R^2)\leq C_2r_2 n.
\end{align*}
This already satisfies the hypothesis assumption for the application of Lemma~\ref{lem:stageII} to $R^2$. It remains to check if the parameters $\epsilon_2,\delta_2, p_2. r_2. q_2. \lambda_2$ satisfy the numerical hypothesis of the statement in~\eqref{eq:parametersIII}. Indeed, by~\eqref{eq:parameters_k_t},~\eqref{eq:b_i},~\eqref{eq:paramters_eps_del_lam} and~\eqref{eq:parameter_p_i}, we immediately obtain that $\log n\leq k\leq n^{0.01}$, $\epsilon_2\leq k^{-2}$, $\delta_{2}\geq C_2^3/k^2$ and $p_2=b_2\geq 1-p\geq 0.01$. For the condition concerning $r_2$, note that $r_2=2c^{(1)}_{\star}\delta_1p_1=2c^{(1)}_{\star}k^{-10}q_1<k^{-6}$. Moreover, by Proposition~\ref{prop:q_well_defined} we have $q_1\geq (1-4k^{-8})p$. Together with our choice of $k$ and the fact that $p\geq n^{-1+\alpha}$, we have
\begin{align}\label{eq:even_r_1}
    r_2=2c^{(1)}_\star k^{-10} q_1\geq pk^{-11}\geq n^{-1+\alpha/4}\geq \frac{k^7\log n}{n}.
\end{align}
For the condition concerning $q_2$, note that by~\eqref{eq:parametersIV}
\begin{align*}
    q_2=k^{-8}q_1=\frac{k^2r_2}{2c^{(1)}_{\star}}\geq C_2^4 k r_2\log^3 n,
\end{align*}
since $k$ is polynomial in $n$. Finally, $\lambda_2=5\log^2 n/n$ and therefore $C_2\log n/n<\lambda_2<r_2/3k^2$. This concludes the first case of the induction.

We now split the argument into two cases depending on the parity of the step $i$. Note by the way that we are inducting, we only need to verify until step $2t-1$.

\vspace{0.2cm}

\noindent \underline{Case 1:} $3\leq i\leq 2t-1$ is odd.

\vspace{0.2cm}

The computation here is very similar to the case $i=1$. Suppose that we have a $3$-graph $R^{i}$ satisfying the hypothesis of Lemma~\ref{lem:global_covering} and all the other parameters are checked. Thus, an application of the lemma outputs a $3$-graph $R^{i+1}:=\widetilde{R}^i$ such that
\begin{align*}
    \deg_{R^{i+1}}(u)=c^{(i)}_{\star}\delta_ip_i n^2\pm 5\lambda_i \quad \text{and} \quad C_i^3\delta_ip_in,
\end{align*}
for every $u\in [n]$, where $c^{(i)}_{\star}$ is a constant bounded both above and below by a function of $C_i$. Therefore, for sufficiently large $C_{i+1}$, $r_{i+1}=2c^{(i)}_{\star}\delta_ip_i$ and our definition of $\lambda_{i+1}$ it holds that
\begin{align*}
\deg_{R^{i+1}}(u)=\frac{r_{i+1}n^2}{2}\pm \lambda_{i+1} n^2 \quad \text{and} \quad \Delta_2(R^{i+1})\leq C_{i+1}r_{i+1} n.
\end{align*}
This already satisfies the hypothesis assumption of Lemma~\ref{lem:stageII} regarding $R^{i+1}$. Checking if the parameters $\epsilon_{i+1},\delta_{i+1}, p_{i+1}, r_{i+1}, q_{i+1}$ and $\lambda_{i+1}$ satisfy~\eqref{eq:parametersIII} is mostly similar to the previous case. The main difference is in our treatment of the inequality concerning $r_{i+1}$, $q_{i+1}$ and $\lambda_{i+1}$. For the inequalities related to $r_{i+1}$, note that by~\eqref{eq:parametersIV},~\eqref{eq:paramters_eps_del_lam} and~\eqref{eq:parameter_p_i}
\begin{align*}
    r_{i+1}=2c^{(i)}_{\star} \delta_ip_i=2c^{(i)}_{\star}k^{-3/2}b_iq_i\leq 4c^{(i)}_{\star}k^{-19/2}q_1<k^{-6},
\end{align*}
for $i\geq 3$ odd. For the lower bound, we have by Proposition~\ref{prop:q_well_defined},~\eqref{eq:parameters_k_t} and the fact that $b_i>1-p\geq 0.01$ that
\begin{align}\label{eq:even_r_2}
    r_{i+1}=2c^{(i)}_{\star}k^{-3/2}b_iq_i\geq k^{-2}q_{2t}\geq k^{-2}n^{-1+\alpha/4}\geq \frac{k^7\log n}{n}.
\end{align}
To check the inequality concerning $q_{i+1}$, note that by~\eqref{eq:parametersIV}
\begin{align*}
    q_{i+1}=k^{-1/4}q_i>C_{i+1}^4\cdot c^{(i)}_{\star}k^{-1/2}q_i\log^3 n=C_{i+1}^4kr_{i+1}\log^3 n.
\end{align*}
Finally, for $\lambda_{i+1}$, just note that $\lambda_{i+1}>\lambda_1=\log^2 n/n>C_{i+1}\log n/n$ and that
\begin{align*}
    \lambda_{i+1}=5^i\lambda_1<\frac{5^{2t}\log^2 n}{n}<\frac{k^5\log n}{3n}\leq \frac{r_{i+1}}{3k^2},
\end{align*}
since $t$ is a constant depending only on $\alpha$ and $k$ is a polynomial in $n$. This concludes this step of the induction.

\vspace{0.2cm}

\noindent \underline{Case 2:} $2\leq i\leq 2t-2$ is even.

\vspace{0.2cm}

Suppose that we have a $3$-graph $R^i$ and all the other parameters satisfying the hypothesis of Lemma~\ref{lem:stageII}. Since we are in an outcome in $\cE_2^{\mathsf{c}}$, our application of Lemma~\ref{lem:stageII} succeeded and produced a $3$-graph $R^{i+1}:=\widetilde{R}^{i}$ such that $R^{i+1}\subseteq G^i$ and 
\begin{align*}
    \deg_{R^{i+1}}(u)=\frac{r_{i+1}n^2}{2}\pm 4\lambda_i n^2 \quad \text{and} \quad \deg_{R^{i+1}}(u,v)=(1\pm \delta_i)r_{i+1}n<C_{i+1}r_{i+1}n,
\end{align*}
for all $u,v \in [n]$, where $r_{i+1}=p_iq_i+r_i-2|\cK_i|/n^2$ and $C_{i+1}$ is a sufficiently large constant. By our definition of $\lambda_{i+1}$ in~\eqref{eq:paramters_eps_del_lam}, we have that $R^{i+1}$ satisfies the hypothesis of Lemma~\ref{lem:global_covering}. It remains to check if $\epsilon_{i+1}, \delta_{i+1}, r_{i+1}, p_{i+1}$ and $\lambda_{i+1}$ satisfy the numerical hypothesis of the statement. The inequality concerning the sizes of $\epsilon_{i+1}$, $\delta_{i+1}$ and $\lambda_{i+1}$ follows immediately from~\eqref{eq:paramters_eps_del_lam}. Moreover, the last assertion on~\ref{stageII4} gives us that $p_iq_i-2r_ik<r_{i+1}<p_{i}q_{i}$. For the remaining inequality, by the definition of the parameters~\eqref{eq:b_i} and by the fact that $q_i<2k^{-8}$ for $i\geq 2$, we have
\begin{align}\label{eq:b_iq_i-to_i+1}
    b_{i+1}q_{i+1}=b_i(1-q_{i})q_{i+1}=2b_i(1-q_i)q_i\geq b_iq_i.
\end{align}
Hence, the upper bound follows from~\eqref{eq:b_iq_i-to_i+1} and the fact that $i$ is even
\begin{align*}
    r_{i+1}<p_iq_i=b_iq_i\leq b_{i+1}q_{i+1}=p_{i+1}.
\end{align*}
To see the lower bound, note by induction that $q_i\geq C_i^4 kr_i\log^3 n$, $p_i=b_i\geq 1-p\geq 0.01$ and that $r_i\geq k^{-2}n^{-1+\alpha/4}$ by our bounds for even $r_i$ obtained in~\eqref{eq:even_r_1} and~\eqref{eq:even_r_2}. Therefore,
\begin{align*}
    r_{i+1}> p_iq_i-2r_ik\geq 0.01\cdot (C_i^4 \log^3 n -2) r_i k\geq k^{-2}n^{-1+\alpha/4}\log ^3 n\geq \frac{C_{i+1}\log n}{\epsilon_{i+1}^2n},
\end{align*}
since $\epsilon_{i+1}^2=k^{-30}>n^{-\alpha/30}$. This concludes the proof of the induction and of the claim.

\vspace{0.2cm}

The proof of the proposition now follows since $\cE_{H}\cap \cE_2^{\mathsf{c}}$ holds with high probability.
\end{proof}

\subsection{Proof of Lemma~\ref{lem:mainiterative}}

We are now ready to prove the main lemma of this section.

\begin{proof}[Proof of Lemma~\ref{lem:mainiterative}]
    Let $k$, $t$ and $\{q_i\}_{1\leq i \leq 2t}$ be the parameters defined in~\eqref{eq:parameters_k_t} and~\eqref{eq:parametersIV}. Let $\{A^i\}_{1\leq i \leq 2t}$ and $\{G^i\}_{1\leq i \leq 2t}$ be the sequence of $3$-graphs defined in Subsection~\ref{subsec:multiple_exposure}. By our exposure procedure, it holds that $H\sim H_{n,p}^{(3)}$ has the same distribution as the union of $G^i$ as in~\eqref{eq:coupling}. Let $H$ be one instance of $H_{n,p}^{(3)}$ satisfying Proposition~\ref{prop:iteration_succeeds}, i.e., an instance where the algorithm described in Subsection~\ref{subsec:iterative} succeeded. We are also going to require that $H$ satisfies Theorem~\ref{thm:deggap} and Proposition~\ref{prop:max_G1}. Since all these statements hold with high probability, such an instance $H$ exists with high probability. Let $R:=R^{2t+1}$ be the last remainder hypergraph obtained in Algorithm~\ref{alg:iterative} and $\cK=\bigcup_{i=1}^{2t} \cK_i$. By construction, it holds that $R=H\setminus\bigcup_{K\in \cK}K$. We claim that $H$, $R$, $\cK$ satisfy conditions~\ref{K1}--~\ref{K3} of Lemma~\ref{lem:mainiterative}.

    To check~\ref{K1}, note that $R:=R^{2t+1}$ is the $3$-graph obtained by the last application of Stage II (Lemma~\ref{lem:stageII}) and is such that
    \begin{align*}
        \deg_{R}(u)=\frac{r_{2t+1}n^2}{2}\pm 4\lambda_{2t}n^2 \quad \text{and} \quad \deg_R(u,v)=(1\pm \delta_{2t})r_{2t+1}n
    \end{align*}
    for $u,v\in [n]$, where the parameters $\delta_{2t}$, $\lambda_{2t}$ and $r_{2t+1}$ were introduced in the last subsection. Therefore, by the last sentence of~\ref{stageII4} and Proposition~\ref{prop:q_well_defined} we have
    \begin{align*}
        r_{2t+1}< p_{2t}q_{2t}\leq q_{2t}\leq n^{-1+\alpha/2}.
    \end{align*}
    Moreover, by~\eqref{eq:paramters_eps_del_lam}, we have $\lambda_{2t}\leq 5^{2t}\log^2 n/n< n^{-1+\alpha/2}$. Putting it all together, we obtain that
    \begin{align*}
        \Delta_1(R)\leq 5n^{1+\alpha/2}<n^{1+\alpha} \quad \text{and} \quad \Delta_2(R)\leq 2n^{\alpha/2}<n^{\alpha}.
    \end{align*}
    This concludes the proof of~\ref{K1}.

    The proof of~\ref{K2} comes from the fact that both conditions~\ref{global2} and~\ref{stageII2} imply that for both Stages I and II it holds that $K\in \cK_i$ satisfies $[n]\setminus V(K)\subseteq S_i$. However, by our construction and Proposition~\ref{prop:max_G1}, the set $S_i$ is equal to the set $S$ of the $n/100$ vertices of smallest degree in $G^1$. In particular, by~\ref{M:samemax}, the vertex of largest degree $x_{\star}$ in $H$ is not contained in $S$. Hence, we have that $x_{\star}\in V(K)$ for every $K\in \cK$.

    Finally, we verify~\ref{K3}. We split into two cases. If $y\notin S$, then by the same argument as in the previous paragraph, we have that $y\in V(K)$ for every $K\in \cK$. Since $R=H\setminus \bigcup_{K\in \cK}K$, this implies by Theorem~\ref{thm:deggap} that
    \begin{align*}
        \deg_R(x_{\star})-\deg_R(y)=\deg_H(x_{\star})-\deg_H(y)\geq \frac{\sqrt{pn^2}}{\log^2 n}.
    \end{align*}

    Otherwise, if $y\in S$, then by~\ref{global3} for every odd $i$, the vertex $y\in S$ is missing in $C_ip_in$ matchings from $\cK_i$. In the same spirit, from~\ref{stageII3} for every even $i$, the vertex $y\in S$ is missing in $C_ikr_in$ matchings from $\cK_i$. Therefore, the total number of matchings that $y$ misses in $\cK$ is at most
    \begin{align*}
        \sum_{i=1}^tC_{2i-1}p_{2i-1}n+\sum_{i=1}^t C_{2i}kr_{2i}n\leq 2C_{2t}tpn,
    \end{align*}
    where we use that $p_{2i-1}=b_{2i-1}q_{2i-1}\leq p$ for $1\leq i \leq t$ and $r_{2i}=2c^{(2i-1)}_{\star}\delta_{2i-1}p_{2i-1}<k^{-1}p$ by the induction hypothesis in the proof of Proposition~\ref{prop:iteration_succeeds} and our choice of $\delta_{2i-1}$. Therefore, by property~\ref{M:deggap} and the fact that $t$ is a constant not depending on $n$, it holds that
    \begin{align*}
        \deg_R(x_{\star})-\deg_R(y)\geq\deg_H(x_{\star})-\deg_H(y)-2C_{2t}tpn\geq \frac{\sqrt{pn^2\log n}}{100}-2C_{2t}pn\geq \frac{\sqrt{pn^2}}{\log^2 n}
    \end{align*}
    for every $y\in S$. This concludes the proof of the lemma.
\end{proof}

\section{Greedy recoloring and proof of Theorem~\ref{thm:main}}\label{sec:recoloring}

In this section we conclude the proof of Theorem~\ref{thm:main} by implementing the strategy already highlighted in the overview. The only missing part in the proof is the following greedy recoloring that, given a sufficiently large reservoir of edges incident to a vertex $x$, allows us to rearrange matchings not containing $x$ to matchings containing it.

\begin{proposition}\label{prop:recoloring}
    Let $G$ be a $3$-graph on $n$ vertices such that the following holds. There exists a vertex $x\in [n]$ and a partition $G=E\cup F_1\cup\ldots \cup F_t$ such that
    \begin{enumerate}[label=(\roman*)]
        \item\label{recoloring1} Every edge $e \in E$ contains the vertex $x$,
        \item\label{recoloring2} $F_i$ is a matching and $x\notin V(F_i)$ for $1\leq i \leq t$,
        \item\label{recoloring3} $|E|>3(t+\Delta_2(G))$,
    \end{enumerate}
    Then there exists a partition of $G=F_1'\cup\ldots \cup F_s'$ into matchings such that $x\in V(F_i')$ for every $1\leq i \leq s$.
\end{proposition}

\begin{proof}
The proof of the proposition is done by greedily rearranging the matchings to contain the vertex $x$. We claim that for any fixed matching $F_i$, there exist three edges $E_i=\{e^{(i)}_1,e_2^{(i)},e_3^{(i)}\} \subseteq E$ such that $F_i\cup E_i=F_{i,1}\cup F_{i,2}\cup F_{i,3}$ where each $F_{i,j}$ is a matching. We do it inductively on the number of matchings. Suppose that for $1\leq j\leq t$, we already rearranged the first $j-1$ matchings and used the disjoint set of edges $E_1,\ldots,E_{j-1} \subseteq E$ to do that. 

For the next matching $F_j$, we are going to choose three edges from the set $E\setminus (E_1\cup\ldots \cup E_{j-1})$. Since each $E_i$ has size three, we obtain that $|E\setminus (E_1\cup\ldots \cup E_{j-1})|\geq |E|-3(j-1)\geq 3\Delta_2(G)+3$. Choose an arbitrary edge of this set to be $e_{1}^{(j)}$. Since $x \notin V(F_j)$, the edge $e_{1}^{(j)}$ is disjoint from at least $|F_j|-2$ edges from $F_j$. Let $\widetilde{F}_{j}\subseteq F_j$ be this set of edges. We may assume that there are exactly two edges $f_{2}^{(j)}$ and $f_3^{(j)}$ intersecting $e_1^{(i)}$ (the case with one or no edges follows similarly). Since $\deg_G(x,z)\leq \Delta_2(G)$ for every $z\in f$, we obtain that at most $3\Delta_2(G)$ edges of $E$ intersect with $f$. Therefore, because $|E\setminus (E_1\cup\ldots\cup E_{j-1})|\geq 3\Delta_2(G)+3$, we can greedily choose two edges $e_2^{(j)}$ and $e_3^{(j)}$ that are disjoint from  $f_{2}^{(j)}$ and $f_3^{(j)}$. Hence, the matchings $F_{i,1}=\widetilde{F}_j\cup \{e_1^{(j)}\}$, $F_{i,2}=\{e_2^{(j)},f_2^{(j)}\}$ and $F_{i,3}=\{e_3^{(j)},f_3^{(j)}\}$ all contain the vertex $x$. This finishes the proof of the claim.

The proof of the proposition now follows by relabelling these matchings and adding the remaining edges in $E$ as matchings consisting of a single edge.
\end{proof}

We now use Proposition~\ref{prop:recoloring} to finish the proof of the theorem.

\begin{proof}[Proof of Theorem~\ref{thm:main}]
    Set $\alpha=1/200$. Let $H\sim H_{n,p}^{(3)}$ be an instance of the random $3$-graph satisfying the statements of Theorem~\ref{thm:deggap} and Lemma~\ref{lem:mainiterative}. Since both statements hold with high probability, then such an instance exists with high probability. Let $x_{\star}$ be the vertex of maximum degree of $H$, let $\cK$ be the collection of matchings obtained in Lemma~\ref{lem:mainiterative} and $R$ the remainder $3$-graph. Our goal is to use Proposition~\ref{prop:recoloring} to obtain a partition of $R$ into a collection of matchings all of them containing the vertex $x_{\star}$.

    Let $E\subseteq R$ be a set of edges incident to $x_{\star}$ of size $\frac{\sqrt{pn^2}}{2\log^2n}$ and let $\widetilde{R}=R\setminus E$. As a consequence of condition~\ref{K3}, we have that
    \begin{align*}
        \deg_{\widetilde{R}}(x_{\star})-\deg_{\widetilde{R}}(y)\geq\deg_{R}(x_{\star})-\deg_R(y)-\frac{\sqrt{pn^2}}{2\log^2n}\geq \frac{\sqrt{pn^2}}{2\log^2 n}
    \end{align*}
    for every $y\neq x_{\star}$. In particular, $x_{\star}$ is still the vertex of maximum degree in $\widetilde{R}$. Moreover, by condition~\ref{K1} the $3$-graph $\widetilde{R}$ satisfies
    \begin{align*}
        \Delta_1(\widetilde{R})\leq\Delta_1(R)\leq n^{1+\alpha} \quad \text{and} \quad \Delta_2(\widetilde{R})\leq \Delta_2(R)\leq n^{\alpha}.
    \end{align*}
    Therefore, an application of Theorem~\ref{thm:kangnibble} with $\eta=1/100<1/96$ gives us that
    \begin{align*}
        \chi'(\widetilde{R})\leq \Delta_1(\widetilde{R})+n^{2/3+\alpha-1/100}=\Delta_1(\widetilde{R})+n^{2/3-1/200}
    \end{align*}
    That is, there is a covering of $\widetilde{R}$ using at most $\Delta_1(\widetilde{R})+n^{2/3-1/200}$ matchings. Let $\cF_{\star}$ be the collection of edge disjoint matchings containing the vertex of maximum degree $x_\star$ and let $F_1\cup \ldots \cup F_t$ be the matchings not containing $x_{\star}$. One can check by the bound above that $t\leq n^{2/3-1/200}$. Finally, let $G=E\cup F_1\cup \ldots \cup F_t$.

    We claim that $G$ and $x_{\star}$ satisfy the hypothesis of Proposition~\ref{prop:recoloring}. Hypotheses~\ref{recoloring1} and~\ref{recoloring2} follow immediately from the construction of the set $E$ and the matchings $F_1,\ldots, F_t$. To check~\ref{recoloring3} note that by~\ref{K1} we have that
    \begin{align*}
        \Delta_2(G)\leq \Delta_2(R)\leq n^{\alpha}=n^{1/200}.
    \end{align*}
    Therefore, since $p\geq n^{-2/3-1/200}$ we obtain
    \begin{align*}
        |E|\geq \frac{\sqrt{pn^2}}{2\log^2 n}\geq \frac{n^{2/3-1/400}}{2\log^2 n}>3n^{2/3-1/200}+3n^{1/200}\geq 3(t+\Delta_2(G)).
    \end{align*}
    Hence, we can apply Proposition~\ref{prop:recoloring} to obtain a partition of $G$ into a collection of edge disjoint matchings $\cF$. Take $\cM=\cK\cup \cF_{\star}\cup \cF$. By Lemma~\ref{lem:mainiterative}, the collection of matchings $\cK$ covers all the edges of $H\setminus R$. The last paragraphs just showed that $R\subseteq \cF\cup \cF_{\star}$. Moreover, by~\ref{K2} and Proposition~\ref{prop:recoloring}, every matching in $\cK\cup \cF\cup \cF_{\star}$ contains the vertex $x_{\star}$ and they are all disjoint. Thus, $|\cM|=\Delta_1(H)$, which concludes the proof.
\end{proof}

\section{Concluding Remarks}

In this paper, we showed that with high probability the random $3$-uniform hypergraph $H_{n,p}^{(k)}$ is of Class 1 for a broad range of $p$. This raises the natural question if the same result can be extended to higher uniformities in an appropriate setting. As shown in
Proposition~\ref{prop:no_vizing_4}, the theorem does not hold when $k$ does not divide $n$. More precisely, when $k\nmid n$ and
$p\geq C\log n/n^{k-3}$, the random hypergraph
$H\sim H_{n,p}^{(k)}$ satisfies $\chi'(H)>\Delta(H)$ with high
probability. Is the above divisibility remark the only obstruction for extending Theorem~\ref{thm:main}? This motivates the following question.

\begin{question}
Let $k\geq 4$ be fixed and suppose that $k\mid n$. For which ranges of $p=p(n)$ does the random $k$-uniform hypergraph $H\sim H_{n,p}^{(k)}$ satisfy $\chi'(H)=\Delta(H)$ with high probability?
\end{question}

We remark that for $p=1$ the above question reduces to the classical Baranyai's theorem on the decomposition of a complete hypergraph into perfect matchings.

\subsection*{Statement on AI use}

The work leading to the proof of Theorem~\ref{thm:main} was done by the authors without the use of any AI tools. However, the authors used AI to assist with editing and proofreading the arguments in the manuscript. 

\bibliographystyle{abbrv}
\bibliography{References}

\appendix

\section{Proof of Propositions \ref{prop:random_subbipartite_pseudo} and \ref{prop:random_hypergraph_pseudo}}\label{appendix:pseudo}

In this appendix we restate and prove Propositions~\ref{prop:random_subbipartite_pseudo} and~\ref{prop:random_hypergraph_pseudo}.

\bipartitepseudo*

\begin{proof}
    We first check \ref{item:random_subbipartite_deg}. Note that for $z\in V(B)$, we can write $\deg_{B_{\phi}}(z)=\sum_{\{z,w\}\in B}\1_{\{z,w\}\in B_\phi}$, where $\1_{\{z,w\}\in B_{\phi}}$ are independent $0/1$-random variables with mean $\phi(z,w)$. Thus for $\gamma\leq 1$, we have
    \begin{align*}
        \E(\deg_{B_{\phi}}(z))=\sum_{\{z,w\}\in B}\phi(z,w)=(1\pm \eta)q\deg_{B}(z)=(\gamma \pm 3\delta/4)pqn.
    \end{align*}
    Hence, by the Chernoff bound (Theorem \ref{chernoff}) and our hypothesis on $p$, $q$ and $\epsilon$, it holds that
    \begin{align*}
        \PP\left(|\deg_{B_{\phi}}(z)-\E(\deg_{B_{\phi}}(z))|>\frac{\delta pqn}{4}\right)\leq 2\exp\left(-\frac{\delta^2 pqn}{48(1+3\delta/4)^2}\right)\leq n^{-C/2}.
    \end{align*}
    A union bound over all the vertices of $B$ shows that the probability that \ref{item:random_subbipartite_deg} does not hold is at most $n\cdot n^{-C/2}\leq n^{-C/4}< n^{-\sqrt{C}}$.

    To prove \ref{item:random_subbipartite_pseudo}, note that we already verified in \ref{item:random_subbipartite_deg} (by taking $\gamma=1$) that the probability that $B_{\phi}$ does not satisfy \ref{P1:degree} is at most $n^{-C/4}$. It remains to check that $B_{\phi}$ satisfies properties \ref{P2:largesets}, \ref{P3:smallS} and \ref{P4:smallT}. Note that for $S\subseteq X$, $T\subseteq Y$, we can write $e_{B_{\phi}}(S,T)=\sum_{s\in S, t\in T}\1_{\{s,t\}\in B_{\phi}}$, where again $\1_{\{s,t\}\in B_{\phi}}$ are independent $0/1$-random variables with mean $\phi(s,t)$. Therefore, we have $\E(e_{B_\phi}(S,T))=(1\pm\eta)q\cdot e_B(S,T)$. Since $|S|, |T|\geq \delta n>\epsilon n$ and $B$ is $(\epsilon,p,n)$-pseudorandom, we have by property \ref{P2:largesets} that 
    \begin{align*}
        \E(e_{B_\phi}(S,T))=(1\pm\eta)(1\pm \epsilon)pq|S||T|=(1\pm3\delta/4)pq|S||T|
    \end{align*}
    Thus, the Chernoff bound implies
    {\footnotesize
    \begin{align}\label{eq:chernoff_ii}
        \PP\left(|e_{B_{\phi}}(S,T)-\E(e_{B_{\phi}}(S,T))|>\frac{\delta pq|S||T|}{4}\right)\leq 2\exp\left(-\frac{\delta^2 pq|S||T|}{48(1+3\delta/4)^2}\right)\leq\exp\left(-\frac{\delta^4pqn^2}{200}\right)\leq 2^{-Cn}.
    \end{align}
    }
    Since there are at most $2^{2n}$ choices of $S$ and $T$, by a union bound we obtain that the probability that \ref{P2:largesets} fails is at most $2^{2n}\cdot2^{-Cn}\leq n^{-C/4}$.

    To verify \ref{P3:smallS}, we split the argument into two cases depending on the size of $|S|$

    \vspace{0.2cm}

    \noindent \underline{Case 1:} $\epsilon n\leq |S|<\delta n$ and $|T|=2|S|$

    \vspace{0.2cm}

    Since $|T|\geq 2\epsilon n>\epsilon n$ and $B$ is $(\epsilon,p,n)$-pseudorandom, then by property \ref{P2:largesets} and the same calculation as in (\ref{eq:chernoff_ii}), we obtain that with probability $1-n^{-C/4}$
    \begin{align*}
        e_{B_{\phi}}(S,T)\leq (1\pm \delta)pq|S||T|\leq 4\delta pq n|S|
    \end{align*}
    for every $S, T$ with $\epsilon n\leq |S|<\delta n$ and $|T|=2|S|$.

     \vspace{0.2cm}

    \noindent \underline{Case 2:} $|S|<\epsilon n$ and $|T|=2|S|$

    \vspace{0.2cm}

    In this case, note that by property \ref{P3:smallS} of $B$ we have that $e_B(S,T)\leq 4\epsilon pn |S|$. Hence,
    \begin{align*}
        \PP\left(e_{B_{\phi}}(S,T)> 4\delta pq n|S|\right)&\leq \sum_{\substack{F\subseteq E_B(S,T)\\ |F|=4\delta pq n|S|}}\prod_{f\in F}\phi(f)\leq \binom{e_B(S,T)}{4\delta pq n|S|}\left((1+\eta)q\right)^{4\delta pq n|S|}\\
        &\leq\left(\frac{e (1+\eta)q \cdot e_B(S,T)}{4\delta pq n|S|}\right)^{4\delta pq n|S|}\leq  \exp\left(-C|S|\log n\right)
    \end{align*}
    There are at most $\binom{n}{|S|}\cdot\binom{n}{|T|}\leq \left(\frac{en}{2|S|}\right)^{4|S|}$ choices of $S$ and $T$ with the size of $S$ fixed and $|T|=2|S|$. Thus, by a union bound we obtain that the probability that \ref{P3:smallS} fails is at most
    \begin{align*}
        \sum_{|S|=1}^{\epsilon n} \left(\frac{en}{2|S|}\right)^{4|S|}\cdot \exp(-C|S|\log n)\leq \sum_{|S|=1}^{\epsilon n}n^{-C|S|/2}\leq n^{-C/4}
    \end{align*}

    The computations for property \ref{P4:smallT} hold similarly. Therefore, we conclude that with probability $1-4n^{-C/4}\geq 1-n^{-\sqrt{C}}$ the graph $B_{\phi}$ satisfies properties \ref{P1:degree}--\ref{P4:smallT} and consequently $B_\phi$ is $(\delta,pq,n)$-pseudorandom.
    \end{proof}

\hypergraphpseudo*

\begin{proof}
        The proof is very similar to the one in Proposition \ref{prop:random_subbipartite_pseudo}. We verify that with high probability $H_q$ satisfies properties \ref{Q1:codegree}--\ref{Q4:smallZ} of Definition \ref{def:hypergraph_pseudorandom}. For property \ref{Q1:codegree}, note that since $H$ is $(\epsilon,p,n)$-pseudorandom, we have $\deg_H(x,y)=(1\pm \epsilon)pn$ for every $x,y \in V$. Moreover, we can write $\deg_{H_q}(x,y)=\sum_{\{x,y,z\}\in H}\1_{\{x,y,z\}\in H_q}$ as the sum of independent $0/1$-random variables with mean
        \begin{align*}
            \E(\deg_{H_q}(x,y))=q\deg_H(x,y)=(1\pm \epsilon)pqn=(1\pm \delta/2)pqn.
        \end{align*}
        Hence, by Chernoff (Theorem \ref{chernoff})
        \begin{align*}
            \PP\left(|\deg_{H_q}(x,y)-\E(\deg_{H_q}(x,y))|>\frac{\delta pqn}{2}\right)\leq 2\exp\left(-\frac{\delta^2pqn}{12(1+\delta/2)^2}\right)\leq n^{-C/2}.
        \end{align*}
        Thus, after a union bound, the probability that $H_q$ fails \ref{Q1:codegree} is at most $\binom{n}{2}\cdot n^{-C/2}\leq n^{-C/4}$.

        To check property \ref{Q2:largesets}, note that $|M|$, $|Z|\geq \delta n>\epsilon n$ implies that $e_H(M,Z)=(1\pm \epsilon)p|M||Z|$. Moreover, we can write $e_{H_q}(M,Z)=\sum_{\{x,y,z\} \in E_H(M,Z)}\1_{\{x,y,z\}\in H_q}$, where $\1_{\{x,y,z\}\in H_q}$ are independent $0/1$-random variables and
        \begin{align*}
            \E(e_{H_q}(M,Z))=q\cdot e_H(M,Z)=(1\pm \epsilon)pq|M||Z|=(1\pm \delta/2)pq|M||Z|.
        \end{align*}
        Thus, by Chernoff,
        {\footnotesize
        \begin{align}\label{eq:chernoffiii}
            \PP\big(|e_{H_q}(M,Z)-\E(e_{H_q}(M,Z))|>\frac{\delta pq|M||Z|}{2}\big)\leq 2\exp\left(-\frac{\delta^2pq|M||Z|}{12(1+\delta/2)^2}\right)\leq \exp\left(-\frac{\delta^4 pqn^2}{50}\right)\leq n^{-Cn/2}.
        \end{align}
        }
        Since every matching corresponds to a permutation of $[n]$, there are at most $n!$ ways of choosing a matching $M\subseteq V^{(2)}$. Hence, there are at most $n!\cdot 2^n\leq n^{2n}$ ways to choose the pair $(M,Z)$. Consequently, the probability that $H_q$ fails \ref{Q2:largesets} is at most $n^{2n}\cdot n^{-Cn/2}<n^{-C/4}$.

        To verify \ref{Q3:smallM}, we split the calculation into two cases depending on the size of $M$

        \vspace{0.2cm}

        \noindent \underline{Case 1}: $\epsilon n\leq |M|<\delta n$ and $|Z|=2|M|$.

        \vspace{0.2cm}

        Since $|Z|\geq 2\epsilon n>\epsilon n$ and $H$ is $(\epsilon,p,n)$-pseudorandom, we obtain that $e_H(M,Z)=(1\pm \epsilon)p|M||Z|$. By the same calculation as in (\ref{eq:chernoffiii}) we have that with probability at least $1-n^{-C/4}$
        \begin{align*}
            e_{H_q}(M,Z)\leq (1+\delta)pq|M||Z|\leq 4\delta pq n|M|.
        \end{align*}

        \vspace{0.2cm}

        \noindent \underline{Case 2:} $|M|<\epsilon n$ and $|Z|=2|M|$.

        \vspace{0.2cm}

        By property \ref{Q3:smallM} on $H$, we have that $e_H(M,Z)\leq 2\epsilon pn|M|$. Thus,
        \begin{align*}
            \PP\left(e_{H_q}(M,Z)> 4\delta pq n|M|\right)&\leq \sum_{\substack{F\subseteq E_H(M,Z)\\|F|=4\delta pqn|M|}}q^{|F|}=q^{4\delta pq n|M|}\binom{e_H(M,Z)}{4\delta pq n|M|}\\&
            \leq \left(\frac{eq\cdot e_H(M,Z)}{4\delta pq n|M|}\right)^{2\delta pqn|M|}\leq \left(\frac{e}{4}\right)^{4\delta pq n|M|}\leq n^{-C|M|}.
        \end{align*}
        There are at most $n^{2|M|}\cdot \binom{n}{|Z|}\leq n^{4|M|}$ ways to choose a pair $(M,Z)$ with $M$ of fixed size and $|Z|=2|M|$. Hence, the probability that $H_q$ fails property \ref{Q3:smallM} is at most
        \begin{align*}
            \sum_{|M|=1}^{\epsilon n} n^{4|M|}\cdot n^{-C|M|}\leq \sum_{|M|=1}^{\epsilon n} n^{-C|M|/2}\leq n^{-C/4}.
        \end{align*}

        An analogous computation shows that the probability that property \ref{Q4:smallZ} does not hold for $H_q$ is also at most $n^{-C/4}$. Thus, we obtain that $H_q$ is $(\delta,pq,n)$-pseudorandom with probability at least $1-4n^{-C/4}\geq 1-n^{-\sqrt{C}}$.
    \end{proof}

\section{Proof of Lemma \ref{lem:kangoreslemma}}\label{appendix:kang}

In this appendix we show that Definition~\ref{def:bipartite_pseudorandom} implies the definition of pseudorandomness proposed in~\cite{kang2023thresholds}. As a consequence, we deduce Lemma~\ref{lem:kangoreslemma} from their original version. To start with, we introduce their exact notion of pseudorandomness in~\cite{kang2023thresholds}.

\begin{definition}\label{def:quasirandomkang}
    Let $\epsilon, p \in [0,1]$ and let $B$ be a bipartite graph with parts $X, Y$ with $|X|=|Y|=n$. We say that $B$ is \emph{$(\delta,p)$-quasirandom} if the following holds for some $\delta' \in [0, \delta]$:
    \begin{enumerate}[label=(QR\arabic*)]
        \item\label{QR1:degree} For every $z\in V(B)$, we have $\deg_B(z)=(1\pm\delta)pn$;
        \item\label{QR2:largesets} For every $S\subseteq X$, $T\subseteq Y$ such that $|S|$, $|T|\geq \delta' n$, it holds
        \begin{align*}
            e_B(S,T) \geq (1- \epsilon)p|S||T|;
        \end{align*}
        \item\label{QR3:smallS} For every $S\subseteq X$, $T\subseteq Y$ such that $|S|<\delta' n$ and $|T|<2|S|$, it holds
        \begin{align*}
            e_B(S,T)\leq pn|S|/3;
        \end{align*}
        \item\label{QR4:smallT} For every $S\subseteq X$, $T\subseteq Y$ such that $|T|<\delta' n$ and $|S|<2|T|$, it holds
        \begin{align*}
            e_B(S,T)\leq  pn |T|/3.
        \end{align*}
    \end{enumerate}
\end{definition}

We can now state the lemma precisely as it is in the original paper:

\begin{lemma}[Lemma 4.4 in \cite{kang2023thresholds}]
Let $1/n_0 \ll \delta \ll1$, let $n \geq n_0$, let $p \in [0,1]$, and let $\gamma \geq 0$. Let $H \subset K_{n,n}$ be $(\delta,p)$-quasirandom. If $L$ is a $(\gamma\pm \delta)pn$-regular spanning subgraph of $K_{n,n}-E(H)$, then there exists $R \subset H \cup L$ such that $R \supseteq L$ is spanning and regular.
\end{lemma} 


We now show that Definition~\ref{def:bipartite_pseudorandom}
implies the quasirandomness conditions in
Definition~\ref{def:quasirandomkang}.

\begin{proposition}\label{prop:pseudo_implies_kang}
    Let $0<\epsilon\leq 1/12$ and $p\in[0,1]$.
    If $B$ is $(\epsilon,p,n)$-pseudorandom, then
    $B$ is $(\epsilon,p)$-quasirandom.
\end{proposition}

\begin{proof}
    We prove that $B$ is $(\epsilon, p)$-quasirandom by checking conditions~\ref{QR1:degree}--~\ref{QR4:smallT}. Note that conditions~\ref{QR1:degree} and~\ref{QR2:largesets}
    follow immediately from conditions~\ref{P1:degree} and~\ref{P2:largesets}, respectively. Indeed, condition~\ref{P2:largesets} is a two-sided version of~\ref{QR2:largesets}. Next we verify condition~\ref{QR3:smallS} for $\delta'=\epsilon$. Let $S\subseteq X$ and $T\subseteq Y$ be subsets satisfying $|S|<\epsilon n$ and $|T|<2|S|$. Since $2|S|<2\epsilon n\leq n$, we may choose
    a superset $T\subseteq T'\subseteq Y$ such that $|T'|=2|S|$. Thus, by condition~\ref{P3:smallS} and the fact that $\epsilon\leq 1/12$, we have
    \begin{align*}
        e_B(S,T)\leq e_B(S,T')\leq 4\epsilon pn|S|\leq \frac{pn|S|}{3},
    \end{align*}
    This proves condition~\ref{QR3:smallS}. The proof of condition~\ref{QR4:smallT} is analogous.
\end{proof}

\section{Properties of the maximum degree of $H_{n,p}^{(3)}$}\label{appendix:max_deg}

In this appendix, we give a proof of Theorem \ref{thm:deggap}. For completeness, we present the argument for the binomial random $k$-uniform hypergraph $H_{n,p}^{(k)}$ with $k\geq 3$. For brevity, we introduce some notation that will be used throughout the section. For integers $0 \leq d \leq N$ and a binomial random variable $\Bin(N,p)$, let
\begin{align}\label{eq:bshortcut}
    b(N,d) = \Prob\left(\Bin(N,p) = d\right) \quad \text{and} \quad B(N,d) = \Prob\left(\Bin(N,p) \geq d\right) = \sum_{i=d}^N b(N,i).
\end{align}
We use the following sharp standard estimates of the binomial distribution (see, e.g., \cite{randomgraphs}).

\begin{theorem}[\cite{randomgraphs}, Theorems 1.2, 1.5, and 1.6]\label{thm:binom}
Suppose $pN \geq 1$ and $1 \leq t \leq (1-p)N/3$.
\begin{enumerate}
    \item[$(i)$] If $d = pN + t$, then
    \begin{align*}
        b(N,d) < \frac{1}{\sqrt{2\pi p(1-p)N}}\exp\left(-\frac{t^2}{2p(1-p)N}+\frac{t}{(1-p)N}+\frac{t^3}{p^2N^2}\right). 
    \end{align*}
    Moreover, it also holds that
    \begin{align*}
        b(N,d) > \frac{1}{\sqrt{2\pi p(1-p)N}}\exp\left(-\frac{t^2}{2p(1-p)N}-\frac{t^3}{2(1-p)^2N^2}-\frac{t^4}{3p^3N^3}-\frac{t}{2pN}-\beta\right),
    \end{align*}
    where $\beta = \frac{1}{12d} + \frac{1}{12(N-d)}$.
    \item[$(ii)$] If $d = pN + t$ and $t=x\sqrt{p(1-p)N}=o\left(\left(p(1-p)N\right)^{2/3}\right)$, then
    \begin{align*}
        B(N,d) = (1+o(1))\frac{1}{x\sqrt{2\pi}}\exp\left(-\frac{x^2}{2}\right)
    \end{align*}
    for $x \to \infty$.
\end{enumerate}
\end{theorem}

As a consequence of Theorem \ref{thm:binom}, one can obtain relatively sharp bounds on the binomial probabilities when $t$ is close to the standard deviation.

\begin{corollary}\label{cor:binom}
Suppose $pN/\log^3N \to \infty$ as $N \to \infty$, $p\leq 0.99$, and $t = \sqrt{2\alpha p(1-p)N\log N}$ for $\alpha:=\alpha(N)>0$ such that $\lim_{N\rightarrow \infty}\alpha(N)=C$ for some positive constant $C$. Then the following holds for $d = pN + t$:
\begin{align}\label{eq:binom2}
   b(N,d) &= \frac{1}{\sqrt{2\pi p(1-p)N}} \exp\left(-\alpha\log N + O(1)\right),\\
   B(N,d) &=(1+o(1))\frac{1}{\sqrt{4\pi \alpha \log N}} \exp\left(-\alpha\log N \right). \nonumber
\end{align}
\end{corollary}

We remark that Corollary \ref{cor:binom} is stronger than the Chernoff bounds (Theorem \ref{chernoff}) in the sense that it provides a sharper upper bound and a lower bound. Unfortunately, it comes with the restriction that it is only applicable for small deviations and certain regimes of $p$. The next result provides a lower bound on the vertex of maximum degree.

\begin{lemma}\label{lem:maxdeg1}
Let $N = \binom{n-1}{k-1}$ and let $p := p(n)$ satisfy $pn^{k-1}/\log^3 n \to \infty$ as $n \to \infty$ and $p\leq 0.99$. Let $Y$ be the set defined by
\begin{align*}
    Y := \left\{v \in V(H_{n,p}^{(k)}): \deg(v) \geq pN + \sqrt{\frac{2}{k-1}\left(1 - \frac{\log\log n}{\log n}\right)p(1-p)N\log N}\right\}.    
\end{align*}
Then, whp, the set $Y$ is non-empty. In particular, whp,
\begin{align*}
    \Delta(H_{n,p}^{(k)}) \geq pN + \sqrt{\frac{2}{k-1}\left(1 - \frac{\log\log n}{\log n}\right)p(1-p)N\log N}.
\end{align*}
\end{lemma}

\begin{proof}
The proof follows from an application of the second moment method. For convenience, let
\begin{align*}
    t:=\sqrt{\frac{2}{k-1}\left(1 - \frac{\log\log n}{\log n}\right)p(1-p)N\log N}\quad \text{and} \quad \alpha:=\frac{1}{k-1}\left(1-\frac{\log \log n}{\log n}\right).
\end{align*}
For a vertex $i \in V(H_{n,p}^{(k)})$, let $Y_i$ be the indicator random variable of whether $i$ is an element of the set $Y$. A computation together with Corollary \ref{cor:binom} shows that
\begin{align}\label{eq:Ey}
    \E(Y)&=nB(N,pN+t)= (1+o(1))\frac{n}{\sqrt{4\pi\alpha\log N}} \exp\left(-\frac{1}{k-1}\left(1-\frac{\log \log n}{ \log n}\right)\log N\right)\nonumber \\
    &= (1+o(1))\frac{n}{\sqrt{4\pi\log n}} \exp\left(\log \log n -\log n +O(1)\right)>c\sqrt{\log n}.
\end{align}
To compute $\Var(Y)$, we will need the following claim.

\begin{claim}\label{clm:jointYij}
Let $i,j \in V(H^{(k)}_{n,p})$ be two distinct vertices. Then,
\begin{align*}
    \Prob(Y_i=Y_j=1)\leq(1+o(1))\frac{\E(Y)^2}{n^2}.
\end{align*}
\end{claim}

\begin{proof}
Define
\begin{align}\label{eq:d*}
    d^* := \max\left\{100\log n,\, p\binom{n-2}{k-2} + 3\sqrt{p\binom{n-2}{k-2}\log n}\right\}.
\end{align}
Let $N(\{i,j\})$ be the neighborhood of the pair $\{i,j\}$, and let $E = \{|N(\{i,j\})| < d^*\}$ denote the event that the pair degree of $\{i,j\}$ is smaller than $d^*$. We claim that $\Prob(E^{\mathsf{c}}) = o(n^{-2})$. Indeed, if $p\binom{n-2}{k-2} \leq 9\log n$, then by the fact that $|N(\{i,j\})|$ follows the distribution $\Bin\left(\binom{n-2}{k-2}, p\right)$, we have
\begin{align*}
    \Prob(E^{\mathsf{c}}) \leq \Prob(|N(\{i,j\})| > 100\log n) < \binom{\binom{n-2}{k-2}}{100\log n}p^{100\log n}\leq \left(\frac{ep\binom{n-2}{k-2}}{100\log n}\right)^{100\log n}= o(n^{-2}).
\end{align*}
Otherwise, if $p\binom{n-2}{k-2} > 9\log n$, applying Theorem \ref{chernoff} with $\epsilon = 3\sqrt{\frac{\log n}{p\binom{n-2}{k-2}}} < 1$, we obtain
\begin{align*}
    \Prob(E^{\mathsf{c}}) \leq \Prob\left(|N(\{i,j\})| > p\binom{n-2}{k-2} + 3\sqrt{p\binom{n-2}{k-2}\log n}\right) < \exp(-3\log n) = o(n^{-2}).
\end{align*}

Using the estimate above, we derive
\begin{align}\label{eq:Yij1}
    \Prob(Y_i = Y_j = 1) &= \Prob(Y_i = Y_j = 1 \mid E)\Prob(E) + \Prob(Y_i = Y_j = 1 \mid E^{\mathsf{c}})\Prob(E^{\mathsf{c}})\nonumber\\
    &\leq \Prob(Y_i = Y_j = 1 \mid E) + o(n^{-2}).
\end{align}
To estimate $\Prob(Y_i = Y_j = 1 \mid E)$, we introduce some notation. For any pair $\{i,j\}\in [n]^{(2)}$, there are $\binom{n-2}{k-1}$ potential edges containing $i$, but avoiding $j$. Likewise, there are $\binom{n-2}{k-2}$ potential edges containing both $i$ and $j$. Denote by $N'=\binom{n-2}{k-1}$ and $N''=\binom{n-2}{k-2}$. Standard binomial identities give us that $N=N'+N''$. Hence, by only looking at the potential edges that are exclusive to~$i$ or $j$, we have
\begin{align}\label{eq:Yij2}
    \Prob(Y_i = Y_j = 1 \mid E) &\leq B\left(N', pN + t - d^*\right)^2.
\end{align}

To finalize the computation, note that
\begin{align}\label{eq:d*andt}
    d^* - p\binom{n-2}{k-2} = o\left(\frac{t}{\log n}\right).
\end{align}
Indeed, if $d^* = 100\log n$, then $d^* - p\binom{n-2}{k-2} \leq 100\log n = o(t/\log n)$ for $pn^{k-1}/\log^3 n\rightarrow \infty$ as $n\rightarrow \infty$. Otherwise, we have by~\eqref{eq:d*} that
\begin{align*}
    d^*-pN''=d^* - p\binom{n-2}{k-2} = 3\sqrt{p\binom{n-2}{k-2}\log n} = O\left(\frac{t}{\sqrt{n}}\right)=o\left(\frac{t}{\log n}\right).
\end{align*}
Therefore, by Corollary \ref{cor:binom}, along with (\ref{eq:Ey}), (\ref{eq:Yij2}), (\ref{eq:d*andt}), the facts that $N=N'+N''$ and $N'=\binom{n-2}{k-1} = (1 + O(1/n))N$, we obtain
\begin{align*}
    \Prob(Y_i = Y_j = 1 \mid E) &\leq B\left(N', pN' + t + o(t/\log n)\right)^2\\
    &=(1+o(1))\frac{1}{4\pi \alpha \log N}\exp\left(-\frac{2}{k-1}\left(1 - \frac{\log\log n}{\log n}\right)\log N + o(1)\right)\\
    &=(1+o(1))\frac{\E(Y)^2}{n^2}
\end{align*}
The result now follows by (\ref{eq:Yij1}) and the fact that $\E(Y)>c\sqrt{\log n}$.
\end{proof}
Applying Claim \ref{clm:jointYij}, we obtain that
\begin{align*}
    \Var(Y)= \E(Y)+\sum_{i\neq j}\Prob(Y_i=Y_j=1)-\E(Y)^2\leq \E(Y)+o(\E(Y)^2)
\end{align*}
Therefore, by Chebyshev's inequality and $\E(Y)>c\sqrt{\log n}$ we have
\begin{align*}
    \Prob(Y=0)\leq \frac{\Var(Y)}{\E(Y)^2}=\frac{1}{\E(Y)}+o(1)=o(1).
\end{align*}
This concludes that whp the set $Y$ is non-empty.
\end{proof}

We are now ready to state and prove the main result of this section, which in particular implies Theorem \ref{thm:deggap}.

\begin{theorem}\label{thm:deggap_k}
   For $p:=p(n)$ satisfying $\log^4 n/n^{k-1} \leq p \leq 0.99$, a random $k$-graph $H:=H_{n,p}^{(k)}$ whp satisfies the following. If $x^\star \in V(H)$ is the vertex of maximum degree, then
   \begin{align*}
       \deg(x^\star)-\deg(y) \geq \frac{\sqrt{pn^{k-1}}}{\log^2 n}
   \end{align*}
   for every $y \in V(H) \setminus \{x^\star\}$.
\end{theorem}

\begin{proof}
The proof follows from a similar calculation as in the last lemma. Set
\begin{align}\label{eq:Ns}
    N = \binom{n-1}{k-1}, \quad N' = \binom{n-2}{k-1}, \quad \text{and} \quad N'' = \binom{n-2}{k-2}.
\end{align}
Note that $N = N' + N''$. Moreover, let
\begin{align}\label{eq:tstandard}
    t:= \sqrt{\frac{2}{k-1}\left(1 - \frac{\log\log n}{\log n}\right)p(1-p)N\log N} \qquad \text{and} \qquad 
    L:= \frac{\sqrt{pn^{k-1}}}{\log^2 n}.
\end{align}
Let $Y$ be the set defined in Lemma \ref{lem:maxdeg1}, and let
\begin{align*}
    A_{x,y} := \{\deg(x) \geq pN + t \text{ and } \deg(x) - L <\deg(y) \leq \deg(x)\}
\end{align*}
be the event that the vertices $x$ and $y$ are such that $x \in Y$ and the vertex $y$ has a gap smaller than $L$ to the degree of $x$. Set
\begin{align*}
    \widetilde{d}(x) =\left|\left\{e\in H:\: x\in e,\, y\notin e\right\}\right| \quad \text{and} \quad \widetilde{d}(y) =\left|\left\{e\in H:\: x\notin e,\, y\in e\right\}\right|.
\end{align*} 
Note that these two variables are independent and have distribution $\Bin(N',p)$. 

To compute the probability that the event $A_{x,y}$ happens, let $E_{x,y}=\{\deg(x,y)<d^*\}$ be the event that the pair degree of $\{x,y\}$ is smaller than $d^{*}$ as defined in~\eqref{eq:d*}. As noted before in Claim~\ref{clm:jointYij}, we have that $\PP(E_{x,y}^{\mathsf{c}})=o(n^{-2})$. Therefore,
\begin{align*}
    \PP(A_{x,y})\leq \PP(A_{x,y}\cap E_{x,y})+o(n^{-2}).
\end{align*}
It remains to estimate $\Prob(A_{x,y}\cap E_{x,y})$. For that, note that by (\ref{eq:d*andt}) and (\ref{eq:tstandard})
\begin{align}\label{eq:ugly_estimate}
    pN + t - d^* - L = pN' + t + O(t/\log n)>pN'.
\end{align}
This in particular implies that $b(N',\ell)$ is decreasing for $\ell\geq pN+t-d^*-L$. Thus, a computation shows that
\begin{align*}
    \Prob(A_{x,y}\cap E_{x,y}) &= \sum_{i=0}^{d^*-1} \Prob(\deg(\{x,y\}) = i) \left( \sum_{j=pN+t-i}^{N'} \Prob(\widetilde{d}(x) = j) \left( \sum_{\ell=j-L+1}^j \Prob(\widetilde{d}(y) = \ell) \right) \right)\\
    &= \sum_{i=0}^{d^*-1} b(N'',i) \left( \sum_{j=pN+t-i}^{N'} b(N',j) \left( \sum_{\ell=j-L+1}^j b(N',\ell) \right) \right)\\
    &\leq \sum_{i=0}^{d^*-1} b(N'',i) \left( \sum_{j=pN+t-i}^{N'} L \cdot b(N',j-L)^2 \right)\\
    &\leq L \sum_{i=0}^{d^*-1} b(N'',i) b(N',pN+t-i-L) B(N',pN+t-i-L)\\
    &\leq L\cdot b(N',pN'+t+O(t/\log n)) B(N',pN'+t+O(t/\log n)).
\end{align*}
Therefore, by Corollary \ref{cor:binom}, we obtain for our choice of $p$ that
\begin{align*}
    \Prob(A_{x,y}) &\leq L \cdot b(N',pN'+t+O(t/\log n)) B(N',pN'+t+O(t/\log n)) + o(n^{-2})\\
    &=O\left(\frac{L}{\pi\sqrt{8p(1-p)N\log n}}\exp\left(-\frac{2}{k-1}\left(1-\frac{\log\log n}{\log n}\right)\log N\right)\right)+o(n^{-2})\\
    &=O\left(\frac{L\log^{3/2}n}{n^2\sqrt{pn^{k-1}}}\right)=o(n^{-2}).
\end{align*}
Hence, the probability that there exist two vertices $\{x,y\}$ satisfying $A_{x,y}$ is 
\begin{align*}
    \Prob\left(\bigcup_{x,y}A_{x,y}\right) = O\left(\frac{L\log^{3/2}n}{\sqrt{pn^{k-1}}}\right) = o(1).
\end{align*}
Therefore, whp no vertex in $Y$ has a degree gap less than $\sqrt{pn^{k-1}}/\log^2 n$ to vertices of smaller degree. In particular, by Lemma \ref{lem:maxdeg1}, this is true for the vertex of maximum degree $x^{\star}\in Y$.
\end{proof}

\end{document}